\documentclass[10pt]{article}
\usepackage{latexsym,amsmath,amssymb,graphics,amscd}
\usepackage{multirow}
\usepackage{booktabs}
\usepackage{color}
\usepackage{tabularx}
\usepackage[hang,footnotesize]{caption}
\usepackage[pdftex]{graphicx}
\usepackage{graphicx}
\addtocounter{footnote}{1}
\makeatletter
\@addtoreset{table}{bsection}
\def\thetable{\thesection.\@arabic\c@table}
\def\fps@table{h, t}
\@addtoreset{equation}{section}

\makeatother

\newtheorem{theorem}{Theorem}[section]
\newtheorem{definition}[theorem]{Definition}
\newtheorem{lemma}[theorem]{Lemma}
\newtheorem{remark}[theorem]{Remark}
\newtheorem{proposition}[theorem]{Proposition}

\newcommand{\bfi}{\bfseries\itshape}

\newsavebox{\savepar}

\makeatletter

\begin{document}

\title{\textbf{Stability of Hamiltonian relative equilibria in symmetric magnetically confined rigid bodies}}
\author{Lyudmila Grigoryeva$^{1}$, Juan-Pablo Ortega$^{2}$, and Stanislav Zub$^{3}$}
\date{}
\maketitle

\begin{abstract}
This work studies the symmetries, the associated momentum map, and  relative equilibria of a mechanical system consisting of a small axisymmetric magnetic body-dipole in an  also axisymmetric external magnetic field that additionally exhibits a mirror symmetry; we call this system the ``orbitron". We study the nonlinear stability of a branch of equatorial quasiorbital relative equilibria using the energy-momentum method and we provide sufficient conditions for their $\mathbb{T}^2$--stability that complete partial stability relations already existing in the literature. These stability prescriptions are explicitly written down in terms of the some of the field parameters, which can be used in the design of stable solutions. We propose new linear methods to determine instability regions in the context of relative equilibria that we use to conclude the sharpness of some of the nonlinear stability conditions obtained. 
\end{abstract}

\bigskip


\makeatletter
\addtocounter{footnote}{1} \footnotetext{%
Laboratoire de Math\'{e}matiques de Besan\c{c}on, Universit\'{e} de Franche-Comt\'{e}, UFR des
Sciences et Techniques. 16, route de Gray. F-25030 Besan\c{c}on cedex. France; Taras Shevchenko National University of Kyiv. 64, Volodymyrs'ka.
01601 Kyiv. Ukraine. {\texttt{Lyudmyla.Grygoryeva@univ-fcomte.fr} }}
\makeatother
\makeatletter
\addtocounter{footnote}{1} \footnotetext{%
Corresponding author. Centre National de la Recherche Scientifique, Laboratoire de Math\'{e}matiques de Besan\c{c}on, Universit\'{e} de Franche-Comt\'{e}, UFR des
Sciences et Techniques. 16, route de Gray. F-25030 Besan\c{c}on cedex.
France. {\texttt{Juan-Pablo.Ortega@univ-fcomte.fr} }}
\makeatother
\makeatletter
\addtocounter{footnote}{1} \footnotetext{%
Taras Shevchenko National University of Kyiv. 64, Volodymyrs'ka.
01601 Kyiv. Ukraine. {\texttt{stah\_z@yahoo.com} }}
\makeatother

\section{Introduction}

Many physical systems exhibit symmetries. A number of techniques have been developed during the last two centuries to take advantage of the conservation laws that are usually associated to these invariance properties to simplify or {\it reduce} those systems in order to make easier the computation of their solutions. The presence of symmetries also creates natural dynamical features that generalize distinguished solutions of their non-symmetric counterparts like the so called {\it relative equilibria} or {\it relative periodic orbits}; relative equilibria are solutions of a symmetric system that coincide with one-parameter group orbits of the action that leaves that system invariant. The justification of this denomination lies in the fact that relative equilibria are equilibria for the reduced Hamiltonian system~\cite{mwr} constructed  with the momentum map associated to the action, provided that this object exists. Regarding the stability of these solutions, the degeneracies caused by the presence of symmetries in a system cause drift phenomena that make non-evident the selection of a stability definition. A very reasonable choice is the concept of stability relative to a subgroup introduced in~\cite{smodg, thesis:patrick} for which a number of energy-momentum based sufficient conditions have been formulated in the literature under different assumptions and levels of generality on the group actions involved and the momentum values at which the relative equilibrium in question takes place~\cite{smodg, drift, thesis:patrick, mrs, mre, persistence:periodic, molecules, Ortega1999a, singreleq, po, pos, thesis, Roberts2002, patrick:roberts:wulff:2002, Ortega2004, MontaldiOlmos2011}. These methods have been used, for example, in the study of the stability of relative equilibria present in different configurations of  rigid bodies~\cite{lrsm, llt, thesis:patrick}, Riemann ellipsoids~\cite{Fasso2001, Rodriguez-Olmos2008}, underwater vehicles~\cite{leonard:1997, leonard:marsden:1997}, vortices~\cite{Pekarsky1998, laurent:nonlinearity, vortices:lim, Laurent-Polz2011}, and molecules~\cite{molecules}.

{\it In this paper we use these methods to establish sufficient conditions for the stability of various branches of relative equilibria present in a mechanical system consisting of a small axisymmetric magnetic body-dipole in an also axisymmetric external magnetic field that additionally exhibits a mirror symmetry}. When the external  field is created by two magnetic poles modeled by two distant ``charges"~\cite{Smythe1939} we call this system the  {\bfi   standard orbitron}; the setup involving  arbitrary external fields exhibiting the above mentioned symmetries  will be referred to as the {\bfi  generalized orbitron}. The generic term {\bfi  orbitron} will refer simultaneously to both the standard and the generalized orbitrons. This problem has been studied for a long time already:  the model was introduced in the 1970s, and the first theoretical and experimental results were presented  in~\cite{Kozorez1974, Kozorez1981};  numerical simulations were carried out later on in~\cite{sszub08, Grygoryeva2009, sszub10}. Unfortunately, these works do not contain a complete mathematical proof of nonlinear stability due to the limitations of the  classical Lyapunov type approach that was followed. In the following pages we will see how in this case, the methods of geometric mechanics and symmetry-based stability analysis are capable of providing sets of sufficient  conditions that complement the partial results already existing in the literature and that ensure nonlinear stability.

Geometric mechanical methods have already been applied in the context of two systems involving spatially extended magnetic bodies, namely the {\it levitron} and the {\it magnetic dumbbells}. The levitron~\cite{Harrigan} is a  magnetic spinning top in the presence of gravitation that can levitate in the air repelled by a base magnet. The stability of this dynamical phenomenon has been explored with the tools of geometric mechanics in~\cite{Dullin1999, dullin:2004, levitronMarsden}. Unfortunately, in this system there are not sufficient conserved quantities available to conclude nonlinear stability using energy-momentum methods and only linear stability estimates are available. The  magnetic dumbbells~\cite{Kozorez1974} are two axisymmetric magnetic rigid bodies in space interacting contactlessly with each other; this system exhibits stable regular relative equilibria for which stability conditions have been found using the energy-momentum method in~\cite{sszub12}.

The results obtained in this paper are also of much interest at the time of clarifying several misconceptions in the physics literature that incorrectly state that purely magnetic systems cannot exhibit stable behavior in the absence of other long-range forces. This belief goes back to Earnshaw's Global Instability Theorem~\cite{Earnshaw1842}, one of the most profusely cited results in the physics literature concerning  stability in magnetic systems~\cite{Bassani2006}. Earnshaw's theory concerns mainly point particles and it was generalized during the last 170 years to a large variety of systems dealing with pure/combined confinement \cite{Simon2001} exhibiting both static and dynamic~\cite{Ginzburg1947, Tamm1979}  solutions. Such extensions were in many occasions not rigorously proved and hence experimental results in the last eighty years by Meissner~\cite{Meissner1933}, Braunbeck~\cite{Braunbek1939}, Arkadiev-Kapitsa~\cite{Arkadiev1947} (levitation with Type I superconductors), Brandt~\cite{Brandt1989}, \cite{Brandt1990} (levitation with Type II superconductors), or Harrigan~\cite{Harrigan} (levitron), raised questions as to the universal applicability of  Earnshaw's theory.

The positive stability results obtained in this paper for dynamic solutions of the orbitron lead us to believe that other similar configurations that have been experimentally observed to be stable could be rigorously proved to have this property despite widespread beliefs in the opposite direction. An important example are the 1941--1947 results by Tamm and Ginzburg that claim that in the case of two interacting magnetic dipoles, orbital motion is impossible using both  classical and  quantum mechanical descriptions~\cite{Ginzburg1947}; nevertheless, there exist experimental prototypes where a small permanent magnet exhibited  quasiorbital motion around another fixed permanent magnet for up to six minutes~\cite{Kozorez1974}. We plan to tackle these questions with methods similar to those put at work in this paper for the orbitron in a forthcoming publication.

The paper is organized as follows: in Section~\ref{The orbitron section} we present the Hamiltonian description of the orbitron by including a detailed geometric description of its phase space, equations of motion, symmetries, and associated momentum map. Section~\ref{Relative equilibria of the orbitron} contains a characterization of the relative equilibria of the orbitron that is obtained out of the critical points of the augmented Hamiltonian, constructed using the momentum map associated to the toral symmetry of this system spelled out in the preceding section. Section~\ref{Stability analysis of the relative equilibria of the orbitron} is dedicated to the stability analysis of two branches of equatorial relative equilibria introduced in Section~\ref{Relative equilibria of the orbitron}. One of these branches is singular, in the sense that it exhibits nontrivial isotropy group, and the other one is regular. The stability study is carried out for both the standard and the generalized orbitrons using the energy--momentum method, which yields in this case a set of conditions whose joint satisfaction is sufficient for the toral stability of the regular relative equilibria. Concerning the singular relative equilibria, none of these solutions  can be proved to be stable using the energy--momentum method for the standard orbitron, while in the generalized case we are able to specify sufficient conditions involving both the design parameters of the external magnetic field and the dynamical features of the system that guarantee its nonlinear stability. In the second part of Section~\ref{Stability analysis of the relative equilibria of the orbitron} we introduce new linear methods to assess the sharpness of the stability conditions; more specifically, we show that the spectral instability of a natural linearized Hamiltonian vector field that can be associated to any relative equilibrium, ensures its nonlinear instability. This result is very instrumental in our setup since it allows us, for example, to prove the nonlinear instability of the singular branch of relative equilibria of the standard orbitron and the sharpness of some  of the nonlinear stability conditions obtained in the regular case. In order to improve the readability of the paper, most proofs of the results in the paper  and a number of technical details about the geometry of the system that are used in those proofs, have been included  in  appendices  at the end of the paper (Section~\ref{Appendices}).

\medskip

\noindent {\bf Acknowledgments:} The authors thank the Fields Institute and the organizers of the Marsden Memorial Program on Geometry, Mechanics, and Dynamics that made possible the collaboration that lead to this work.  LG acknowledges financial support from the Faculty for the Future Program of the Schlumberger Foundation.

\section{The orbitron}
\label{The orbitron section}

\begin{figure}[h]
\centering
\includegraphics[scale=0.3]{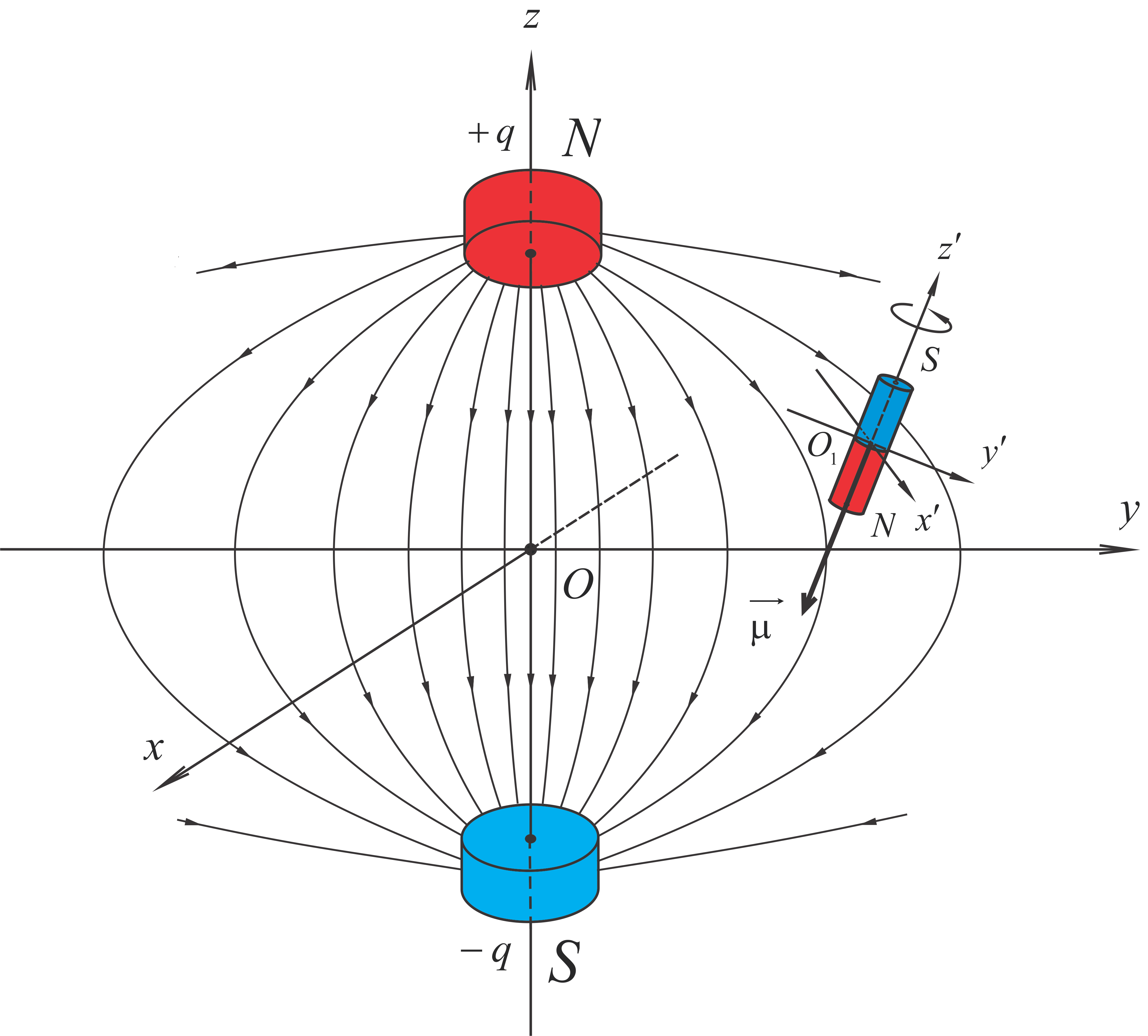}
\caption{Schematic representation of the standard orbitron. The magnetic rigid body interacts exclusively with the fixed magnetic poles represented in the picture. The opposite poles of each fixed magnet are assumed to be very distant in comparison with the dimensions of the small rigid body; therefore, their influence is  negligible  and they are hence not represented .}
\label{fig:orbitronSketch}
\end{figure}

The standard orbitron is a a small axisymmetric magnetized rigid body (for example a small permanent magnet or a current-carrying loop) with magnetic moment $\mathbf{\mu} $, in the permanent magnetic field created by two fixed magnetic poles modeled by opposite charges placed at distance $h$~\cite{Smythe1939} in the absence of gravity (see Figure \ref{fig:orbitronSketch}); in this definition the adjective ``small'' refers to the size of the body in comparison with the distance $2h$ between the magnetic poles. In this section we provide the Hamiltonian description of this physical system. 

\medskip

\noindent {\bf Phase space.} The configuration space of the orbitron is the special Euclidean group in three dimensions $SE(3) = SO(3) \times \mathbb{R}^3$. 
The $\mathbb{R}^3 $ factor of $SE(3) $ accounts for the position of the center of mass in space of the rigid body and  $SO(3)$ specifies its orientation with respect to a fixed initial frame. The orbitron is a simple mechanical system in the sense that its Hamiltonian function is of the form kinetic plus potential energy and that its phase space is the cotangent bundle $T ^\ast SE (3) $  of its configuration space $SE(3) $ endowed with the canonical symplectic structure $\omega $ obtained as minus the differential of the corresponding Liouville one form. 

As the cotangent bundle of any Lie group,  $T ^\ast SE (3) $ can be right or left trivialized in order to obtain the so called space or body coordinates, respectively (see Appendix~\ref{details phase space}), of the phase space. These trivializations provide an identification of the bundle $T ^\ast SE (3) $ with the product $SE(3)\times \mathfrak{se}(3)^\ast $, where the symbol $\mathfrak{se}(3)^\ast $ stands for the dual of the Lie algebra $\mathfrak{se}(3)$ of $SE (3) $. 

In this paper we will work in body coordinates unless it is specified otherwise. Using this representation, we denote by  $(A, \mathbf{x} )$ the elements of  $SE(3)=SO(3) \times \mathbb{R}^3$ and by $((A, \mathbf{x} ),({\bf \Pi}, \mathbf{p} )) $ those of $T ^\ast SE(3) \simeq SE(3)\times \mathfrak{se}(3)^\ast $ using body coordinates. 

\medskip

\noindent {\bf Equations of motion.} The Hamiltonian of the orbitron is given by the sum of the kinetic $T( \boldsymbol{\Pi} , \mathbf{p} ) $ and the potential $V \left(A, \mathbf{x} \right) $ energy, that is,
\begin{equation} 
\label{Hamiltonian orbitron}
h(\left(A, \mathbf{x} \right) , \left(\boldsymbol{\Pi},\mathbf{p}\right) ) = T( \boldsymbol{\Pi} , \mathbf{p} ) + V \left(A, \mathbf{x} \right).
\end{equation}
The expression of the kinetic energy is:
\begin{equation}
\label{potential and kinetic}
T(\boldsymbol{\Pi}, {\bf p} ): = \frac{1}{2} \boldsymbol{\Pi}^T  \mathbb{I}_{ref}^{-1} \boldsymbol{\Pi}  + \frac{1}{2M} \| \mathbf{p}  \|^2, 
\end{equation}
where $M$ is the mass of the axisymmetric magnetic body and  $\mathbb{I}_{ref}$ the reference inertia tensor $\mathbb{I}_{ref} = {\rm diag}(I_1,I_1,I_3)$. The coincidence between the first two principal moments of inertia is related to an axial symmetry  with respect to the third coordinate that we assume in the body.
The potential energy is given by
\begin{equation}
V \left(A, \mathbf{x} \right) := -\mu \langle \mathbf{B}(\mathbf{x}), A \mathbf{e}_3 \rangle,
\end{equation}
where $\mathbf{x} = (x, y, z)\in \mathbb{R} ^3$,  $\mu $ is the magnetic moment of the axisymmetric rigid body/dipole, and $ \mathbf{B} (\mathbf{x} )$ is the strength of the magnetic field created by two magnetic poles/``charges" $\pm q$ placed at the points $(0,0,h)$ and $(0,0,-h)$, $h>0$, that is,
\begin{equation}\label{magnetic_field}
\mathbf{B} (\mathbf{x} )= \dfrac{\mu _0 q}{4\pi}\left( \dfrac{x}{D(\mathbf{x} )_+^{3/2}} - \dfrac{x}{D(\mathbf{x} )_-^{3/2}}, \dfrac{y}{D(\mathbf{x} )_+^{3/2}}-\dfrac{y}{D(\mathbf{x} )_-^{3/2}},\dfrac{z-h}{D(\mathbf{x} )_+^{3/2}}-\dfrac{z+h}{D(\mathbf{x} )_-^{3/2}}\right),
\end{equation}
with $D(\mathbf{x} )_+=x ^2 + y ^2 + (z-h) ^2$, $D(\mathbf{x} )_-=x ^2 + y ^2 + (z+h) ^2$, and $\mu _0$ the magnetic permeability of vacuum. A small axisymmetric magnetized rigid body subjected to a external magnetic field of the form specified in~(\ref{magnetic_field}) will be called a  {\bfi  standard orbitron}.

As we will see later on, most of the results that we present in this paper hold for systems with external magnetic  fields that share the following symmetry properties presented by $\mathbf{B}$ in~(\ref{magnetic_field}), namely:
\begin{description}
\item [(i)] Equivariance with respect to rotations $R _{\theta_S} ^Z$ around the $OZ$ axis:
\begin{equation} \label{mag_equivariant}
\mathbf{B}(R_{\theta_S}^Z \mathbf{x} ) = R_{\theta_S}^Z\mathbf{B}(\mathbf{x}) \enspace {\rm for} \enspace \theta _S \in \mathbb{R}.
\end{equation}
\item [(ii)] Behavior with respect to the mirror transformation
\begin{equation}
\label{mirror transformation}
(x,y,z) \longmapsto (x,y,-z)
\end{equation}
according to the prescription
\begin{align}\label{Bx_transformation}
&B _x (x,y,z)=-B _x (x,y,-z),\\
\label{By_transformation}
&B _y (x,y,z)=-B _y (x,y,-z),\\
\label{Bz_transformation}
&B _z (x,y,z)=B _z (x,y,-z).
\end{align}
\end{description}
Consider an arbitrary magnetic field in the magnetostatic approximation in a domain free of other magnetic sources that satisfies these symmetry properties. A small axisymmetric magnetized rigid body subjected to the influence of such an external  field   will be called a {\bfi  generalized orbitron}. The generic term {\bfi  orbitron} will refer simultaneously to both the standard and the generalized orbitrons.

The equations of motion of the orbitron are determined by Hamilton's equations:
\begin{equation}
\label{Hamiltonian equations}
\mathbf{i}_{X _h} \omega= \mathbf{d} h,
\end{equation}
where $\mathbf{i}$ denotes the interior derivative, $\mathbf{d} $ is Cartan's exterior derivative, and  $X _h \in \mathfrak{X}(T ^\ast SE (3))$ the Hamiltonian vector field associated to $h \in C^{\infty} (T ^\ast SE (3))$. It can be proved (see Appendix~\ref{Equations of motion of the orbitron}) that in body coordinates, Hamilton's equations~(\ref{Hamiltonian equations}) amount to the set of differential equations:
\begin{align}
&\dot{A} = A \widehat{ \mathbb{I}_{ref}^{-1} \boldsymbol{\Pi}},\label{equations motion 1}\\
&\dot{\mathbf{x} } = \dfrac{1}{M}A {\bf p},\label{equations motion 2}\\
&\dot{\boldsymbol{\Pi} } = \boldsymbol{\Pi} \times  \mathbb{I}_{ref}^{-1} \boldsymbol{\Pi} + A^{-1} \mathbf{B}(\mathbf{x} )\times \mathbf{e}_3,  \label{equations motion 3}\\
&\dot{\mathbf{p} } = \mathbf{p} \times  \mathbb{I}_{ref}^{-1} \boldsymbol{\Pi} + \mu A^{-1} D\mathbf{B}( \mathbf{x} ) ^T A \mathbf{e}_3. \label{equations motion 4}
\end{align}
The symbol $\widehat{ \mathbb{I}_{ref}^{-1} \boldsymbol{\Pi}} $ stands for the antisymmetric matrix associated to the vector $ \mathbb{I}_{ref}^{-1} \boldsymbol{\Pi} \in \mathbb{R} ^3$ via the Lie algebra isomorphism\ \   $\widehat{} :  \left( \mathbb{R}^3, \times \right) \longrightarrow \left( \mathfrak{so}(3), [\cdot,\cdot] \right) $ introduced in Appendix~\ref{details phase space} and $D$ for the differential.

\medskip

\noindent {\bf Toral symmetry of the orbitron and associated momentum map.} The axial symmetry of the magnetic rigid body and the rotational spatial symmetry of the magnetic field created by the two poles with respect to rotations around the OZ axis endow this system with a toral symmetry which is obtained as the cotangent lift of the following action on the configuration space:
\begin{equation}
\label{toral action config space}
\begin{array}{cccc}
\Phi : &(\mathbb{T}^2=S^1\times S^1)\times SE(3)& \longrightarrow &SE(3)\\
	&\left(\left( e^{i\theta_S},e^{i\theta_B}\right) ,(A,\mathbf{x} )\right)&\longmapsto&
({R}^{Z}_{\theta_S} A {R}^{Z}_{-\theta_B},{R}^{Z}_{\theta_S}\mathbf{x} ),
\end{array}
\end{equation}
where $R _\theta ^Z$ denotes the rotation matrix around the third axis by an angle $\theta $. The first circle action involving ${R}^{Z}_{\theta_S} $ implies a spatial rotation of the center of mass of the body and the second one, given by ${R}^{Z}_{\theta_B}$, accounts for a rotation of the magnetic body around its symmetry axis. In Appendix~\ref{The toral action on phase space} we show that the cotangent lift, also denoted by $\Phi $, is a canonical symmetry given by
\begin{equation}
\label{toral action cotangent}
\begin{array}{ccccc} 
&\Phi: & (\mathbb{T}^2=S^1\times S^1)\times T^*SE(3)  &  \longrightarrow  & T^*SE(3) \\ 
&&\left(\left( e^{i\theta_S},e^{i\theta_B}\right) ,((A,\mathbf{x} ), (\boldsymbol{\Pi} , \mathbf{p} ))\right)&\longmapsto&\left( ( R_{\theta _S } A R_{-\theta _B }, R_{\theta _S } \mathbf{x} ) , (  R_{\theta _B } \boldsymbol{\Pi}, R_{\theta _B } \mathbf{p} ) \right).
\end{array} 
\end{equation}
that has an invariant momentum map associated $\boldsymbol{J} :T^*SE(3) \longrightarrow \mathfrak{t}^*$  given by:
\begin{align}
\label{momentum map toral action}
&\boldsymbol{J} \left( \left(A, \mathbf{x} \right) , \left(\boldsymbol{\Pi},\mathbf{p}\right) \right) = \left( \langle A \boldsymbol{\Pi} + \mathbf{x} \times A \mathbf{p} ,   \mathbf{e}_3 \rangle , -\langle \boldsymbol{\Pi},  \mathbf{e}_3 \rangle \right).
\end{align}
A straightforward computation shows that the Hamiltonian of the orbitron is invariant with respect to the action~(\ref{toral action cotangent}), that is,
\begin{equation*}
h \circ \Phi_{\left(e^{i\theta_S},e^{i\theta_B}\right)}=h, \quad \mbox{for any} \quad \left(e^{i\theta_S},e^{i\theta_B}\right) \in \mathbb{T}^2,
\end{equation*}
which, by Noether's Theorem~\cite[Theorem 4.2.2]{Abraham1978}, allows us to conclude that the level sets of the momentum map~(\ref{momentum map toral action}) are preserved by the associated Hamiltonian dynamics, that is, if $F _t $ is the flow of the vector field $X _h$ then $ \mathbf{J} \circ F _t = \mathbf{J} $ for any $t$.

The action~(\ref{toral action cotangent}) has two isotropy subgroups, namely, the identity $\{e\} $ and the diagonal circle $H:=\left\{\left(e^{i\theta},e^{i\theta}\right)\mid e^{i\theta} \in S ^1\right\} $. The orbit type submanifold $\left(T^*SE(3)\right) _H $  is given by
\begin{equation}
\label{singular orbit type}
\left(T^*SE(3) \right)_H = \left\{ \left(\left(R _\theta^Z, a \mathbf{e} _3 \right),\left(b \mathbf{e} _3, c \mathbf{e} _3 \right)  \right)\mid \theta, a, b, c \in \mathbb{R}\right\},
\end{equation}
and $\left(T^*SE(3) \right)_{\{e\}} =T^*SE(3) \setminus \left(T^*SE(3)\right) _H$. The bifurcation lemma (see for instance~\cite[Proposition 4.5.12]{Ortega2004}) guarantees  that the restriction of the momentum map to the regular isotropy type $T^*SE(3) _{\{e\}}$ is a submersion and that it has rank one at points in the isotropy type $T^*SE(3) _H $.

\section{Relative equilibria of the orbitron}
\label{Relative equilibria of the orbitron}

In this section we specify the equations that characterize the relative equilibria of the orbitron with respect to its toral symmetry. 

\medskip

\noindent {\bf Relative equilibria: setup and background.} Consider a vector field $X \in \mathfrak{X}(M) $ on a manifold $M$ that is equivariant with respect to action of a Lie group $G$ on it. We say that the point $m \in M $ is a relative equilibrium with velocity $\xi\in \mathfrak{g} $ if the value of vector field  at that point coincides with the infinitesimal generator $\xi_M$ associated to $\xi $, that is,
\begin{equation}
\label{relative equilibrium condition}
X _h (m)= \xi_M(m).
\end{equation}  
The Lie algebra element $\xi\in \mathfrak{g} $ is called the velocity of the relative equilibrium.
This defining property is equivalent to saying that the  flow $F _t $ associated to the vector field $X$ at the point $ m\in M $ coincides with the one-parameter Lie subgroup of $G$ generated by $\xi\in \mathfrak{g}$, that is,
\begin{equation}
\label{relative equilibrium with flow}
F _t (m)=\exp t \xi \cdot m,
\end{equation}
where $\exp $ is the Lie group exponential map $\exp : \mathfrak{g} \rightarrow G$. In the Hamiltonian setup, relative equilibria have a very convenient characterization that uses the critical points of a function instead of the equilibria of the vector field $X- \xi_M $, as in~(\ref{relative equilibrium condition}). Indeed, consider now a symmetric Hamiltonian system $(M, \omega, h, G, \mathbf{J}:M \rightarrow \mathfrak{g}^\ast) $ and assume that the momentum map $\mathbf{J} $ is coadjoint equivariant; it can be shown~\cite{Abraham1978} that the point $m \in M $ is a relative equilibrium of the Hamiltonian vector field $X _h $ with velocity $\xi\in \mathfrak{g} $ if and only if 
\begin{equation}
\label{relative critical points}
\mathbf{d} \left(h- \mathbf{J}^\xi\right)(m)=0,
\end{equation}
where $\mathbf{J}^\xi:= \langle\mathbf{J}, \xi\rangle$. The combination $h- \mathbf{J}^\xi$ is usually referred to as the {\bfi  augmented Hamiltonian}. If the relative equilibrium $m\in M$ is such that  $\mathbf{J}(m)= \mu \in \mathfrak{g}^\ast $ and we denote its isotropy subgroup with respect to the $G$ action by $G _m$, the law of conservation of the isotropy~\cite{Ortega2004} and Noether's Theorem imply~\cite[Theorem 2.8]{Ortega1999a} that $\xi \in {\rm Lie} \left(N_{G _\mu}(G _m)\right) $, where $G _\mu $ is the coadjoint isotropy of $\mu \in \mathfrak{g}^\ast$ and $N_{G _\mu}(G _m) $  is the normalizer group of $G _m $ in $G _\mu $ (note that $G _m \subset G _\mu $ necessarily due to the equivariance of the momentum map). Finally, notice that the velocity of a relative equilibrium with nontrivial isotropy is not uniquely defined; indeed, it is clear in~(\ref{relative equilibrium with flow}) that if $\xi\in \mathfrak{g}$ is a velocity for the relative equilibrium $m $, then so is $\xi+ \eta $ for any $\eta \in {\rm Lie} \left(G _m\right) $.

\medskip

\noindent {\bf Relative equilibria equations of the orbitron.} The next proposition, proved in Appendix~\ref{proof proposition relative}, specifies the critical point equations~(\ref{relative critical points}) in the case of the orbitron and shows the existence of  branches of relative equilibria whose stability we will study in the next section.

\begin{proposition}
\label{prop_rel_eq}
Consider the orbitron system introduced in Section~\ref{The orbitron section} whose Hamiltonian function is given by~(\ref{Hamiltonian orbitron}) and let 
$\mathbf{z}=\left(\left(A, \mathbf{x}\right), \left(\boldsymbol{\Pi}, {\bf p}\right)\right) \in T ^\ast SE(3) $. Then:
\begin{description}
\item [(i)] The point $\mathbf{z}$ is a relative equilibrium of the orbitron with velocity $(\xi_1, \xi_2)\in \mathbb{R} ^2$ with respect to the toral symmetry introduced in Section~\ref{The orbitron section} if and only if the following identities are satisfied:
\begin{align}
\label{equationderivative1s}
&\mu \left[  \mathbf{B} (\mathbf{x} ) \times A  \mathbf{e}_3 \right]  + \xi _1  \left[ A \mathbf{p} \times (\mathbf{x} \times \mathbf{e}_3 ) -  A \boldsymbol{\Pi} \times \mathbf{e}_3 \right] = 0, \\
\label{equationderivative2s}
&- \mu   D\mathbf{B} (\mathbf{x} ) ^T (A  \mathbf{e}_3)  - \xi _1 \left( A \mathbf{p} \times \mathbf{e}_3 \right)  = 0,\\
\label{equationderivative3s}
&\mathbb{I}_{ref}^{-1} \boldsymbol{\Pi} + \xi _2 \mathbf{e}_3 - \xi _1 A^{-1} \mathbf{e}_3 =0,\\
\label{equationderivative4s}
& \frac{1}{M} \mathbf{p} - \xi _1 A^{-1} \left(  \mathbf{e}_3 \times \mathbf{x} \right) = 0. 
\end{align}
\item [(ii)] Consider now $A_0=R_{\theta_0}^Z$,  $\mathbf{x} _0= \left( x, y, 0\right) $, $\boldsymbol{\Pi} _0 = I_3 \left( \xi _1 - \xi _2 \right) \mathbf{e}_3 $ and $\mathbf{p} _0=M \xi _1  A_0^{-1}\left( -y,x, 0\right) $. The point $\mathbf{z} _0 = \left( \left(A_0, \mathbf{x}_0 \right) , \left(\boldsymbol{\Pi}_0,\mathbf{p}_0\right) \right) $ is a relative equilibrium of the standard orbitron with velocity  $ \left( \xi _1 , \xi _2 \right) $, where $\xi _2 $ is an arbitrary real number and $\xi _1 $ is either arbitrary when $\mathbf{x}_0 =\mathbf{0} $ or 
\begin{equation}
\label{xx1 regular}
\xi _1  = \pm \left( - \dfrac{3h\mu q \mu _0 }{2 \pi MD(\mathbf{x}_0 )^{5/2}}\right)^{1/2},
\end{equation}
when $\mathbf{x}_0 \neq \mathbf{0}$. In view of the expression~(\ref{xx1 regular}) for the spatial velocity $\xi_1$, the existence of the latter relative equilibrium is only guaranteed when  $\mu q <0 $.
\item[(iii)]
The conclusions in the previous part also hold for the generalized orbitron. 
In this situation $B _z (x,y,z) = f(x ^2 + y ^2,z )$ for some $f \in C^\infty(\mathbb{R} ^2 )$, and the spatial velocity $\xi _1 $ of the relative equilibria with $\mathbf{x}_0 \neq \mathbf{0} $ is given by 
\begin{equation}
\label{xx1 generalized}
\xi _1 = \pm \left( -\dfrac{2}{M} \mu f'_1\right) ^{1/2},
\end{equation}
where $f'_1 = \left.\dfrac{\partial f(v,z )}{\partial v}\right|_{v=x ^2+ y ^2, z = 0}$.
In view of the expression of the spatial velocity $\xi_1$ in~(\ref{xx1 generalized}), the existence of this relative equilibrium is only guaranteed when  $\mu f _1 ' <0 $.
\end{description}
 The relative equilibria in these statements  for which $\mathbf{x}_0 \neq \mathbf{0} $ (respectively $\mathbf{x}_0 = \mathbf{0} $) have trivial (respectively nontrivial $H$) isotropy and hence belong  to the orbit type $\left(T^*SE(3)\right) _{\{e\}} $ (respectively $\left(T^*SE(3)\right) _{H} $); we will refer to them as {\bfi  regular relative equilibria} (respectively {\bfi  singular relative equilibria}).

\end{proposition}
\begin{figure}[!htp]
\centering
\includegraphics[scale=0.3]{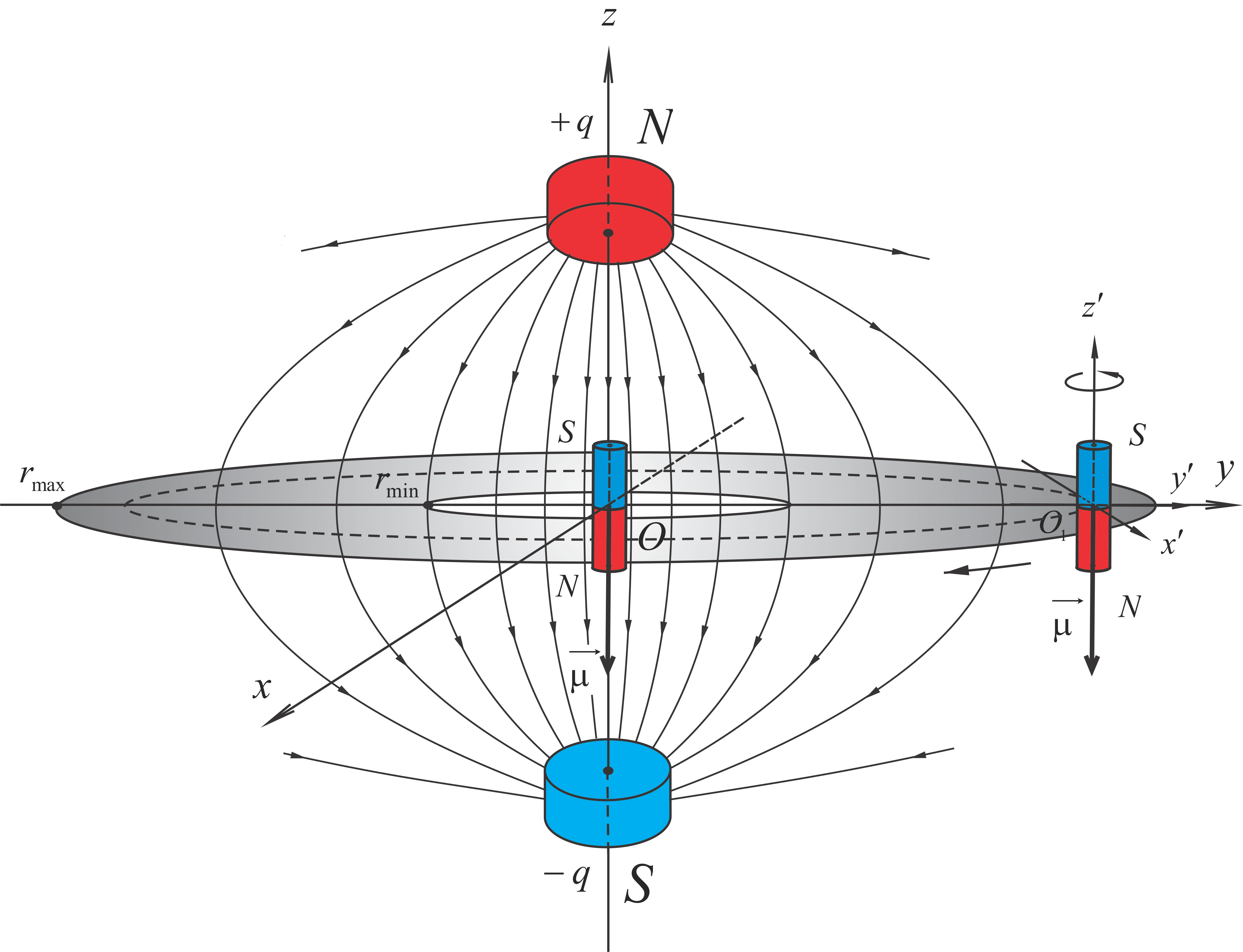}
\caption{Regular and singular relative equilibria of the standard orbitron. The symbols $r_{min} $ and $r_{max}$ represent the stability region in configuration space determined by the conditions that will be presented later on in part {\bf (i)} of Theorem~\ref{kozorez relations}.}
\label{fig:relEqulibriaSketch}
\end{figure}

\section{Stability analysis of the relative equilibria of the orbitron}
\label{Stability analysis of the relative equilibria of the orbitron}

In this section we study the stability properties of the branches of relative equilibria of the orbitron introduced in the second and third parts of Proposition~\ref{prop_rel_eq}.

\medskip

\noindent {\bf The energy--momentum method.} As we already explained in the introduction, the degeneracies present in symmetric systems cause various drift phenomena that complicate the selection of a stability criterion. The most natural and fruitful choice is that of stability relative to a subgroup, introduced in~\cite{smodg} for relative equilibria and in~\cite{po} for relative periodic orbits.
\begin{definition}
\label{stability relative to a subgroup}
Let $X\in\frak{X}(M)$ be a $G$--equivariant vector field on the $G$--manifold $M$ and let $G'$ be a subgroup of $G$.
A relative equilibrium $m\in M$ of $X$, is called
$G'${\bfi--stable}, or {\bfi stable modulo} $G'$, 
if for any $G'$--invariant open neighborhood $V$ of the 
orbit $G'\cdot m$, there is an open neighborhood $U\subset  V$ of
$m$, such that if $F_t$ is the flow of the 
vector field $X$ and $u\in U$, then $F_t (u)\in V$ 
for all $t\geq 0$.
\end{definition}

In the Hamiltonian setup there exists a variety of Dirichlet type results that provide sufficient conditions for the $G _\mu $--stability of a given relative equilibrium, where $\mu \in \mathfrak{g}^\ast $ is the momentum value in which it is sitting and $G _\mu $ is its coadjoint isotropy. The reason why the subgroup $G _\mu$ arises naturally is clear if we look at the stability problem from the symplectic reduction point of view; more explicitly, consider a symmetric Hamiltonian system $(M, \omega, h, G, \mathbf{J}: M \rightarrow \mathfrak{g}^\ast)$ that exhibits a relative equilibrium at the point $m \in M $ such that $\mathbf{J}(m)= \mu $. Suppose that the momentum map $\mathbf{J}$ is coadjoint equivariant and that the coadjoint isotropy $G _\mu$ acts freely and properly on the momentum fiber $\mathbf{J}^{-1}(\mu) $; in these conditions, the quotient space $ \mathbf{J}^{-1}(\mu) / G _\mu$ is naturally symplectic~\cite{mr,mwr} and the $G$-equivariant Hamiltonian vector field associated to $h$ projects onto another Hamiltonian vector field in which the relative equilibrium $m$ becomes a standard equilibrium. The importance of this construction in our context comes from the fact that the standard Lyapunov stability of the reduced equilibrium is equivalent to the $G _\mu $--stability of the relative equilibrium.

The following result, known as the {\bfi  energy--momentum method}, provides a sufficient condition for the $G _\mu $--stability of a given relative equilibrium. This result has been introduced at different levels of generality in~\cite{smodg, Ortega1999a, patrick:roberts:wulff:2002, MontaldiOlmos2011}.
\begin{theorem}[Energy-momentum method]
\label{energy momentum statement}
Let $(M,\omega, h)$ be a symplectic Hamiltonian system with a 
symmetry given by the Lie group $G$ acting
properly on $M$ with an
associated coadjoint equivariant momentum map $\mathbf{J} :M\rightarrow
\mathfrak{g}^\ast$. Let $m\in M$ be a relative equilibrium such that 
$\mathbf{J}(m)=\mu\in\mathfrak{g}^\ast $ and assume that the coadjoint isotropy subgroup $G _\mu $ is compact. Let $\xi\in
{\rm Lie}N_{G _\mu}(G _m)$ be a velocity of the relative equilibrium. If the quadratic form
\begin{equation}
\label{quadratic form stability}
\mathbf{d}^2(h-\mathbf{J}^\xi)(m)|_{W\times W}
\end{equation}
is definite for some (and hence for any) subspace $W$ such that
\[
\label{definitionofwre}
\ker T_m \mathbf{J}=
 W\oplus T_m(G_\mu\cdot
m),
\]
then $m$ is a $G_\mu$--stable relative equilibrium.
If $\dim W=0$, then $m$ is 
always a $G_\mu$--stable relative equilibrium. The quadratic 
form $\mathbf{d}^2(h-\mathbf{J}^\xi)(m)|_{W\times W}$, will be
called the {\bfi stability form}
of the relative equilibrium $m$ and $W$ a {\bfi  stability space}. 
\end{theorem}

\begin{remark}
\normalfont
Even though we work exclusivey in the Hamiltonian setup, this criterion has elaborate counterparts in the Lagrangian side~\cite{slm, Rodriguez-Olmos2006}.
\end{remark}

\begin{remark}
\normalfont
The statement in Theorem~\ref{energy momentum statement} can be generalized to the context of Hamiltonian actions on Poisson manifolds and can be stated so that one can take advantage of existing Casimirs or other non-symmetry related conserved quantities in order to prove the stability of a given relative equilibrium~\cite[Theorem 4.8]{pos}. More explicitly, if in the conditions of Theorem~\ref{energy momentum statement} there exists a set
of $G_\mu$--invariant conserved quantities
$C_1,\ldots,\,C_n:M\rightarrow\mathbb{R}$ for which
\[
\mathbf{d}(h-\mathbf{J}^\xi+C_1+\cdots+C_n)(m)=0,
\]
and
\[
\mathbf{d}^2(h-\mathbf{J}^\xi+C_1+\cdots+C_n)(m)|_{W\times W}
\]
is definite for some (and hence for any) $W$ such that
\[
\ker\mathbf{d} C^1(m)\cap\ldots\cap\ker\mathbf{d} C^n(m)\cap 
\ker T_m \mathbf{J}=
 W\oplus T_m(G_\mu\cdot
m),
\]
then $m$ is a $G_\mu$--stable relative equilibrium.
\end{remark}

\medskip

\noindent {\bf Nonlinear stability of the orbitron relative equilibria.} The application of the energy-momentum method to the relative equilibria of the orbitron introduced in Proposition~\ref{prop_rel_eq} makes possible the determination of sizeable regions in parameter space for which those solutions are $G _\mu $-stable ($\mathbb{T}^2 $-stable in this case). We spell this out in the statement of the following theorem whose proof is provided in the Appendix~\ref{proof kozorez}.

\begin{theorem}
\label{kozorez relations}
Consider the relative equilibria introduced in Proposition~\ref{prop_rel_eq}. Then:
\begin{description}
\item[(i)] The regular relative equilibria of the standard orbitron in part {\bf (ii)} of Proposition~\ref{prop_rel_eq}, that is, those for which $\mathbf{x}_0 \neq \mathbf{0} $, are $\mathbb{T}^2 $--stable whenever the following three inequalities are satisfied:
\begin{equation}
\label{kozoriez} 
\dfrac{2}{3}<\dfrac{r^2}{h^2} < 4,
\end{equation} 
\begin{equation}\label{signxi2}
{\rm sign}(\xi _1^0 ) I_3 \xi _2<-\mid \xi_1^0  \mid \left( I_1-I_3+\dfrac{2}{3} M \dfrac{(r^2+h^2)h^2}{3 r^2-2 h^2}\right),
\end{equation}
where $r ^2=\|\mathbf{x}_0\| ^2 $,  $\xi _1 ^0 = \pm \left( - \dfrac{3h\mu q \mu _0 }{2 \pi MD(\mathbf{x}_0 )^{5/2}}\right)^{1/2} $, and $\mu q <0$. 

The singular relative equilibria ($\mathbf{x}_0 = \mathbf{0} $) are always formally unstable, in the sense that the stability  form~(\ref{quadratic form stability})  exhibits a nontrivial signature.
\item[(ii)]
\label{kozorez relations generalized}
The regular relative equilibria of the generalized orbitron in part {\bf (iii)} of Proposition~\ref{prop_rel_eq} are $\mathbb{T}^2 $--stable whenever the following conditions hold: 
\begin{align}
& \mu  f'_1< 0,\label{f11}\\
&\mu\left(2 f_{1}'+r^2f_{1}''\right)<0,\label{kozorez generalized1}\\
&\mu f _2'' <0,\label{f22}\\
&{\rm sign}(\xi _1^0 ) I_3 \xi _2<-  |\xi_1^0| \left(  (I_1-I_3)+\dfrac{1}{2}M\left(\dfrac{f _0}{f'_1}+4 r ^2\dfrac{{f_{1}' }}{f_{2}'' } \right) \right),\label{cond3_xi2}
\end{align}
where $r ^2=\|\mathbf{x}_0\| ^2 $, $f\in C^{\infty}(\mathbb{R} ^2)$ is the function such that $B _z (x,y,z) = f(r ^2,z )$, $f _0= f (r ^2, 0)$, $f'_1 = \left. \dfrac{\partial f(v,z)}{\partial v}\right|_{v=r ^2, z= 0}$, $f''_1 = \left. \dfrac{\partial ^2 f(v,z)}{\partial v ^2}\right|_{v=r ^2, z= 0}$,  $f''_2 =  \left. \dfrac{\partial ^2 f(v,z)}{\partial z ^2}\right|_{v=r ^2, z= 0}$, and $\xi _1 ^0 = \pm{ \left( -\dfrac{2}{M} \mu f'_1\right)}^{1/2}$.

The singular branch ($\mathbf{x}_0= {\bf 0}$) is $\mathbb{T} ^2$--stable if the following conditions are satisfied:
\begin{align}
&\mu f' _1<0,\label{f11_gen_singular_cond}\\
&\mu f'' _2< 0,\label{f22_gen_singular_cond}\\
& \xi_1^2< -\dfrac{2}{M} \mu f'_1,\label{xi1_gen_singular_cond}\\
& {\rm sign}(\xi_1)\Pi_0> \dfrac{I _1\xi_1^2- \mu f_0}{|\xi_1|},\label{Pi0_gen_singular_cond}
\end{align}
where $\Pi_0= I _3\left( \xi_1- \xi_2\right)$ and we use the same notation as above for $f _0$, $f' _1$, and $f'' _2$, replacing $v = r ^2$ by $v = 0$. When $\mu f _0<0$ and $\dfrac{f _0}{f' _1}<\dfrac{2}{M} I _1$, the conditions \eqref{xi1_gen_singular_cond} and \eqref{Pi0_gen_singular_cond} can be replaced by the following single $\xi_1$--independent  optimal condition:
\begin{equation}
\label{Pi0_gen_singular_optimal_cond}
|\Pi_0|>2 \sqrt{- \mu f_0 I _1}.
\end{equation}
This optimal condition is achieved by using the spatial velocities $\xi_1=\pm  \left(-\dfrac{1}{I _1} \mu f_0\right)^{1/2}$; the positive (respectively negative) sign for the velocity corresponds to positive (respectively negative) values of $\Pi_0 $.
\end{description}
\end{theorem}

\begin{remark}
\normalfont
The right inequality in~(\ref{kozoriez}) was already known by Kozorez~\cite{Kozorez1981} but it does not ensure by itself the nonlinear stability of this symmetric configuration. We will refer to this inequality as the {\bfi  Kozorez condition}. The extension of this inequality in the context of the generalized orbitron is given by \eqref{kozorez generalized1}.
\end{remark}
\begin{remark}
\normalfont
The formal instability of the singular branch of the standard orbitron is not informative about its actual nonlinear stability or instability. This point is determined via a complementary spectral stability analysis of the linearized system that we carry out later on in Theorem~\ref{unstable branch} and that allows us to conclude the nonlinear instability of this  singular branch of relative equilibria.
\end{remark}
\begin{remark}
\label{why gaussian elimination}
\normalfont
The proof of the theorem presented in  Appendix~\ref{proof kozorez} consists of studying the definiteness of the stability form~(\ref{quadratic form stability}) introduced in Theorem~\ref{energy momentum statement}. A quick dimension count shows that the stability spaces corresponding to the regular and singular branches of relative equilibria are eight and ten dimensional, respectively. The need of determining the sign of the eigenvalues of stability forms in high dimensions like ours has motivated the introduction in the literature of various block diagonalizations for it based on arguments of  dynamic~\cite{slm, Rodriguez-Olmos2006} or kinematic~\cite{Ortega1999a} nature. An elementary but important observation that we point out in the proof of this theorem is that in order to ensure the stability of the relative equilibrium in question there is no need to compute the eigenvalues of the stability form but only to determine its signature; the relevance of this statement lies in the fact that by Sylvester's Law of Inertia, the signature is invariant by conjugation with respect to invertible matrices and hence can be read out of the pivots of the matrix obtained by performing Gaussian elimination on the stability form. Unlike the situation faced when computing eigenvalues, Gaussian elimination can be carried out formally and not just numerically in virtually any dimension. This remark is of much importance for non-simple mechanical systems for which dynamic block diagonalizations similar to those cited above are rarely available.
\end{remark}
\begin{remark}
\normalfont
Conditions \eqref{f11_gen_singular_cond}--\eqref{Pi0_gen_singular_cond} can be used in the design of magnetic fields capable of confining magnetic rigid bodies that do not exhibit spatial rotation. This is the working principle of devices such as magnetic contactless flywheels or levitrons. In the case of flywheels, up until now only actively controlled  versions have been developed; as to the levitron, the potentials that have been considered so far~\cite{Dullin1999, dullin:2004, levitronMarsden} do not allow to conclude nonlinear stability using the methods put at work in Theorem~\ref{kozorez relations} and only the spectral stability of the corresponding linearized systems has been considered. We plan to explore in detail these systems in a future publication.
\end{remark}
\medskip

\noindent {\bf Linear stability and instability analysis tools for relative equilibria.} The use of the energy-momentum method provides sufficient but not necessary nonlinear stability conditions. More specifically, there is no guarantee that the stability regions determined by the inequalities in the statement of Theorem~\ref{kozorez relations} are optimal in the sense that as soon as those conditions are violated stability disappears. In the context of stability studies for standard equilibria one usually proceeds by examining the spectral stability of the linearization at the equilibrium of the vector field in question, that is, when the sufficient stability conditions obtained via a Dirichlet type criterion are violated, one looks for eigenvalues of the linearization that exhibit a nonzero real part, whose existence would imply the nonlinear instability of the equilibrium of the original vector field. 

This way to proceed can be extended in the context of regular relative equilibria by looking at the spectral stability of the linearization of the reduced Hamiltonian vector field at the equilibrium corresponding to the relative equilibrium in the symplectic Marsden--Weinstein reduced space~\cite{mwr}. Even though in the singular case, there exist reduced spaces that generalize the Marsden--Weinstein reduced space~\cite{sl, or2006a, or2006}, the equivalence between $G _\mu $-stability of a relative equilibrium and standard nonlinear stability of the corresponding reduced equilibrium does not hold anymore, which makes necessary the formulation of a criterion that, as the energy-momentum method in Theorem~\ref{energy momentum statement}, provides a linear stability analysis tool for relative equilibria whose formulation does not need reduction; such a statement is provided in the next proposition, whose proof can be found in the appendix, and we will apply it later on to the branches introduced in Proposition~\ref{prop_rel_eq} whose nonlinear stability was studied in Theorem~\ref{kozorez relations}. In order to fix the notation and to make the presentation self contained, we start by recalling the notion of linearization of a vector field at an equilibrium point.

\begin{definition}
Let $X \in \mathfrak{X}(M)$ be a vector field on the manifold $M$ and let $m_0 \in  M$ be an equilibrium point, that is, $X(m_0) = 0$. The {\bfi  linearization} $X'$ of $X$ at the point $m_0$ is a vector field $X' \in \mathfrak{X}(T_{m_0}M)$ on the vector space $T_{m_0}M $, defined by
\begin{align*}
\begin{array}{ccc} 
X': T_{m_0}M & \longrightarrow & T_{m_0} M \times T_{m_0} M\\ 
v &\longmapsto&\left( v, \left. \frac{d}{d\lambda } \right|_{\lambda  = 0} T_{m_0} F_{\lambda } \cdot v\right),
\end{array} 
\end{align*}
where $F_{\lambda } $ is the flow of $X$. The eigenvalues of the linear map $\Pi _2 \circ X': T_{m_0}M \longmapsto T_{m_0}M$ are called the {\bfi  characteristic exponents} of $X$ at $m_0$. The map $\Pi _2:  T_{m_0} M \times T_{m_0} M \rightarrow  T_{m_0} M $ is the projection onto the second factor.
\end{definition}

\begin{proposition}
\label{linear tools for instability}
Let $G$ be a Lie group acting canonically and properly on the symplectic manifold $(M , \omega ) $ and suppose that there exists a coadjoint equivariant momentum map $ \mathbf{J}:M \rightarrow \mathfrak{g}^\ast$ that can be associated to it. Let $h \in C^\infty(M)^{G} $ be a $G$--invariant Hamiltonian and let $m \in M $ be a relative equilibrium of the corresponding $G$--equivariant Hamiltonian vector field $X _h $ with velocity $\xi\in \mathfrak{g}$. Consider a $G _m$--invariant stability space $W$ such that 
\begin{equation*}
\ker T _m\mathbf{J}=W\oplus T _m \left(G _\mu \cdot m\right),
\end{equation*}
with $\mu := \mathbf{J} (m) $ and $G _\mu \subset G $ the coadjoint isotropy of $\mu\in \mathfrak{g}^\ast$. Then:
\begin{description}
\item [(i)] $ \left(W, \omega _W\right) $ with $\omega_W:= \omega(m)| _W $ is a symplectic vector subspace of $\left(T _mM, \omega (m)\right) $.
\item [(ii)] There exists a symplectic slice $ \left(S, \omega_S\right) $ at $m\in M$ such that $\left(T _mS, \omega_S (m)\right)= \left(W, \omega _W\right) $.
\item [(iii)] The Hamiltonian vector field $X_{h _S ^\xi}\in \mathfrak{X}(S) $ in $S$ associated to the Hamiltonian function $h _S^\xi:= \left.\left(h- \mathbf{J}^\xi\right)\right | _S $ exhibits an equilibrium at the point $m \in S \subset M $.
\item [(iv)] The linearization $X'_{h _S ^\xi}\in \mathfrak{X}(T _m S) = \mathfrak{X} (W)$ of $X_{h _S ^\xi} $ at $m\in S$ coincides with the linear Hamiltonian vector field $X _Q $ on $(W, \omega _W)$ that has as Hamiltonian vector field the stability form 
\begin{equation*}
Q (w):= \mathbf{d}^2\left(h- \mathbf{J}^\xi\right)(m)(w,w), \qquad w \in W.
\end{equation*}
\item [(v)] Suppose that the two tangent spaces $T _m \left( G_{\mu } \cdot m\right) $ and  $T _m \left( G \cdot m\right) $ coincide. Then
\begin{equation} \label{TmM}
T _m M = W \oplus W ^\omega. 
\end{equation}
Additionally, let $h ^\xi:=h -\mathbf{J} ^\xi \in C^{\infty}(M) $ be the augmented Hamiltonian and let $X_{h ^\xi}' \in \mathfrak{X}(T _m M ) $ be the linearization of the Hamiltonian vector field $X_{h ^\xi}$ at $m$. 
Then
\begin{equation} \label{XQ}
X_Q = \mathbb{P}_W X_{h^ \xi }'  \boldsymbol{i} _W,
\end{equation}
where $\boldsymbol{i}_W: W \hookrightarrow T _m M$ is the inclusion, $\mathbb{P}_W: T_mM \longrightarrow W$ is the projection according to \eqref{TmM}, and $X_{h^ \xi }'$ is the linearization of $X_{h^ \xi }$ at $m$.
\item [(vi)] If the linear vector field $X _Q $ is spectrally unstable in the sense that it exhibits eigenvalues with a nontrivial real part, then the relative equilibrium $m \in M $ of $X _h $ is nonlinearly $K $--unstable, for any subgroup $K \subset G $.
\end{description}
\end{proposition}

We now provide a result that spells out how to compute the linearization of a Hamiltonian vector field at an equilibrium for systems whose phase space is the cotangent bundle of a Lie group. The following proposition  expresses the linearization that we need in terms of a linear map on the Euclidean vector space formed by the direct product of the Lie algebra and its dual.

\begin{proposition}
\label{linearization for t stars g}
Let $G$  be a Lie group with Lie algebra $ \mathfrak{g} $ and let $T ^\ast G $ be its cotangent bundle endowed with the canonical symplectic form. Consider now the body coordinates (left trivialized) expression $G \times \mathfrak{g}^\ast $ of $T ^\ast G $ and let $h \in C^{\infty}(G \times \mathfrak{g}^\ast) $ be a Hamiltonian function whose associated Hamiltonian vector field $X _h  $  exhibits an equilibrium at the point $(g, \mu) \in G \times \mathfrak{g}^\ast $. Then:
\begin{description}
\item [(i)]  Let $\varphi : G \times \left(G \times \mathfrak{g}^\ast\right) \rightarrow G \times \mathfrak{g}^\ast $ be the cotangent lift  of the action by left translations of $G$ on $G$ expressed in body coordinates. Let $h ^g:=h \circ \varphi_g $; the Hamiltonian vector field $X_{h ^g} $ exhibits an equilibrium at the point $(e,\mu)$.
\item [(ii)] Let $\Phi _g:= T_{(e, \mu)} \varphi_g: T_{(e, \mu)}\left(G \times  \mathfrak{g}^\ast\right)\simeq \mathfrak{g}\times \mathfrak{g}^\ast \longrightarrow T_{(g, \mu)}\left(G \times  \mathfrak{g}^\ast\right) $ and let $Q \in C^{\infty} \left(T_{(g, \mu)}\left(G \times  \mathfrak{g}^\ast\right)\right)$ (respectively $Q ^g \in C^{\infty}(\mathfrak{g}\times \mathfrak{g}^\ast))$ be the quadratic form associated to the second derivative of $h$ at $(g, \mu) $ (respectively of $h ^g $ at $(e, \mu)$). Then 
\begin{equation}
\label{relation quadratic forms}
Q ^g=Q \circ \Phi_g
\end{equation}
and the associated linear Hamiltonian vector fields considered as linear maps satisfy:
\begin{equation}
\label{relation linearizations}
\Phi_g\circ X_{Q ^g}= X _Q \circ \Phi_g.
\end{equation}
\item [(iii)] The linearization $X_{Q ^g}: \mathfrak{g}\times \mathfrak{g}^\ast \rightarrow \mathfrak{g}\times \mathfrak{g}^\ast $ is given by:
\begin{equation}
\label{expression linearization at e}
X_{Q ^g}(\xi, \tau)= \left(\pi_{\mathfrak{g}^\ast} \left({\rm Hess}(\xi, \tau)\right), -\pi_{\mathfrak{g}} \left({\rm Hess}(\xi, \tau)\right)+ {\rm ad}^\ast _{\pi_{\mathfrak{g}^\ast } {\rm Hess}(\xi, \tau)}\mu\right), \quad \mbox{for any} \quad \left(\xi, \tau\right)\in \mathfrak{g}\times \mathfrak{g}^\ast,
\end{equation}
where $\pi_{\mathfrak{g}} : \mathfrak{g} \times \mathfrak{g}^\ast \rightarrow \mathfrak{g} $, $\pi_{\mathfrak{g}^\ast } : \mathfrak{g} \times \mathfrak{g}^\ast \rightarrow \mathfrak{g}^\ast  $ are the canonical projections and $ {\rm Hess}:   \mathfrak{g} \times \mathfrak{g}^\ast \rightarrow \mathfrak{g} \times \mathfrak{g}^\ast $ is the linear map associated to the Hessian of $h ^g $ at $(e, \mu) $ by the relation
\begin{equation*}
\langle{\rm Hess}(\xi, \tau), (\eta, \rho)\rangle= \mathbf{d}^2h ^g(e, \mu)((\xi, \tau), (\eta, \rho)), \quad (\xi, \tau), (\eta, \rho) \in \mathfrak{g}\times \mathfrak{g}^\ast.
\end{equation*}
\end{description}
\end{proposition}

We now implement the expression for the linearization of a Hamiltonian vector field obtained in this proposition, in the particular case of the  cotangent bundle $T ^\ast SE(3)$. Let $h \in C^{\infty}(T ^\ast (SE(3))) $ be a Hamiltonian function and let $X _h $ be the corresponding Hamiltonian vector field that we assume has an equilibrium at the point $\mathbf{z} _0= \left((A _0, \mathbf{x}_0), (\boldsymbol{\Pi}_0, {\bf p} _0)\right)$, that is, $\mathbf{d} h(\mathbf{z}_0) =0$.  Let $g= (A _0, \mathbf{x} _0)\in SE(3) $ and let $\mathbf{z}= \left((I,{\bf 0}),( \boldsymbol{\Pi}_0, {\bf p} _0)\right)$; using the notation in the previous proposition, it is clear that $\mathbf{z} _0= \varphi_g (\mathbf{z}) $. Let  ${\rm Hess}(\mathbf{z} ): \mathfrak{se}(3) \times \mathfrak{se}(3)^\ast  \rightarrow  \mathfrak{se}(3) \times \mathfrak{se}(3)^\ast$ be the linear map associated to the Hessian of $h \circ \varphi_g $ at $\mathbf{z}$, that is,   for any $\mathbf{v}, \mathbf{w} \in T_{\mathbf{z}} \left( T ^\ast SE(3) \right) \simeq  \mathfrak{se}(3) \times \mathfrak{se}(3)^\ast$,
\begin{equation*}
\langle \mathbf{v}, {\rm Hess} (\mathbf{z}) \mathbf{w}\rangle = \mathbf{d} ^ 2 (h\circ \varphi _g) (\mathbf{z} ) (\mathbf{v},\mathbf{w}). 
\end{equation*}
Now, given $\mathbf{v} =\left(\delta A , \boldsymbol{\delta}\mathbf{x} , \boldsymbol{\delta}\boldsymbol{\Pi} , \boldsymbol{\delta}\mathbf{p} \right) \in \mathfrak{se}(3) \times \mathfrak{se}(3)^\ast$, define the projection (also available also for the $\boldsymbol{\delta}\mathbf{x} $, $\boldsymbol{\delta}\boldsymbol{\Pi} $, $\boldsymbol{\delta}\mathbf{p}$ components):
\begin{equation}
\begin{array}{ccc}
\boldsymbol{\pi}_{\delta A }: \mathfrak{se}(3) \times \mathfrak{se}(3)^\ast&\longrightarrow&\mathbb{R}^{3}\\
\left(  \delta A  ,  \boldsymbol{\delta}\mathbf{x} , \boldsymbol{\delta}\boldsymbol{\Pi} , \boldsymbol{\delta}\mathbf{p} \right)&\longmapsto& \delta A  
\end{array}
\end{equation}
By relations~(\ref{relation linearizations}) and~(\ref{expression linearization at e}) in Proposition~\ref{linearization for t stars g}, and the expression~(\ref{ad star}), the linearization  $X'_h$ of $X _h $ at $\mathbf{z} _0 $ is given by 
\begin{equation}
\label{expression of linearization 1}
X _h'= \Phi_g \circ X _{h ^g}' \circ \Phi_{g ^{-1}},
\end{equation}
where $X _{h ^g}':\mathfrak{se}(3) \times \mathfrak{se}(3)^\ast\simeq \mathbb{R}^{12} \rightarrow \mathfrak{se}(3) \times \mathfrak{se}(3)^\ast\simeq \mathbb{R}^{12} $ is the linear map determined by the twelve by twelve matrix
\begin{equation}
\label{expression of linearization}
X'_{h^g} = \left(
\begin{array}{c}
 \boldsymbol{\pi}  _{\boldsymbol{\delta}\boldsymbol{\Pi} } {\rm Hess} (\mathbf{z}_0   ) \\
 \boldsymbol{\pi} _{\boldsymbol{\delta}\mathbf{p} } {\rm Hess}(\mathbf{z}_0   ) \\
- \boldsymbol{\pi} _{\delta A } {\rm Hess} (\mathbf{z}_0 ) + \widehat{ \boldsymbol{\Pi} }_0 \boldsymbol{\pi} _{\boldsymbol{\delta}\boldsymbol{\Pi}} {\rm Hess}(\mathbf{z}_0   ) + \widehat{ \mathbf{p} }_0  \boldsymbol{\pi} _{\boldsymbol{\delta}\mathbf{p} } {\rm Hess}(\mathbf{z}_0  ) \\
-\boldsymbol{\pi} _{\boldsymbol{\delta}\mathbf{x} } {\rm Hess}(\mathbf{z}_0  )+ \widehat{ \mathbf{p} }_0  \boldsymbol{\pi} _{\boldsymbol{\delta}\mathbf{\Pi} } {\rm Hess}(\mathbf{z}_0  )
\end{array}
\right). 
\end{equation}
This expression should be understood as a vertical concatenation of four matrices with three rows and twelve columns each. More explicitly, given that for any  $\mathbf{v} =\left(\widehat{\delta A} , \boldsymbol{\delta}\mathbf{x} , \boldsymbol{\delta}\boldsymbol{\Pi} , \boldsymbol{\delta}\mathbf{p} \right) \in \mathfrak{se}(3) \times \mathfrak{se}(3)^\ast$, $\Phi _g (\mathbf{v})=\left(A _0\widehat{\delta A} , A _0\boldsymbol{\delta}\mathbf{x} , \boldsymbol{\delta}\boldsymbol{\Pi} , \boldsymbol{\delta}\mathbf{p} \right)\in T _{\mathbf{z} _0} \left(T ^\ast SE(3)\right)$, we can write
\begin{equation*}
X _h'\left(A _0\widehat{\delta A} , A _0\boldsymbol{\delta}\mathbf{x} , \boldsymbol{\delta}\boldsymbol{\Pi} , \boldsymbol{\delta}\mathbf{p} \right)= \left(A _0 X _A, A _0 X _{\mathbf{x}}, X_{\boldsymbol{\Pi}}, X_{\boldsymbol{p}}\right),
\end{equation*}
where $\left( X _A,  X _{\mathbf{x}}, X_{\boldsymbol{\Pi}}, X_{\boldsymbol{p}}\right) $ is the image by~(\ref{expression of linearization}) of the vector $\left(\delta A , \boldsymbol{\delta}\mathbf{x} , \boldsymbol{\delta}\boldsymbol{\Pi} , \boldsymbol{\delta}\mathbf{p} \right) $.

\medskip

\noindent {\bf Linear stability and instability of the orbitron relative equilibria.} The results presented in the previous paragraph provide all the necessary tools to carry out the linear stability analysis of the relative equilibria of the standard and generalized orbitron introduced in the  parts {\bf (ii)} and {\bf (iii)}  of Proposition~\ref{prop_rel_eq}. We will proceed by using expressions~(\ref{expression of linearization 1}) and~(\ref{expression of linearization}) in order to compute the linearization at the relative equilibria of the Hamiltonian vector fields associated to the augmented Hamiltonians constructed with the appropriate relative equilibrium velocities. We subsequently use part {\bf (v)} of Proposition~\ref{linear tools for instability} in order to write down the linearization of the restriction of this vector field to the tangent space to a symplectic slice (equivalently, a stability space); finally, we use the last part of this result in order to search for instability regions by looking for eigenvalues of this linearization that exhibit a nontrivial real part and determine how sharp the nonlinear sufficient stability conditions in Theorem~\ref{kozorez relations} are; more specifically, we will see that  there might exist relative equilibria that are nonlinearly stable even though the conditions in Theorem~\ref{kozorez relations} are not satisfied. A detailed description of this implementation is provided in Appendix~\ref{Linear stability and instability of the orbitron relative equilibria}. The following result, formulated using the terminology introduced in Proposition~\ref{prop_rel_eq}, summarizes the results of the linear analysis.

\begin{theorem}
\label{unstable branch}
Consider the relative equilibria introduced in Proposition~\ref{prop_rel_eq}. Then:
\begin{description}
\item [(i)] Regarding the relative equilibria of the standard orbitron in part {\bf (ii)} of the proposition:
\begin{description}
\item [(a)] The regular relative equilibria that do not satisfy the Kozorez relation ($r^2/h^2 < 4$) are unstable and  this stability condition is consequently sharp. The other two stability conditions in~(\ref{kozoriez}) and~(\ref{signxi2}) are not sharp, that is, there are regions in parameter space that do not satisfy them and where the linearized system is spectrally stable.
\item [(b)] The singular relative equilibria of the standard orbitron are nonlinearly unstable.
\end{description} 
\item [(ii)] Regarding the relative equilibria of the generalized orbitron in  part {\bf (iii)} of the proposition:
\begin{description}
\item [(a)] The regular relative equilibria that do not satisfy the generalized Kozorez relation~(\ref{kozorez generalized1}), namely,  
$\mu \left(2f _1 ' + r ^2f _2''\right)<0$, 
are unstable and  this stability condition is consequently sharp. The remaining stability conditions~(\ref{f11}),~(\ref{f22}), and~(\ref{cond3_xi2}) are not sharp, that is, there are regions in parameter space that do not satisfy them and where the linearized system is spectrally stable.
\item [(b)] The spectral stability of the singular relative equilibria of the generalized orbitron is equivalent to the following three conditions:
\begin{align}
& \mu f _1'<0,\label{sing_gen_f11}\\
&\mu f _2''<0,\label{sing_gen_f22}\\
&\Pi_0^2> -4 \mu f _0I _1,\label{sing_gen_Pi0}
\end{align}
where $\Pi _0=I _3(\xi_1- \xi_2)$. This statement implies that the nonlinear stability conditions~(\ref{f11_gen_singular_cond}) and~(\ref{f22_gen_singular_cond}) are sharp and that the remaining conditions are not.
\end{description} 
\end{description}
\end{theorem}

\begin{figure}[!htp]
\includegraphics[scale=0.3]{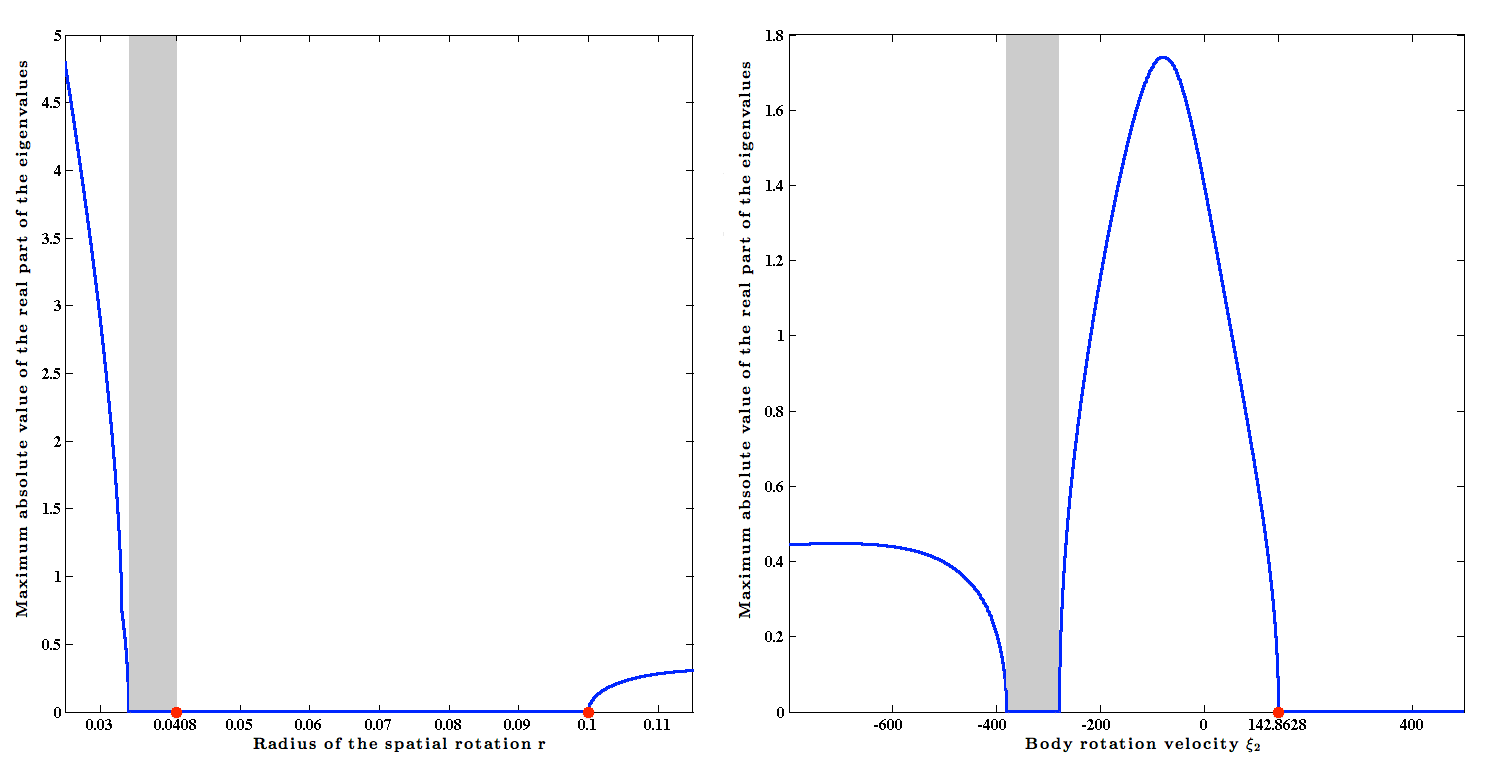}
\caption{Spectral stability study for the relative equilibria of a standard orbitron with $h = 0.05$ m, $M = 0.0068$ kg, $\mu_0 = 4\pi\cdot10^{-7}$~N$\cdot$A$^{-2}$, $\mu = -0.18375$ A$\cdot$m$^2$, $q = 17.58$ A$\cdot$m, $I_1 = 0.17\cdot10^{-6}$ kg$\cdot$m$^2$, $I_3=0.1\cdot10^{-6}$ kg$\cdot$m$^2$. The position of the red bullets indicates the critical values of $r$ (m) and $\xi_2 $ (rad$\cdot $s$^{-1}$) determined by the stability conditions~(\ref{kozoriez}) and~(\ref{signxi2}), respectively. The grey bands correspond to the stability gaps discussed in the proof of Theorem~\ref{unstable branch}, part {\bf (i)} in which the system is spectrally stable while the stability form exhibits a nontrivial signature.}
\label{fig:stabilityGaps}
\end{figure}

\noindent\textbf{Proof.\ \ (i) Part {\bf (a)}} The linearization $X _Q $ at the regular relative equilibria of the Hamiltonian vector field in the stability space associated to the augmented Hamiltonian is provided in the expression~(\ref{linearization regular}). This matrix is block diagonal and the top two by two block has as eigenvalues 
\begin{equation*}
\lambda_{\pm}=\pm \xi_1^0\sqrt{\frac{r^2-4h^2}{r^2+h^2}},
\end{equation*}
which are real whenever $r^2/h^2 >4$, that is, when the Kozorez relation is violated. In conclusion, the part {\bf (vi)} of Proposition~\ref{linear tools for instability} ensures that as soon as the Kozorez relation is violated the relative equilibria cease to be stable. The lack of sharpness of the two other stability conditions in~(\ref{kozoriez}) and~(\ref{signxi2}) is observed by studying the spectrum of the remaining six by six block of the linearization $X _Q$ which may be purely imaginary in regions of the parameter space in which those conditions are violated. The expressions corresponding to those six eigenvalues are very convoluted and we therefore do not include them in the paper; in turn, we illustrate this phenomenon in Figure~\ref{fig:stabilityGaps}, in which we plot the maximum absolute value of the real part of the eigenvalues of the linearization versus the radius of spatial rotation $r$ and the body rotation velocity $\xi_2$, respectively, when all the system parameters specified in the caption remain constant. The graph on the left hand side shows that when the radius goes beyond the critical value  stipulated by the left inequality in~(\ref{kozoriez}) the spectrum of the linearization remains purely imaginary for a while and the system is hence potentially stable; it is also visible that, as we proved above, the system becomes spectrally unstable as soon as the Kozorez relation ceases to be satisfied. The lack of sharpness of the condition~(\ref{signxi2}) is illustrated in the right hand side graph and is of a slight different nature; indeed, as soon as the condition is not satisfied, spectral instability appears but if the body rotation velocity is sufficiently decreased the system  becomes again spectrally stable in some interval of the $\xi_2$ parameter space. 

\medskip

\noindent\textbf{(i) Part {\bf (b)}} The corresponding linearization $X_Q$ at the singular relative equilibria  is described  in~(\ref{xq singular}). Its spectrum includes the two following eigenvalues:
\begin{eqnarray*}
\lambda _1 &=& \dfrac{1}{h ^2 }\sqrt{-\dfrac{3 \mu _0 \mu q }{M \pi }},\\
\lambda _2 &=& \sqrt{ - \left( \xi _1 - \dfrac{1}{h ^2 }\sqrt{-\dfrac{3\mu _0 \mu q}{2 M \pi }} \right) ^2 }.
\end{eqnarray*}
The eigenvalue $\lambda _1 $ can only be purely imaginary when $\mu q > 0$. This in turns implies that the term $\sqrt{-\dfrac{3 \mu _0 \mu q}{2M \pi }}$ in $\lambda _2 $ is purely imaginary and prevents the eigenvalue to be purely imaginary unless ${-\dfrac{3 \mu _0 \mu q}{2M \pi }}$ is zero. 

\medskip

\noindent {\bf (ii) Part (a)} Analogously to the situation in the proof of {\bf (i)}\ Part {\bf (a)}, the linearization $X _Q$ at the regular relative equilibria of the generalized orbitron exhibits the following two eigenvalues:
\begin{equation*}
\lambda_{\pm}= \pm 2\sqrt{\dfrac{1}{M}\mu\left(2 f _1'+ r ^2 f _1 ''\right)},
\end{equation*}
which are obviously purely imaginary if and only if the generalized Kozorez relation~(\ref{kozorez generalized1}) holds. The lack of sharpness in the remaining relations follows from the fact that they contain as particular cases the stability conditions for the standard orbitron that, as we illustrated in Figure~\ref{fig:stabilityGaps}, are not necessary for the spectral stability of $X _Q $. 

\medskip

\noindent {\bf (ii) Part (b)} The linearization $X _Q $ at the singular relative equilibria of the generalized orbitron is provided in~(\ref{xq singular generalized}) and its spectrum is made up by the following ten eigenvalues:
\begin{align}
\lambda_1 ^{\pm} &= \pm \sqrt{\dfrac{1}{M} \mu f ''_2},\label{sing_gen_lambda1}\\
\lambda _{2,\pm} ^{\pm} &= \pm \sqrt{-\dfrac{1}{M} \left( \xi _1 \sqrt{M} \pm \sqrt{-2 \mu f ' _1 }\right) ^2  },\label{sing_gen_lambda2}\\
\lambda _{3,\pm} ^{\pm} &= \pm \dfrac{1}{2}\sqrt{-\dfrac{1}{I _1 } \left( \left( 2 \xi _1 I _1 - \Pi _0 \right) \pm \sqrt{4 \mu f _0 I _1 + \Pi _0 ^2 }   \right) ^2 }.\label{sing_gen_lambda3}
\end{align}
The eigenvalues $\lambda_1 ^{\pm} $ can be purely imaginary only when $\mu f '' _2 <0$.  In order for the four eigenvalues $\lambda _{2,\pm} ^{\pm}$ to have the same property, the term $\sqrt{-2 \mu f ' _1 }$ has to be necessarily a real number, which yields the condition $\mu f' _1 <0$. These two relations obviously imply that the nonlinear stability conditions~(\ref{f11_gen_singular_cond}) and~(\ref{f22_gen_singular_cond}) are sharp.
Finally, the remaining four eigenvalues $\lambda _{3,\pm} ^{\pm} $ are purely imaginary whenever the term $\sqrt{4 \mu f _0 I _1 + \Pi _0 ^2 }$ is real, which requires in turn that the relation $\Pi _0 ^2 >4 \mu f _0 I _1$ is satisfied. We note that this relation may hold without~(\ref{xi1_gen_singular_cond}) and~(\ref{Pi0_gen_singular_cond}) or~(\ref{Pi0_gen_singular_optimal_cond}) being satisfied. Indeed, take for example a system for which $\mu f _0<0 $; in that situation, the relation~(\ref{sing_gen_Pi0}) does not impose any constraint on $\Pi_0 $ an hence it is easy to find values for this variable that violate ~(\ref{xi1_gen_singular_cond}) and~(\ref{Pi0_gen_singular_cond}) or~(\ref{Pi0_gen_singular_optimal_cond}). 
\quad $\blacksquare$

\section{Appendices}
\label{Appendices}

\subsection{The geometry of the phase space of the orbitron $(T ^\ast SE (3), \omega) $}
\label{details phase space}

\noindent {\bf Lie group and Lie algebra structure of the configuration space.} The configuration space of the orbitron is the Lie group
$SE(3) = SO(3) \times \mathbb{R}^3$ endowed with the semidirect product structure associated to the composition rule
\begin{equation} 
\label{Lie product}
\begin{array}{cccc}
\Psi:& SE(3)\times SE(3)&\longrightarrow &SE(3) \\
    		&(({A}_1,\mathbf{x}_1),({A}_2,\mathbf{x}_2)) &\longmapsto &({A}_1 {A}_2, {A}_1 \mathbf{x}_2 + \mathbf{x}_1),
\end{array}
\end{equation} 
for which $e=(I,0)$ and  $({A},\mathbf{x})^{-1}= ({A}^{-1},-{A}^{-1}\mathbf{x})$. In order to spell out the Lie algebra structure associated to the Lie product~(\ref{Lie product}) we start by recalling the Lie algebra isomorphism \ \   $\widehat{} :  \left( \mathbb{R}^3, \times \right) \longrightarrow \left( \mathfrak{so}(3), [\cdot,\cdot] \right) $ between  the Lie algebra $(\mathfrak{so}(3),[\cdot ,\cdot ])$ of $SO(3)$ and $(\mathbb{R} ^3, \times) $ endowed with the standard cross product, given by the assignment
\begin{equation*} 
\mathbf{x} = 
(x_1,x_2,x_3)\in  \mathbb{R} ^3 \longmapsto
\widehat{\mathbf{x}}:=\begin{pmatrix}
               0 & -x_3 &  x_2 \\
             x_3 &    0 & -x_1 \\
            -x_2 &  x_1 &    0
\end{pmatrix}.
\end{equation*}
We recall that isomorphism $\enspace \widehat{}\enspace$  satisfies that  $\widehat{\mathbf{x}} \mathbf{w} = \mathbf{x} \times \mathbf{w}$ and that  for any $A \in SO(3)$ and $\mathbf{x} \in \mathbb{R}^3$
\begin{align}
& T_IL_{A} \cdot \widehat{\mathbf{x}} = A \widehat{\mathbf{x}},\label{elementary relation 1}\\
& \mathrm{Ad}_{A} \widehat{\mathbf{x}} = A \widehat{\mathbf{x}} A^{-1} = \widehat{A\mathbf{x}},\label{elementary relation 2}\\
& \mathrm{Ad}_{A} \widehat{\mathbf{x}} = T_I \left( L _A \circ R_{A^{-1}}  \right)  \widehat{\mathbf{x}} = \left. \frac{d}{dt} \right|_0 A \exp  t
\widehat{\mathbf{x}} A^{-1}=A \widehat{ \mathbf{x} } A^{-1}= \widehat{ A \mathbf{x} },\label{elementary relation 3}
\end{align}
where $L _A: SO(3) \rightarrow SO(3) $ (respectively $R _A $) denotes left (respectively right) translations and $ \mbox{\rm Ad} _A : \mathfrak{so}(3) \rightarrow \mathfrak{so}(3) $ is the adjoint representation. The $\enspace \widehat{} \enspace$ isomorphism induces another one 
\begin{align*} 
\enspace \widehat{}\enspace :\  &\mathbb{R}^3 \longrightarrow \mathfrak{so}(3)^*\nonumber \\
&\boldsymbol{\pi}  \longmapsto \widehat{\boldsymbol{\pi} } 
\end{align*}
uniquely determined by the relation $\langle \widehat{\boldsymbol{\pi}  }, \widehat{ \mathbf{x}}  \rangle :=\langle \boldsymbol{\pi}   , \mathbf{x} \rangle _{\mathbb{R}^3} $, with $\langle \cdot    , \cdot  \rangle _{\mathbb{R}^3} $ the Euclidean inner product in $\mathbb{R} ^3$.
Using this isomorphism, we have
\begin{equation} \label{coadj}
\mathrm{Ad}^*_A \widehat{ \boldsymbol{\pi} } = \widehat{ A ^{-1} \boldsymbol{\pi}  }.
\end{equation}
Using this notation, the Lie algebra structure of $ \mathfrak{se}(3) = \mathfrak{so}(3)\times \mathbb{R}^3$ is given by the bracket 
\begin{equation}
\left[ \left( \widehat{ \boldsymbol{\rho}_1 }, \boldsymbol{\tau} _1 \right) , \left( \widehat{\boldsymbol{\rho}_2},\boldsymbol{ \tau}_2  \right) \right] := \left( \widehat{\boldsymbol{ \rho}_1 \times \boldsymbol{\rho}_2} , \boldsymbol{\rho}_1 \times \boldsymbol{\tau}_2 - \boldsymbol{\rho}_2 \times \boldsymbol{\tau}_1 \right). 
\end{equation}
Additionally, for any $(A, \mathbf{x})\in SE(3) $, $\left( \widehat{\boldsymbol{\rho}}, \boldsymbol{\tau} \right), \left( \widehat{\boldsymbol{\rho}} _1 , \boldsymbol{\tau} _1 \right), \left( \widehat{\boldsymbol{\rho} } _2 , \boldsymbol{\tau} _2 \right) \in \mathfrak{se}(3)$, $\left(\widehat{ \boldsymbol{\mu}} , \boldsymbol{\alpha} \right)\in \mathfrak{se}(3)^\ast $, $\boldsymbol{\beta}, \boldsymbol{\gamma} \in \mathbb{R} ^3$, the following relations that we use later on in the paper hold
\begin{align}
T_{\left( I, 0\right)}L_{\left( A, \mathbf{x}  \right) } \cdot \left( \widehat{\boldsymbol{\rho}} , \boldsymbol{\tau} \right) & = \left( A \widehat{\boldsymbol{\rho}}, A \boldsymbol{\tau} \right) \label{leftact}\\
T_{\left( I , 0\right) }R_{\left( A, \mathbf{x} \right) } \cdot \left( \widehat{\boldsymbol{\rho}} , \boldsymbol{\tau} \right) &= \left( \widehat{\boldsymbol{\rho} } A, \boldsymbol{\rho} \times \mathbf{x} + \boldsymbol{\tau}  \right) \label{rightact}\\
{\rm ad} _{\left( \widehat{\boldsymbol{\rho} _1} , \boldsymbol{\tau}_1 \right) }\left( \widehat{\boldsymbol{\rho} } _2 , \boldsymbol{\tau} _2 \right)  &= \left( \widehat{ \boldsymbol{\rho} _1 \times \boldsymbol{\rho} _2 }, \boldsymbol{\rho} _1 \times \boldsymbol{\tau} _2 - \boldsymbol{\rho} _2 \times \boldsymbol{\tau} _1 \right), \\
{\rm ad}^* _{\left( \widehat{ \boldsymbol{\rho} } , \boldsymbol{\tau} \right) } \left(\widehat{\boldsymbol{\mu}} , \boldsymbol{\alpha} \right) &= \left( \widehat{ \boldsymbol{\mu}  \times \boldsymbol{\rho} } + \widehat{ \boldsymbol{\alpha} \times \boldsymbol{\tau} }, \boldsymbol{\alpha} \times \boldsymbol{\rho} \right), \label{ad star} \\
T^*_{\left( I,0\right) } R_{\left( A, \mathbf{x} \right) } \left( \widehat{ \boldsymbol{\beta}} A,\boldsymbol{\gamma}\right) &= \left( \boldsymbol{\beta}+ \mathbf{x} \times \boldsymbol{\gamma}, \boldsymbol{\gamma}\right),  \\
T^*_{\left( I,0\right) } L_{\left( A, \mathbf{x} \right) } \left( A \widehat{ \boldsymbol{\beta}}, \boldsymbol{\gamma}\right) &= \left( \boldsymbol{\beta}, A^{-1}\boldsymbol{\gamma}\right).  \label{t star}
\end{align}
In the last two expressions we have identified $T ^\ast _{(A, \mathbf{x})}SE(3) $ with $T  _{(A, \mathbf{x})}SE(3) $ using the Frobenius norm in the $SO(3) $ part and the Euclidean norm in the $\mathbb{R}^3 $ part.
Using these equalities, it is easy to see that the adjoint and coadjoint actions of $SE(3)$ on its algebra $\mathfrak{se}(3)$ and its dual $\mathfrak{se}(3)^\ast $ are determined by:
\begin{align}
{\rm Ad}_{(A, \mathbf{x}) } \left( {\widehat{ \boldsymbol{\rho} }, \boldsymbol{\tau} } \right) &= \left( {\rm Ad}_A \widehat{\boldsymbol{\rho}}, - ({\rm Ad } _A \widehat{\boldsymbol{\rho} }) \mathbf{x} + A \boldsymbol{\tau} \right) = \left( \widehat{ A \boldsymbol{\rho}}, \mathbf{x} \times A \boldsymbol{\rho} + A \boldsymbol{\tau} \right),   \\
{\rm Ad}^*_{(A, \mathbf{x} )} \left( {\widehat{\boldsymbol{\mu}  }, \boldsymbol{\alpha}  } \right) &= \left( {\rm Ad}^*_A \widehat{\boldsymbol{\mu}  } - \widehat{\left(  A^{-1} \left( \mathbf{x} \times \boldsymbol{\alpha }\right) \right)} , A^{-1} \boldsymbol{\alpha} \right) = \left( \left( \widehat{ A^{-1} ( \boldsymbol{\mu} - ( \mathbf{x} \times \boldsymbol{\alpha} ))} \right), A^{-1} \boldsymbol{\alpha} \right).\label{coadjaction}
\end{align}

\medskip

\noindent {\bf Body and space coordinates for $T ^\ast SE (3) $.} Given an arbitrary Lie group $G$ with Lie algebra $\mathfrak{g} $, we recall (see for example~\cite{Abraham1978}) that the maps 
\begin{equation}
\label{maps space coordinates}
\begin{array}{cccc}
\varrho_1:& TG &\longrightarrow &G \times g \\
	&u_g &\longmapsto &(g,T_g R_{g^{-1}}\cdot u_g) \\
	&T_e R_g\cdot \xi &\longmapsto &(g,\xi).
\end{array}\qquad \mbox{and} \qquad
\begin{array}{cccc}
\varrho_2: &T^*G &\longrightarrow &G \times g^* \\
	&\alpha_g &\longmapsto &(g,T_e^* R_{g}\cdot \alpha_g) \\
	&T^*_g R_{g^{-1}}\cdot \mu  &\longmapsto  &(g,\mu)
	\end{array}
\end{equation}
define trivializations of the tangent $TG$ and cotangent bundles $T ^\ast G $, respectively, that are usually referred to as {\bfi   space coordinates} of these bundles. Notice that if $\varrho _1 \left( u _g \right) = \left( g, \xi \right)$, $\varrho _2 \left( \alpha  _g \right) = \left( g, \mu \right)$, then $\left\langle \alpha _g , u _g \right\rangle = \left\langle \mu , \xi \right\rangle $.

Analogously, the trivializations obtained using left translations instead via the maps 
\begin{equation}
\label{maps body coordinates}
\begin{array}{cccc}
\lambda_1:& TG &\longrightarrow &G \times g \\
	&u_g &\longmapsto &(g,T_g L_{g^{-1}}\cdot u_g) \\
	&T_e L_g\cdot \xi &\longmapsto &(g,\xi).
\end{array}\qquad \mbox{and} \qquad
\begin{array}{cccc}
\lambda_2: &T^*G &\longrightarrow &G \times g^* \\
	&\alpha_g &\longmapsto &(g,T_e^* L_{g}\cdot \alpha_g) \\
	&T^*_g L_{g^{-1}}\cdot \mu  &\longmapsto  &(g,\mu)
	\end{array}
\end{equation}
are usually referred to as {\bfi   body coordinates}.
Notice that if $\lambda  _1 \left( u _g \right) = \left( g, \xi \right)$, $\lambda _2 \left( \alpha  _g \right) = \left( g, \mu \right)$, then $\left\langle \alpha _g , u _g \right\rangle = \left\langle \mu , \xi \right\rangle $.

We now use these maps to establish the relation between the space and body  coordinates $\left( (A, \mathbf{x} ), ( \boldsymbol{\Pi}_S, \mathbf{p}  _S) \right) $ and $\left( (A, \mathbf{x} ), ( \boldsymbol{\Pi}_B, \mathbf{p}  _B) \right) $, respectively, of an arbitrary point in $T^*SE(3)$. Indeed, using~(\ref{maps body coordinates}),~(\ref{maps space coordinates}), and  \eqref{coadjaction}, we have that
\begin{align*}
\left( (A, \mathbf{x} ), (\boldsymbol{\Pi} _B, \mathbf{p} _B)\right)& = \lambda _2 \left( T^*_{(A, \mathbf{x} )} R_{(A, \mathbf{x} )^{-1}} \cdot (\boldsymbol{\Pi} _S, \mathbf{p} _S)\right) = \left( (A, \mathbf{x} ), {\rm Ad}^* _{(A, \mathbf{x} )} (\boldsymbol{\Pi} _S, \mathbf{p} _S)\right) \\
&= \left( (A, \mathbf{x} ), (A^{-1} (\boldsymbol{\Pi} _S - \mathbf{x} \times \mathbf{p} _S), A^{-1} \mathbf{p} _S)\right).
\end{align*}
Consequently,
\begin{align} \label{spacetobody}
&\boldsymbol{\Pi} _B = A^{-1} \left( \boldsymbol{\Pi} _S - \mathbf{x} \times \mathbf{p} _S\right),\nonumber \\
& \mathbf{p} _B = A^{-1} \mathbf{p} _S.
\end{align}
Conversely,
\begin{align} \label{bodytospace}
&\boldsymbol{\Pi} _S = A \boldsymbol{\Pi} _B + \mathbf{x} \times A \mathbf{p} _B, \nonumber \\
& \mathbf{p} _S = A \mathbf{p} _B.
\end{align}

\subsection{Equations of motion of the orbitron}
\label{Equations of motion of the orbitron}
In this section we obtain the equations of motion~(\ref{equations motion 1})-(\ref{equations motion 4}) of the orbitron using body coordinates. We will proceed by writing down first the differential equations that define a Hamiltonian vector field on the left trivialized cotangent bundle $G \times \mathfrak{g}^\ast$ of an arbitrary Lie group $G$ with Lie algebra $\mathfrak{g}$.

\begin{proposition}
Let $G$ be a Lie group with Lie algebra $\mathfrak{g} $ and let $T ^\ast G $ be its cotangent bundle endowed with the canonical symplectic form. Let $\omega _B $ be the corresponding symplectic form on the trivial bundle $G \times \mathfrak{g}^\ast $ obtained out of $T ^\ast G $ by left trivialization (body coordinates) and let $h\in C^{\infty}(G \times \mathfrak{g}^\ast)$ be a Hamiltonian function.  For any $(g, \mu) \in G \times \mathfrak{g}^\ast $, the Hamiltonian vector field $X _h\in \mathfrak{X}(G \times \mathfrak{g}^\ast)$ associated to $h$ is given by 
\begin{equation}
\label{h from partial}
X _h \left( g, \mu \right) = \left( T _I L _g  \cdot X _G \left( g, \mu \right) , X _{\mathfrak{g} ^\ast } \left( g, \mu \right) \right),
\end{equation}
where $X_{G}(g, \mu) \in \mathfrak{g} $ and $ X _{\mathfrak{g} ^\ast } \left( g, \mu \right) \in \mathfrak{g}^\ast $ are determined by
 \begin{eqnarray}
 X_G  \left( g, \mu \right) &=& D_{\mathfrak{g} ^\ast } h  \left( g, \mu \right),\label{xG}\\
 X_{\mathfrak{g} ^\ast }  \left( g, \mu \right) &=& - T_I ^\ast L_g \cdot D_G h  \left( g, \mu \right) + {\rm ad} ^\ast _{D_{\mathfrak{g} ^\ast } h  \left( g, \mu \right)} \mu.\label{xgstar}
 \end{eqnarray}
\end{proposition}

\noindent\textbf{Proof.\ \ }
Using the expression of the canonical symplectic form $\omega_B$ of $T ^\ast G $ in body coordinates (see for instance~\cite[Expression (6.2.5)]{Ortega2004}) it is easy to see that $X _G$, $X _{\mathfrak{g}^\ast} $, and hence $X _h $, are determined by the relation
\begin{multline*}
\omega_B(g, \mu)\left(X _h \left( g, \mu \right), \left(T _I L _g  \cdot \xi _G, \beta\right)\right)=\langle \beta , X _G \left( g, \mu \right)  \rangle - \langle X_{g ^\ast } \left( g, \mu \right), \xi _G \rangle \\
+ \langle \mu , \left[ X _G  \left( g, \mu \right), \xi _G\right] \rangle = D_G h  \left( g, \mu \right) \cdot T _I L _g  \cdot \xi _G + D_{\mathfrak{g} ^\ast } h  \left( g, \mu \right) \cdot \beta,  
\end{multline*}
where  $\xi _G \in \mathfrak{g} $ and $\beta \in \mathfrak{g} ^\ast $ are arbitrary and $D_G h $ and $D_{\mathfrak{g}^\ast}h $ are the partial derivatives of $h$ with respect to $G$ and $\mathfrak{g}^\ast $, respectively.
 Equivalently,
 \begin{eqnarray*}
 X_G  \left( g, \mu \right) &=& D_{g ^\ast } h  \left( g, \mu \right),\\
 X_{g ^\ast }  \left( g, \mu \right) &=& - T_I ^\ast L_g \cdot D_G h  \left( g, \mu \right) + {\rm ad} ^\ast _{D_{g ^\ast } h  \left( g, \mu \right)} \mu,
 \end{eqnarray*}
 as required. \quad $\blacksquare$
 
\medskip 
 
We now consider the case we are interested in, that is, $G = SE(3) = SO(3) \times \mathbb{R} ^3 $ and 
\begin{equation}
\label{hamiltonian function prop}
h \left( (A, \mathbf{x} ), (\boldsymbol{\Pi} , \mathbf{p} )\right) = \dfrac{1}{2} \boldsymbol{\Pi} ^T\mathbb{I}_{ref} ^{-1} \boldsymbol{\Pi} + \dfrac{1}{2M} \mathbf{p}^T \mathbf{p} - \mu \langle \mathbf{B}(\mathbf{x} ), A \mathbf{e}_3  \rangle .
\end{equation}
Let 
 \begin{equation*}
 v _{\left( A, \mathbf{x} \right)} = T_{\left( I, \mathbf{0} \right) } L _{\left( A, \mathbf{x} \right) }\cdot \left( \widehat{ \delta A}, \boldsymbol{\delta } \mathbf{x} \right) = \left( A \widehat{\delta A }, A \boldsymbol{\delta}\mathbf{x} \right) 
 \end{equation*}
 be an arbitrary element of $T_{\left( A, \mathbf{x} \right) } SE(3)$ and $\beta = \left( \boldsymbol{\delta}\boldsymbol{\Pi} , \boldsymbol{\delta}\mathbf{p}  \right) \in \mathfrak{se}(3)^\ast  $. Then,  as
 \begin{multline*}
{\bf d} h \left( (A, \mathbf{x}) , (\boldsymbol{\Pi} , \mathbf{p} )\right) \cdot \left( v_{\left( A, \mathbf{x} \right)  }, \beta\right) =  \left. \dfrac{d}{dt} \right|_0 h \left( \left( (A, \mathbf{x} ) \cdot (\exp t \widehat{ \delta A }, t \boldsymbol{\delta}\mathbf{x}    )\right), (\boldsymbol{\Pi} +t\boldsymbol{\delta}\boldsymbol{\Pi}, \mathbf{p}  + t \boldsymbol{\delta}\mathbf{p} ) \right) = \\
 \langle \mathbb{I}_{ref}^{-1} \boldsymbol{\Pi} ,  \boldsymbol{\delta}\boldsymbol{\Pi}  \rangle + \dfrac{1}{M} \langle \mathbf{p} , \boldsymbol{\delta}\mathbf{p} \rangle  - \mu  \langle DB (\mathbf{x} )^T A \mathbf{e}_3 ,   A \boldsymbol{\delta}\mathbf{x}    \rangle + \langle A (\widehat{ \mathbf{e}_3 \times A^{-1} \mathbf{B}(\mathbf{x} ) }), A \widehat{ \delta A  }  \rangle,
 \end{multline*}
 we can conclude that
 \begin{align}
 D_G h \left( (A, \mathbf{x} ), (\boldsymbol{\Pi} , \mathbf{p} )\right) &=\left( 
\begin{array}{c} 
A \left[ \widehat{ \mathbf{e}_3 \times A^{-1} \mathbf{B} (\mathbf{x} ) }\right]  \\ 
-\mu D\mathbf{B}  (\mathbf{x} )^T A \mathbf{e}_3 
\end{array} \right),\label{DG}\\ 
 D_{g ^\ast } h \left( (A, \mathbf{x} ), (\boldsymbol{\Pi} , \mathbf{p} )\right) &=\left( 
\begin{array}{c} 
\mathbb{I}_{ref}^{-1} \boldsymbol{\Pi}   \\ 
\frac{1}{M} \mathbf{p} 
\end{array} \right).\label{Dgstar} 
\end{align}
Now using \eqref{xG} and \eqref{xgstar}, together with and \eqref{DG}, \eqref{Dgstar},~(\ref{ad star}), and ~(\ref{t star}), we obtain
 \begin{align*}
 X_{g ^\ast } \left( g, \mu \right)  &=\left( 
\begin{array}{c} 
- \mathbf{e}_3 \times A^{-1} \mathbf{B}(\mathbf{x} ) + \boldsymbol{\Pi} \times \mathbb{I}_{ref}^{-1} \boldsymbol{\Pi} \\ 
\mu A ^{-1} D\mathbf{B} (\mathbf{x} )^T A \mathbf{e}_3  + \mathbf{p} \times \mathbb{I}_{ref}^{-1} \boldsymbol{\Pi}
\end{array} \right), \\
 X_G \left( g, \mu \right)   &=\left( 
\begin{array}{c} 
\mathbb{I}_{ref}^{-1} \boldsymbol{\Pi}   \\ 
\frac{1}{M} \mathbf{p} 
\end{array} \right). 
\end{align*}
Consequently, by \eqref{h from partial} we conclude that the equations of motion associated to the Hamiltonian~(\ref{hamiltonian function prop}) are
\begin{align*}
&\dot{A} = A \widehat{ \mathbb{I}_{ref}^{-1} \boldsymbol{\Pi}},\\
&\dot{\mathbf{x} } = \dfrac{1}{M}A {\bf p},\\
&\dot{\boldsymbol{\Pi} } = \boldsymbol{\Pi} \times  \mathbb{I}_{ref}^{-1} \boldsymbol{\Pi} + A^{-1} \mathbf{B}(\mathbf{x} )\times \mathbf{e}_3,  \\
&\dot{\mathbf{p} } = \mathbf{p} \times  \mathbb{I}_{ref}^{-1} \boldsymbol{\Pi}  + \mu A^{-1} D\mathbf{B}( \mathbf{x} ) ^T A \mathbf{e}_3. 
\end{align*}

\subsection{The toral action on phase space $T^*SE(3)$ and the associated momentum map}
\label{The toral action on phase space}

\noindent {\bf The expression of the lifted action in body coordinates.} We start by proving that the cotangent lift of the toral action on $SE (3) $ in~(\ref{toral action config space}) is given  by~(\ref{toral action cotangent}) when using  body coordinates.
Consider $H$ and  $G$  two arbitrary Lie groups and let 
$\Phi : H \times G \longrightarrow G$ be an action of $H$ on $G$. We recall that the lift of this action to the cotangent bundle $T ^\ast G $ of $G$, also denoted by $\Phi$, is given by 
\begin{equation*}
\begin{array}{cccc}
\Phi : &H \times T^*G &\longrightarrow &T^*G\\
	&\left( h, \alpha _g \right) & \longmapsto &T^*_{\Phi_h \left( g\right)  } \Phi _{h^{-1}} \cdot \alpha _g.
\end{array}
\end{equation*}
Using the maps introduced in~(\ref{maps body coordinates}), this action is expressed in body coordinates as:
\begin{equation*}
\Phi \left( h, \left( g, \mu \right) \right) := \lambda  _2 \left( \Phi \left( h, \lambda  _2 ^{-1} \left( g, \mu \right) \right) \right), \quad \mbox{for any} \quad h \in H, g \in G, \text{ and } \mu\in \mathfrak{g}^\ast,
\end{equation*}
or equivalently, 
\begin{equation*}
\Phi _h \left( g, \mu \right) = \left( \Phi _h \left( g \right) , T _e ^* L _{\Phi _h \left( g \right) } \cdot  T^* _{\Phi _h \left( g \right)} \left( L _{g^{-1}} \circ \Phi _{h^{-1}}  \right) \mu \right) = \left( \Phi _h \left( g \right) , T _e ^* \left( L _{g^{-1}} \circ  \Phi _{h^{-1}} \circ L _{\Phi _{h} \left( g\right)}   \right) \mu \right). 
\end{equation*}
In the particular case of $H=\mathbb{T}^2$, $G=SE(3)$, and the toral action introduced in~(\ref{toral action config space}), that is,
\begin{equation*}
\begin{array}{cccc}
\Phi : &(\mathbb{T}^2=S^1\times S^1)\times SE(3)& \longrightarrow &SE(3)\\
	&\left(\left( e^{i\theta_S},e^{i\theta_B}\right) ,(A,\mathbf{x} )\right)&\longmapsto&
({R}^{Z}_{\theta_S} A {R}^{Z}_{-\theta_B},{R}^{Z}_{\theta_S}\mathbf{x} ),
\end{array}
\end{equation*}
we consider $g= \left( A, \mathbf{x}  \right) \in SE(3)$, $\mu = \left( \widehat{\boldsymbol{\Pi}} , \mathbf{p} \right) \in \mathfrak{se}(3)^\ast  $, and $h= \left( e^{i\theta_S},e^{i\theta_B}\right) \in \mathbb{T}^2 $. Then,
\begin{equation}
\label{cotangent lift generic}
\Phi _h \left( g, \mu \right) = \left( \left( R_{\theta _S }^Z A R_{\theta _{-B} }^Z, R _{\theta _S }^Z \mathbf{x} \right), T _e ^* \left( L _{g^{-1}} \circ  \Phi _{h^{-1}} \circ L _{\Phi _{h} \left( g\right)}   \right) \mu \right).
\end{equation}
In order to compute the second part of this expression let $\xi = \left( \widehat{ \boldsymbol{\rho} }, \boldsymbol{\tau} \right) \in \mathfrak{se}(3)$. Then 
\begin{align*}
&\left\langle T _e ^* \left( L _{g^{-1}} \circ  \Phi _{h^{-1}} \circ L _{\Phi _{h} \left( g\right)}   \right) \mu , \xi \right\rangle = \left\langle \mu , T _e  \left( L _{g^{-1}} \circ  \Phi _{h^{-1}} \circ L _{\Phi _{h} \left( g\right)}   \right) \xi \right\rangle \\
&= \left. \frac{d}{dt} \right|_0 \left\langle \left( \widehat{\boldsymbol{\Pi}} , \mathbf{p} \right) , L_{\left( A^{-1}, -A^{-1} \mathbf{x} \right) } \circ \Phi _{\left( e^{-i\theta_S},e^{-i\theta_B}\right) } \circ L_{\left( R_{\theta _S }^Z A R_{-\theta _{B} }^Z, R _{\theta _S }^Z \mathbf{x} \right)} \left( \exp t \widehat{ \boldsymbol{\rho} }, t \boldsymbol{\tau} \right) \right\rangle \\
& =\left. \frac{d}{dt} \right|_0 \left\langle \left( \widehat{\boldsymbol{\Pi}} , \mathbf{p} \right) ,   \left( A^{-1} A R_{-\theta _{B}}^Z  \exp t \widehat{ \boldsymbol{\rho} } R_{\theta _B }^Z, A^{-1} A R_{-\theta _{B}}^Z  t  \boldsymbol{\tau} +A^{-1} \mathbf{x} - A^{-1} \mathbf{x}  \right)  \right\rangle\\
& =\left. \frac{d}{dt} \right|_0 \left\langle \left( \widehat{\boldsymbol{\Pi}} , \mathbf{p} \right) ,   \left( R^Z_{-\theta _{B}}  \exp t \widehat{ \boldsymbol{\rho} } R^Z_{\theta _B }, t R^Z_{-\theta _{B}} \boldsymbol{\tau}  \right)  \right\rangle= \left\langle \left( \widehat{\boldsymbol{\Pi}} , \mathbf{p} \right) ,   \left( \rm{Ad}_{R^Z_{- \theta _B }} \widehat{ \boldsymbol{\rho} }, R^Z_{- \theta _B } \boldsymbol{\tau}    \right)  \right\rangle \\
	&= \left\langle \left( {\rm Ad}^*_{R^Z_{- \theta _B }} \widehat{\boldsymbol{\Pi}}, R^Z_{\theta _B } \mathbf{p} \right), \left( \widehat{ \boldsymbol{\rho} }, \boldsymbol{\tau} \right)   \right\rangle.
\end{align*}
Given that  by~(\ref{elementary relation 3}) ${\rm Ad^*}_{R^Z_{- \theta _B }} \widehat{\boldsymbol{\Pi}} = \widehat{ R^Z_{\theta _B} \boldsymbol{\Pi}} $, the last equality together with~(\ref{cotangent lift generic}) yield the expression~(\ref{toral action cotangent}) of the lifted action in body coordinates, that is,
\begin{equation}
\Phi _{\left( e^{i\theta_S},e^{i\theta_B}\right)} \left( (A, \mathbf{x} ) , ( \boldsymbol{\Pi} ,\mathbf{p} ) \right) = \left( ( R_{\theta _S } A R^Z_{-\theta _B }, R^Z_{\theta _S } \mathbf{x} ) , (  R^Z_{\theta _B } \boldsymbol{\Pi} , R^Z_{\theta _B } \mathbf{p}  ) \right).
\end{equation}

\medskip

\noindent {\bf The infinitesimal generators of the toral action.}
We first show that for any Lie algebra element $\left( \xi _1 , \xi _2 \right) \in \mathbb{R}^2 = {\rm Lie} (\mathbb{T}^2)$ and $\left(A, \mathbf{x} \right) \in SE(3)$,
\begin{align}
(\xi_1,\xi_2)_{SE(3)}\left(A, \mathbf{x} \right)  &=  T_{\left( I, 0 \right) } R_{\left(A, \mathbf{x} \right) }  \left( \widehat{\xi_1 \mathbf{e} _3 -A\xi_2 \mathbf{e} _3 }, A\xi_2 \mathbf{e} _3 \times \mathbf{x} \right)\label{infgenright}\\
	&= T_{\left( I, 0 \right) } L_{ \left( A, \mathbf{x} \right)}  \left( {\rm Ad}_{A^{-1}} \widehat{\xi_1 \mathbf{e} _3} -\widehat{ \xi_2 \mathbf{e} _3 }, A^{-1} \left( \xi_1 \mathbf{e} _3 \times \mathbf{x}\right)  \right).\label{infgenleft1}
\end{align}
We start by proving the first equality
\begin{align*}
&(\xi_1,\xi_2)_{SE(3)}\left(A, \mathbf{x} \right)  = \left. \frac{d}{dt} \right|_0 (\exp t\widehat{\xi_1 \mathbf{e} _3 } A \exp ( - t \widehat{\xi_2 \mathbf{e} _3 }) ,\exp t \widehat{\xi_1 \mathbf{e} _3 } \mathbf{x} )
= (\widehat{\xi_1 \mathbf{e} _3 }A-A\widehat{\xi_2 \mathbf{e} _3 },\widehat{\xi_1 \mathbf{e} _3 } \mathbf{x} ) \\
	&=  (\widehat{\xi_1 \mathbf{e} _3 }A-A\widehat{\xi_2 \mathbf{e} _3 }A^{-1}A,\xi_1 \mathbf{e} _3 \times \mathbf{x} ) 
 =  \left( (\widehat{\xi_1 \mathbf{e} _3 }-A\widehat{\xi_2 \mathbf{e} _3 }A^{-1}) A, (\widehat{\xi_1 \mathbf{e} _3 }-A\widehat{\xi_2 \mathbf{e} _3 }+ A\widehat{\xi_2 \mathbf{e} _3 } ) \times \mathbf{x}  \right) \\
 	&=  \left( (\widehat{\xi_1 \mathbf{e} _3 -A\xi_2 \mathbf{e} _3 })A, (\xi_1 \mathbf{e} _3 -A\xi_2 \mathbf{e} _3) \times \mathbf{x} + (A\xi_2 \mathbf{e} _3 \times \mathbf{x})  \right) = T_{\left( I, 0 \right) } R_{\left(A, \mathbf{x} \right) }  \left( (\widehat{\xi_1 \mathbf{e} _3 -A\xi_2 \mathbf{e} _3 }), A\xi_2 \mathbf{e} _3 \times \mathbf{x} \right),
\end{align*}
where in the last equality we used \eqref{rightact}. Regarding~(\ref{infgenleft1}), note that
\begin{align*}
&(\xi_1,\xi_2)_{SE(3)}\left(A, \mathbf{x} \right)  = \left. \frac{d}{dt} \right|_0 (\exp t\widehat{\xi_1 \mathbf{e} _3 } A \exp ( - t \widehat{\xi_2 \mathbf{e} _3 }) ,\exp t \widehat{\xi_1 \mathbf{e} _3 } \mathbf{x} )
= (\widehat{\xi_1 \mathbf{e} _3 }A-A\widehat{\xi_2 \mathbf{e} _3 },\widehat{\xi_1 \mathbf{e} _3 } \mathbf{x} ) \\
	&=  (A A^{-1} \widehat{\xi_1 \mathbf{e} _3 }A-A\widehat{\xi_2 \mathbf{e} _3 },\xi_1 \mathbf{e} _3 \times \mathbf{x} ) =\left( T_I L_A({\rm Ad}_{A^{-1}} \widehat{\xi_1 \mathbf{e} _3 }-\widehat{\xi_2 \mathbf{e} _3 }), (AA^{-1}(\xi_1 \mathbf{e} _3 \times \mathbf{x}))  \right)\\
	& = T_{\left( I, 0 \right) } L_{ \left( A, \mathbf{x} \right)}  \left( {\rm Ad}_{A^{-1}} \widehat{\xi_1 \mathbf{e} _3} -\widehat{ \xi_2 \mathbf{e} _3 }, A^{-1} \left( \xi_1 \mathbf{e} _3 \times \mathbf{x}\right)  \right),
\end{align*}
where we used \eqref{leftact}.

The infinitesimal generator of the lifted $\mathbb{T}^2$-action on $T^*SE(3)$ in body coordinates is given by
\begin{equation}
\label{infinitgenBody}
(\xi_1,\xi_2)_{T^*SE(3)}\left(A, \mathbf{x}, \boldsymbol{\Pi}  , \mathbf{p}   \right)= \left( A ({\rm Ad}_{A^{-1}} (\widehat{\xi_1 \mathbf{e} _3} )-\widehat{ \xi_2 \mathbf{e} _3 }), \widehat{  {\xi_1 \mathbf{e} _3 } } \mathbf{x} , \widehat{ {\xi_2 \mathbf{e} _3 }}  \boldsymbol{\Pi}  , \widehat{ {\xi_2 \mathbf{e} _3 }} \mathbf{p}  \right)
\end{equation}
Indeed,
\begin{align*}
&(\xi_1,\xi_2)_{T^*SE(3)}\left(A, \mathbf{x}, \boldsymbol{\Pi}  , \mathbf{p}   \right)  = \left. \frac{d}{dt} \right|_0  \exp t (\xi_1,\xi_2) \cdot \left(A, \mathbf{x}, \boldsymbol{\Pi}  , \mathbf{p}   \right)  \\
&= \left. \frac{d}{dt} \right|_0 \left( \exp t\widehat{\xi_1 \mathbf{e} _3 } A \exp (-t\widehat{\xi_2 \mathbf{e} _3 }), \exp t\widehat{\xi_1 \mathbf{e} _3 } \mathbf{x} , \exp t\widehat{\xi_2 \mathbf{e} _3 } \boldsymbol{\Pi} , \exp t\widehat{\xi_2 \mathbf{e} _3 } \mathbf{p}  \right) \nonumber \\
&= \left( AA ^{-1}\widehat{\xi_1 \mathbf{e} _3}A -A\widehat{ \xi_2 \mathbf{e} _3 },  \widehat{\xi_1 \mathbf{e} _3 } \mathbf{x} , \widehat{\xi_2 \mathbf{e} _3 } \boldsymbol{\Pi}  , \widehat{\xi_2 \mathbf{e} _3 } \mathbf{p}  \right) = \left( A \left({\rm Ad}_{A^{-1}} \left(\widehat{\xi_1 \mathbf{e} _3}\right) -\widehat{ \xi_2 \mathbf{e} _3 }\right), \widehat{  {\xi_1 \mathbf{e} _3 } } \mathbf{x} , \widehat{ {\xi_2 \mathbf{e} _3 }}  \boldsymbol{\Pi}  , \widehat{ {\xi_2 \mathbf{e} _3 }} \mathbf{p}  \right).
\end{align*}

\medskip

\noindent {\bf The momentum map of the toral action}
Given a lifted  action of a Lie group $H$ on the cotangent bundle $T ^\ast G $ of a Lie group $G$ endowed with the canonical symplectic form, the map $\mathbf{J}: T ^\ast G\longrightarrow \mathfrak{g}^\ast $ defined by
\begin{equation}
\label{canonical momentum map}
\langle \boldsymbol{J} (\alpha_g),\xi\rangle = \langle\alpha_g,\xi_G(g)\rangle \quad \mbox{for any} \quad g \in G,\, \alpha _g \in T^*G, \text{ and } \xi \in \mathfrak{h},
\end{equation}
is a coadjoint equivariant momentum map for this canonical action (see~\cite[Corollary 4.2.11]{Abraham1978}).
We now study the  particular case  we are interested in, that is, $H = \mathbb{T}^2$, $G=SE(3)$, and consider  an arbitrary point $g=\left(A, \mathbf{x} \right) \in SE (3)$, $\mu= \left(\boldsymbol{\Pi},\mathbf{p}\right) \in \mathfrak{se}(3) ^\ast  $ and $\alpha _g = T^* _g L _{g^{-1}} \cdot \mu \in T ^\ast SE(3)$ the covector that in body coordinates is expressed as $\left(\left(A, \mathbf{x} \right),\left(\boldsymbol{\Pi},\mathbf{p}\right)\right)$. With this notation, the expression in body coordinates of the momentum map $\mathbf{J}: SE (3) \times \mathfrak{se}(3)^\ast \longrightarrow \mathbb{R}^2 $ in~(\ref{canonical momentum map}) is given by
\begin{equation}
\label{momentum map in body coordinates}
\boldsymbol{J} \left( \left(A, \mathbf{x} \right) , \left(\boldsymbol{\Pi} ,\mathbf{p} \right) \right) = \left( \langle A \boldsymbol{\Pi}  + \mathbf{x} \times A \mathbf{p}  ,   \mathbf{e}_3 \rangle , -\langle \boldsymbol{\Pi}  ,  \mathbf{e}_3 \rangle \right). 
\end{equation}
Indeed, for any $(\xi_1, \xi_2)\in \mathbb{R}^2 $, 
\begin{align*}
&\left\langle\boldsymbol{J}  \left( \left(A, \mathbf{x} \right) , \left(\boldsymbol{\Pi} ,\mathbf{p} \right) \right) , \left( \xi _1 , \xi _2 \right) \right\rangle = 
\left\langle  T^*_{\left(A, \mathbf{x} \right) } L_{\left(A, \mathbf{x} \right) ^{-1} }  \left(\boldsymbol{\Pi} ,\mathbf{p} \right) , T_{\left( I, 0 \right) } L_{\left(A, \mathbf{x} \right) }  \left( (\widehat{A^{-1} \xi_1 \mathbf{e} _3 -\xi_2 \mathbf{e} _3 }), A^{-1} (\xi_1 \mathbf{e} _3 \times \mathbf{x} )\right)  \right\rangle \\
&= \langle \boldsymbol{\Pi}  , A^{-1}\xi_1 \mathbf{e} _3 -\xi_2 \mathbf{e} _3 \rangle + \langle \mathbf{p}  , A^{-1}(\xi_1 \mathbf{e} _3 \times \mathbf{x}) \rangle =  \langle \boldsymbol{\Pi}  , A^{-1}\xi_1 \mathbf{e}_3  -\xi_2  \mathbf{e}_3  \rangle + \langle \mathbf{p}  , A^{-1}(\xi_1 \mathbf{e}_3  \times \mathbf{x}) \rangle,
\end{align*}
which proves~(\ref{momentum map in body coordinates}) since $(\xi_1, \xi_2)\in \mathbb{R}^2 $ is arbitrary.

\subsection{Proof of Proposition~\ref{prop_rel_eq}}
\label{proof proposition relative}

\medskip

\noindent {\bf (i)}
Using the statement preceeding~(\ref{relative critical points}) we will specify the relative equilibria of the orbitron by characterizing the points $\mathbf{z}=\left(  \left(A, \mathbf{x} \right) , \left(\boldsymbol{\Pi},\mathbf{p}\right)\right) \in T^*SE(3)$ for which
\begin{equation}
 \mathbf{d} (h - \boldsymbol{J} ^{(\xi_1,\xi_2)}) \left( \left(A, \mathbf{x} \right) , \left(\boldsymbol{\Pi},\mathbf{p}\right) \right) =0
\end{equation}
for some $(\xi_1,\xi_2) \in \mathbb{R}^2$. We start by computing the tangent of the momentum map and the differential of the Hamiltonian. Let $\mathbf{v}=\left( ( \widehat{\delta A  } A , \boldsymbol{\delta}\mathbf{x} ) , (  \boldsymbol{\delta}\boldsymbol{\Pi} , \boldsymbol{\delta}\mathbf{p} )\right) \in T_{\mathbf{z}} \left( T^*SE(3) \right) $ be an arbitrary vector at the point $\mathbf{z}$, then it is easy to check that 
\begin{align}
 \mathbf{d} T \left(\boldsymbol{\Pi},\mathbf{p}\right) \cdot \mathbf{v}  &= \langle \boldsymbol{\Pi}, \mathbb{I}_{ref}^{-1} \boldsymbol{\delta}\boldsymbol{\Pi} \rangle + \frac{1}{M} \langle \mathbf{p}, \boldsymbol{\delta}\mathbf{p} \rangle,\label{diff t}\\
\mathbf{d} V \left(A, \mathbf{x} \right) \cdot \mathbf{v} &=  -\mu \left[ \langle D\mathbf{B}(\mathbf{x})(\boldsymbol{\delta}\mathbf{x}), A   \mathbf{e}_3  \rangle + \langle \mathbf{B} (\mathbf{x} ), {\delta A  } \times A  \mathbf{e}_3 \rangle \right].\label{diff v}
 \end{align}
with $T$ and $V$ the kinetic and potential energies introduced in~(\ref{potential and kinetic}). Additionally,
\begin{align}
&T_{ \left( \left(A, \mathbf{x} \right) , \left(\boldsymbol{\Pi} ,\mathbf{p} \right) \right)} \boldsymbol{J}  \cdot \left( (\widehat{\delta A  } A , \boldsymbol{\delta}\mathbf{x} ) , (\boldsymbol{\delta}\boldsymbol{\Pi}  , \boldsymbol{\delta}\mathbf{p} ) \right) \nonumber \\
&=\Big( \langle \delta A \times A\boldsymbol{\Pi}  +A \boldsymbol{\delta \Pi}  + \boldsymbol{\delta} \mathbf{x} \times A\mathbf{p}  + \mathbf{x} \times ( \delta A \times A \mathbf{p}  ) + \mathbf{x}  \times A \boldsymbol{\delta }\mathbf{p} , \mathbf{e}_3 \rangle,- \langle \boldsymbol{\delta}\boldsymbol{\Pi}  , \mathbf{e}_3 \rangle \Big).\label{differential momentum map}
\end{align}
Consequently, using~(\ref{diff t}),~(\ref{diff v}) and \eqref{differential momentum map} we have, for any $(\xi_1, \xi_2)\in \mathbb{R}^2$
\begin{multline}
\label{critical 1}
\mathbf{d} \left(h - \boldsymbol{J} ^{(\xi_1,\xi_2)}\right) (\mathbf{z}) \cdot \mathbf{v}  = \boldsymbol{\Pi} ^T\mathbb{I}_{ref}^{-1} \boldsymbol{\delta}\boldsymbol{\Pi} + \frac{1}{M} \mathbf{p} \cdot \boldsymbol{\delta}\mathbf{p} -\mu \left[ \langle D\mathbf{B}(\mathbf{x})(\boldsymbol{\delta}\mathbf{x}), A   \mathbf{e}_3  \rangle + \langle \mathbf{B} (\mathbf{x} ), {\delta A  } \times A  \mathbf{e}_3 \rangle \right]  \\
+ \xi _2 \boldsymbol{\delta}\boldsymbol{\Pi}  \cdot  \mathbf{e}_3 - \xi _1  \left( \delta A \times A\boldsymbol{\Pi}+A \boldsymbol{\delta}\boldsymbol{\Pi}  + \boldsymbol{\delta}\mathbf{x} \times A\mathbf{p} +\mathbf{x} \times (\delta A \times A \mathbf{p} ) + \mathbf{x} \times A\boldsymbol{\delta}\mathbf{p}  \right) \cdot \mathbf{e}_3.   
\end{multline}
Therefore, as $\widehat{\delta A  }$, $ \boldsymbol{\delta}\mathbf{x}$, $\boldsymbol{\delta}\boldsymbol{\Pi}$, and $\boldsymbol{\delta}\mathbf{p}$ in this expression are arbitrary, it can be checked that the points $\mathbf{z} \in T^*SE(3)$ for which 
$
 \mathbf{d} (h - \boldsymbol{J} ^{(\xi_1,\xi_2)}) \left(\mathbf{z} \right) =0
$
are characterized by the equations:
\begin{align}
\label{equationderivative1}
&\mu \left[  \mathbf{B} (\mathbf{x} ) \times A  \mathbf{e}_3 \right]  + \xi _1  \left[ A \mathbf{p} \times (\mathbf{x} \times \mathbf{e}_3 ) -  A \boldsymbol{\Pi} \times \mathbf{e}_3 \right] = 0, \\
\label{equationderivative2}
&- \mu   D\mathbf{B} (\mathbf{x} ) ^T (A  \mathbf{e}_3)  - \xi _1 \left( A \mathbf{p} \times \mathbf{e}_3 \right)= 0,\\
\label{equationderivative3}
&\mathbb{I}_{ref}^{-1} \boldsymbol{\Pi} + \xi _2 \mathbf{e}_3 - \xi _1 A^{-1} \mathbf{e}_3 =0,\\
\label{equationderivative4}
& \frac{1}{M} \mathbf{p} - \xi _1 A^{-1} \left(  \mathbf{e}_3 \times \mathbf{x} \right) = 0,
\end{align}
as required. 

\medskip

\noindent {\bf (ii)}
We show that the points $\mathbf{z} _0 = \left( \left(A_0, \mathbf{x}_0 \right) , \left(\boldsymbol{\Pi}_0,\mathbf{p}_0\right) \right) $ of the form specified in the statement of the proposition satisfy equations \eqref{equationderivative1}--\eqref{equationderivative4} and hence constitute a branch of relative equilibria. We proceed by considering $A_0=R_{\theta_0}^Z$ and $\mathbf{x} _0= \left( x, y, 0\right) $ and using equations \eqref{equationderivative1}--\eqref{equationderivative4} to determine $\boldsymbol{\Pi}_0 $, $\mathbf{p}_0 $, and the velocity $ \boldsymbol{\xi } = \left( \xi _1 , \xi _2 \right) $ in the statement.

Notice first that $A_0 \mathbf{e}_3 =  \mathbf{e}_3 $, hence by \eqref{equationderivative4} we have that 
\begin{equation} \label{p0regular}
\mathbf{p} _0=M \xi _1 A_0 ^{-1}\left( -y,x, 0\right),
\end{equation}
necessarily. Now by \eqref{equationderivative3} \\
\begin{equation} 
\boldsymbol{\Pi} _0 =  \mathbb{I}_{ref}\left( \xi _1 - \xi _2 \right) \mathbf{e}_3 = I_3 \left( \xi _1 - \xi _2 \right) \mathbf{e}_3.
\end{equation}
In order to handle \eqref{equationderivative2} we note that $D\mathbf{B}(\mathbf{x})$ is given by the matrix whose components are
 \begin{align*}
& \dfrac{\partial{B_x}}{\partial{x}} =k \left(  \dfrac{D(\mathbf{x} )_+ - 3x^2}{D(\mathbf{x} )_+ ^{5/2}} - \dfrac{D(\mathbf{x} )_- - 3x^2}{D( \mathbf{x} )_- ^{5/2}} \right), \\
&\dfrac{\partial{B_x}}{\partial{y}} =k \left(  \dfrac{- 3xy}{D(\mathbf{x} )_+ ^{5/2}} + \dfrac{3xy}{D(\mathbf{x} )_- ^{5/2}} \right), \\
&\dfrac{\partial{B_x}}{\partial{z}} = k \left( \dfrac{- 3x(z-h)}{D(\mathbf{x} )_+ ^{5/2}} + \dfrac{3x(z+h)}{D(\mathbf{x} )_- ^{5/2}} \right), \\
&\dfrac{\partial{B_y}}{\partial{x}} =k \left(  \dfrac{- 3xy}{D(\mathbf{x} )_+ ^{5/2}} + \dfrac{3xy}{D( \mathbf{x} )_-^{5/2}}\right), \\
&\dfrac{\partial{B_y}}{\partial{y}} = k \left( \dfrac{D(\mathbf{x} )_+ - 3y^2}{D(\mathbf{x} )_+ ^{5/2}} - \dfrac{D( \mathbf{x} )_- - 3y^2}{D( \mathbf{x} )_- ^{5/2}}\right), \\
&\dfrac{\partial{B_y}}{\partial{z}} = k \left( \dfrac{- 3y(z-h)}{D(\mathbf{x} )_+^{5/2}} + \dfrac{3y(z+h)}{D( \mathbf{x} )_- ^{5/2}}\right),
 \end{align*}
  \begin{align*}
&\dfrac{\partial{B_z}}{\partial{x}} = k \left( \dfrac{- 3x(z-h)}{D(\mathbf{x} )_+^{5/2}} + \dfrac{3x(z+h)}{D( \mathbf{x} )_- ^{5/2}}\right), \\
&\dfrac{\partial{B_z}}{\partial{y}} = k \left( \dfrac{- 3y(z-h)}{D(\mathbf{x} )_+ ^{5/2}} + \dfrac{3y(z+h)}{D( \mathbf{x} )_- ^{5/2}}\right), \\
& \dfrac{\partial{B_z}}{\partial{z}} = k \left( \dfrac{D(\mathbf{x} )_+ - 3(z-h)^2}{D(\mathbf{x} )_+^{5/2}} - \dfrac{D( \mathbf{x} )_- - 3(z+h)^2}{D( \mathbf{x} )_-^{5/2}}\right), 
 \end{align*}
where $D(\mathbf{x} )_+=x ^2 + y ^2 + (z-h) ^2$, $D(\mathbf{x} )_-=x ^2 + y ^2 + (z+h) ^2$ and $k = \dfrac{\mu _0 q}{4 \pi }$.
Consequently,
\begin{align*}
D\mathbf{B}({\bf x}_0)=k \left( 
\begin{array}{ccc} 
0 & 0 & \frac{6xh}{D(\mathbf{x} _0)^{5/2}} \\ 
0 & 0 & \frac{6yh}{D(\mathbf{x} _0)^{5/2}}\\ 
\frac{6xh}{D(\mathbf{x} _0)^{5/2}} & \frac{6yh}{D(\mathbf{x} _0)^{5/2}} & 0
\end{array} 
\right),
\end{align*}
where $D(\mathbf{x} _0)=D(\mathbf{x} _0)_+ = D (\mathbf{x} _0)_-$.
Hence 
\begin{equation}
\label{DBx0regular}
D\mathbf{B}\left( \mathbf{x}_0 \right)^T \left( A_0 \mathbf{e}_3 \right)= D\mathbf{B}\left( \mathbf{x}_0 \right)^T \mathbf{e}_3 = \dfrac{6kh}{D(\mathbf{x}_0)^{5/2}} \mathbf{x}_0.
\end{equation}
Note additionally that by \eqref{p0regular} 
\begin{equation}
\label{p0timesAe3}
A_0\mathbf{p}_0 \times \mathbf{e}_3 =  M \xi _1 \mathbf{x}_0. 
\end{equation}
Then by equalities \eqref{DBx0regular} and by \eqref{p0timesAe3}, equation \eqref{equationderivative2} holds whenever $\mathbf{x}_0= \boldsymbol{0} $ or when $\mathbf{x}_0 \neq \boldsymbol{0} $ and $\xi _1 ^2 =  - \dfrac{3h\mu q \mu _0 }{2 \pi MD(\mathbf{x}_0 )^{5/2}} $; we note that in both situations, there are no restrictions on the second component of the velocity $\xi_2$. Finally, it can be readily verified that \eqref{equationderivative1} always holds at the point $\left( \left(A_0, \mathbf{x}_0 \right) , \left(\boldsymbol{\Pi}_0,\mathbf{p}_0\right) \right) $ by using that $\mathbf{B}(\mathbf{x} _0)=-\dfrac{\mu_0qh}{2\pi D(\mathbf{x} _0)^{3/2}} \mathbf{e}_3$.

\medskip

\noindent {\bf (iii)}
Suppose that we are in the presence of a magnetic field $\mathbf{B}$  equivariant with respect to rotations around the $OZ$ axis and that behaves as indicated in  \eqref{Bx_transformation}--\eqref{Bz_transformation} with respect to the mirror transformation~(\ref{mirror transformation}). Notice first that by \eqref{Bx_transformation} and \eqref{By_transformation} 
\begin{equation}\label{Bxy_zero}
B _x (x,y,0) = B _y (x,y,0)
\end{equation}
and hence 
\begin{equation}\label{Bz_e3}
\mathbf{B}(x,y,0)=B _z (x,y,0) \mathbf{e}_3. 
\end{equation}
Additionally, by \eqref{mag_equivariant}, $B_z (x,y,0)$ is rotationally invariant with respect to rotations in the $OXY$ plane, hence
\begin{equation}\label{Bz_func_sqr_xy}
B _z (x,y,0) = f(x ^2 + y ^2 ), \enspace {\rm for \enspace some}  \enspace f \in C^\infty(\mathbb{R} ^2 ).
\end{equation}
Conditions \eqref{equationderivative4s} and \eqref{equationderivative3s} show that if $A _0 = R_{\theta _0 } ^z $ and $\mathbf{x} _0 = (x,y,0)$, then $\boldsymbol{\Pi} _0 = I _3 \left( \xi _1 - \xi _2 \right) \mathbf{e}_3 $ and $\mathbf{p} _0 =  M \xi _1 A _0 ^{-1} \left(-y,x,0\right) $ necessarily. If we use $\mathbf{z} _0 = \left( (A _0 , \mathbf{x} _0 ),(\boldsymbol{\Pi} _0 , \mathbf{p} _0 )\right)  $ and  \eqref{Bz_e3} in the expression \eqref{equationderivative1s}, it can be easily verified that this relation is automatically satisfied.

In order to study the expression \eqref{equationderivative2s}, we take derivatives on both sides of \eqref{Bz_transformation} and obtain that 
\begin{equation*}
\partial _z B _z (x,y,z) = - \partial _z B _z (x,y,-z)
\end{equation*}
which shows that 
\begin{equation}\label{deriv_Bz_zero}
\partial _z B _z (x,y,0) = 0.
\end{equation}
Finally, by \eqref{Bz_e3} and \eqref{deriv_Bz_zero} the relation \eqref{equationderivative2s} amounts to 
\begin{equation*}
- \mu (\partial _x B _z , \partial _y B _z, 0 )=M \xi _1 ^2 \mathbf{x} _0.
\end{equation*}
By \eqref{Bz_func_sqr_xy} this is equivalent to 
\begin{equation*}
-2 \mu f'(x ^2 + y ^2 ) \mathbf{x} _0 = M \xi _1 ^2 \mathbf{x} _0,
\end{equation*}
which guarantees that \eqref{equationderivative2s} is satisfied provided that
\begin{equation}
\xi _1 = \pm \left( -\dfrac{2}{M} \mu f'( x ^2 + y ^2 )\right) ^{1/2},
\end{equation}
as required. \quad $\blacksquare$

\subsection{Proof of Theorem~\ref{kozorez relations}}
\label{proof kozorez}

We will proceed by  using Theorem~\ref{energy momentum statement} in order to determine the regions in parameter space  for which the stability form~(\ref{quadratic form stability}) at the relative equilibria  is definite, which in turn ensures $\mathbb{T}^2 $--stability. 

We start by denoting the augmented Hamiltonian as $h ^{\boldsymbol{\xi}} :=h- \boldsymbol{J}^{\boldsymbol{\xi}}   $, for any $\boldsymbol{\xi}= \left( \xi _1 , \xi _2 \right)\in {\rm Lie} \left(\mathbb{T}^2\right)$. Let  $\mathbf{z}= \left( \left(A, \mathbf{x} \right) , \left(\boldsymbol{\Pi},\mathbf{p}\right) \right) \in T ^\ast \left(SE (3)\right)$ expressed in body coordinates. As we saw in the proof of Proposition~\ref{prop_rel_eq} (see Appendix~\ref{proof proposition relative}), the partial derivatives of $h ^{\boldsymbol{\xi}} $ are given by:
\begin{itemize}
\item
$ h_A^ {\boldsymbol{\xi}}:=D_A h ^{\boldsymbol{\xi}} (\mathbf{z}) = \mu \left[  \mathbf{B} (\mathbf{x} ) \times A  \mathbf{e}_3 \right]  + \xi _1  \left[ A \mathbf{p} \times (\mathbf{x} \times \mathbf{e}_3 ) -  A \boldsymbol{\Pi} \times \mathbf{e}_3 \right]  $,
\item 
$h_{\mathbf{x}} ^{\boldsymbol{\xi}}:=D_{\mathbf{x} } h ^ {\boldsymbol{\xi}} (\mathbf{z}) = - \mu   D\mathbf{B} (\mathbf{x} ) ^T (A  \mathbf{e}_3)  - \xi _1 \left( A \mathbf{p} \times \mathbf{e}_3 \right) $,
\item
$h _{\boldsymbol{\Pi} } ^{\boldsymbol{\xi}}:=D_{\boldsymbol{\Pi} } h ^{\boldsymbol{\xi}} (\mathbf{z}) = \mathbb{I}_{ref}^{-1} \boldsymbol{\Pi} + \xi _2 \mathbf{e}_3 - \xi _1 A^{-1} \mathbf{e}_3$,
\item
$h_{\mathbf{p} }^ {\boldsymbol{\xi}} :=D_{\mathbf{p} } h ^ {\boldsymbol{\xi}} (\mathbf{z}) = \frac{1}{M} \mathbf{p} - \xi _1 A^{-1} \left(  \mathbf{e}_3 \times \mathbf{x} \right)$.
\end{itemize}
In order to compute the Hessian of the augmented Hamiltonian, we write down the derivatives of its partial derivatives in the direction given by the vector
$\mathbf{v}=\left. \dfrac{d}{dt} \right|_0 \left( ( \exp t \widehat{\delta A  } A , \mathbf{x} + t \boldsymbol{\delta}\mathbf{x} ) , ( \boldsymbol{\Pi} + t \boldsymbol{\delta}\boldsymbol{\Pi} , \mathbf{p} + t \boldsymbol{\delta}\mathbf{p} )\right)$. A straightforward computation yields:
\begin{itemize}
\item $\mathbf{d} h_A ^ {\boldsymbol{\xi}} ( \mathbf{z}) \cdot\mathbf{v} =  \mu \left[  (D\mathbf{B}(\mathbf{x} ) \boldsymbol{\delta}\mathbf{x})  \times A \mathbf{e}_3 +  \mathbf{B} (\mathbf{x} ) \times ( \widehat{\delta A  } A  \mathbf{e}_3 )\right]  + \xi _1  \Big[ ( \widehat{\delta A  } A \mathbf{p} + A  \boldsymbol{\delta}\mathbf{p})\times ( \mathbf{x} \times \mathbf{e}_3 ) + A \mathbf{p} \times (\boldsymbol{\delta}\mathbf{x}  \times \mathbf{e}_3 ) - \widehat{\delta A  } A \boldsymbol{\Pi} \times \mathbf{e}_3  - ( A \boldsymbol{\delta}\boldsymbol{\Pi} \times \mathbf{e}_3 ) \Big]$,
\item $\mathbf{d} h_{\mathbf{x} } ^{\boldsymbol{\xi}} ( \mathbf{z} ) \cdot \mathbf{v} =  -\mu \left( T_{\mathbf{x} } \mathbf{F} ( \boldsymbol{\delta}\mathbf{x} )\right) (A \mathbf{e}_3 )-\mu \mathbf{F} (\mathbf{x} ) \left( \delta A \times A \mathbf{e}_3\right) - \xi _1 ( \widehat{ \delta A } A \mathbf{p} \times \mathbf{e}_3 + A \boldsymbol{\delta}\mathbf{p} \times \mathbf{e}_3 )$,
where 
\begin{align*}
\mathbf{F} :\ & \mathbb{R}^3 \longrightarrow M_{3 \times 3}\\
& \mathbf{x} \longmapsto D\mathbf{B}(\mathbf{x} )^T, 
\end{align*}
\item $\mathbf{d} h_{\boldsymbol{\Pi} } ^ {\boldsymbol{\xi}} (\mathbf{z}) \cdot \mathbf{v} =  \mathbb{I}_{ref}^{-1} \boldsymbol{\delta}\boldsymbol{\Pi} + \xi _1 A^T \widehat{ \delta A}  \mathbf{e}_3$,
\item $\mathbf{d} h_{\mathbf{p} } ^ {\boldsymbol{\xi}} ( \mathbf{z}) \cdot \mathbf{v} = \dfrac{\boldsymbol{\delta}\mathbf{p} }{M} - \xi _1 A^T\widehat{ \delta A } ( \mathbf{x} \times \mathbf{e}_3)+ \xi _1 A^T (\boldsymbol{\delta}\mathbf{x} \times \mathbf{e}_3 ) $.
\end{itemize}
Consequently, the matrix expression associated to $\mathbf{d}^2 \left(h- \mathbf{J}^{\boldsymbol{\xi}}\right) (\mathbf{z}) $ is given by:
\begin{equation}
\label{hessian in general}
\left(
\begin{array}{cccc}
- \mu \left[ \widehat{ \mathbf{B} (\mathbf{x} ) } \widehat{ A  \mathbf{e}_3 }\right]  + \xi _1  \Big[ \widehat{ \mathbf{x} \times \mathbf{e}_3 }  \widehat{ A \mathbf{p}} - \widehat{ \mathbf{e}_3 } \widehat{ A \boldsymbol{\Pi}  }\Big] &
-\mu   \widehat{ A \mathbf{e}_3 } \mathbf{F}(\mathbf{x} )^T  - \xi _1  \widehat{ A \mathbf{p} } \widehat{ \mathbf{e}_3 } &
\xi _1 \widehat{\mathbf{e}_3 }A&
-\xi _1 \widehat{  \mathbf{x} \times  \mathbf{e}_3 } A\\
\mu \left[ \mathbf{F} (\mathbf{x} ) \widehat{ A \mathbf{e}_3 }\right] - \xi _1 \widehat{ \mathbf{e}_3 } \widehat{ A \mathbf{p} } &
-\mu  {T} _{\mathbf{x} } \mathbf{F} (\cdot) (A \mathbf{e}_3 )&
0&
\xi _1 \widehat{\mathbf{e}_3 }A\\
- \xi _1 A^T\widehat{ \mathbf{e}_3 } & 0 & \mathbb{I}_{ref}^{-1} &
0\\
\xi _1 A^T \widehat{ \mathbf{x} \times \mathbf{e}_3 } &
- \xi _1 A^T \widehat{\mathbf{e}_3 } &
0&
\dfrac{1}{M} \mathbb{I}_{id}
\end{array}
\right)
\end{equation}
We now compute the value of the Hessian~(\ref{hessian in general}) at the relative equilibria in the second and third parts of Proposition~\ref{prop_rel_eq}, that is, $\mathbf{z}_0= \left( \left(A_0, \mathbf{x}_0 \right) , \left(\boldsymbol{\Pi}_0,\mathbf{p}_0\right) \right) $ with $A_0=R_{\theta}^Z$, $\mathbf{x}_0= \left( x, y, 0\right)^T$, $\boldsymbol{\Pi}_0=I_3 \left( \xi _1^0 - \xi _2\right) \mathbf{e}_3 $, and $\mathbf{p}_0=M \xi _1 ^0A _0^{-1} \left( -y, x, 0\right)^T  $. We start by noticing that 
\begin{equation}\label{TxF_is_Hess}
T _{\mathbf{x} } \mathbf{F} (\cdot)(\mathbf{e}_3 ) = {\rm  Hess} \left( {B_z} \right)(\mathbf{x} ).
\end{equation}
Indeed, for any $\boldsymbol{\delta}\mathbf{x} \in T_{\mathbf{x}}\mathbb{R}^3 $
\begin{multline}\label{TxF_is_Hess_proof}
T _{\mathbf{x} }\mathbf{F} (\boldsymbol{\delta}\mathbf{x} )(\mathbf{e}_3 )= \left. \dfrac{d}{dt} \right|_0 \mathbf{F} \left( \mathbf{x} + t \boldsymbol{\delta}\mathbf{x} \right) (\mathbf{e}_3 ) = \left. \dfrac{d}{dt} \right|_0 D\mathbf{B} \left( \mathbf{x} + t \boldsymbol{\delta}\mathbf{x} \right)^T \mathbf{e}_3 \\
= \left. \dfrac{d}{dt} \right|_0 \left( \begin{array}{ccc} \dfrac{\partial{B_z}}{\partial{x}} \left( \mathbf{x} + t \boldsymbol{\delta}\mathbf{x} \right)\\ \\\dfrac{\partial{B_y}}{\partial{z}}\left( \mathbf{x} + t \boldsymbol{\delta}\mathbf{x} \right)\\\\ \dfrac{\partial{B_z}}{\partial{z}} \left( \mathbf{x} + t \boldsymbol{\delta}\mathbf{x} \right) \end{array} \right) = {\rm Hess}  \left( {B_z} \right)(\mathbf{x} ) \cdot \boldsymbol{\delta}\mathbf{x}. 
\end{multline}
Therefore, the matrix expression associated to $\mathbf{d}^2 \left(h- \mathbf{J}^{\boldsymbol{\xi}}\right) (\mathbf{z}_0) $ is given by:
\begin{equation}
\label{hessian in particular}
\left(
\begin{array}{cccc} 
- \mu \left[ \widehat{ \mathbf{B} (\mathbf{x}_0 ) } \widehat{\mathbf{e}_3 }\right]  + \xi _1  \Big[   \widehat{ \mathbf{x} _0 \times \mathbf{e}_3  } \widehat{ \mathbf{p} _0 } - \widehat{\mathbf{e}_3} \widehat{ \boldsymbol{\Pi} _0 }\Big]
& -\mu\widehat{\mathbf{e}_3 } \mathbf{F}(\mathbf{x}_0 )^T   - \xi _1  \widehat{ \mathbf{p}_0 } \widehat{ \mathbf{e}_3}
& \xi _1 \widehat{ \mathbf{e}_3 }
& -\xi _1 (\widehat{  \mathbf{x}_0 \times  \mathbf{e}_3})\\ 
\mu\mathbf{F} (\mathbf{x}_0 ) \widehat{ \mathbf{e}_3 }- \xi _1 \widehat{  \mathbf{e}_3} \widehat{  \mathbf{p}_0}  
&-\mu {\rm Hess}  \left( {B_z} \right)(\mathbf{x}_0 ) 
&0
&\xi _1 \widehat{\mathbf{e}_3 }\\
- \xi _1 \widehat{\mathbf{e}_3 }
&0
&\mathbb{I}_{ref}^{-1}
&0\\
\xi _1( \widehat{ \mathbf{x} _0 \times \mathbf{e}_3 }) 
&- \xi _1 \widehat{ \mathbf{e}_3 }
&0
&\dfrac{1}{M} \mathbb{I}_{id}
\end{array} \right).
\end{equation}
In order to construct stability forms for the regular and singular branches, we  now determine stability spaces $W$ to which we will restrict the Hessian~(\ref{hessian in particular}). 

\medskip

\noindent {\bf A stability space for the regular branch ($r>0 $).} In this case, the kernel of the derivative of the momentum map is given by:
\begin{align*}
\ker T_{{\bf z}_0} \boldsymbol{J}& = \left\{ v= \left( ( \widehat{\delta A   } , \boldsymbol{\delta}\mathbf{x} ) , \left( \boldsymbol{\delta}\boldsymbol{\Pi} , \boldsymbol{\delta}\mathbf{p} \right) \right) \in T_{\mathbf{z}_0} \left( SE(3) \times \mathfrak{se}(3)^*\right) \mid T_{\mathbf{z}_0 } \boldsymbol{J} \cdot v = 0 \right\}\\
&= \left\{ v \in T_{\mathbf{z}_0 } \left( SE(3) \times \mathfrak{se}(3)^*\right) \mid \boldsymbol{\delta}\boldsymbol{\Pi} \cdot \mathbf{e}_3 =0, \delta p _2 =- M \xi _1 ^0  \delta x_1 \right\}, 
\end{align*}
and using~(\ref{infinitgenBody}), the tangent space $\mathfrak{t}^2 \cdot {\bf z}_0:=T_{{\bf z}_0} \left(\mathbb{T}^2\cdot {\bf z}_0\right)$ to the toral orbit that goes through the relative equilibrium ${\bf z}_0 $ can be characterized as:
\begin{align*}
\mathfrak{t}^2 \cdot {\bf z}_0&= \left\{ \left( \xi _1 , \xi _2 \right)_{T^*SE(3)}(\mathbf{z}_0 )\mid \xi_1, \xi_2 \in \mathbb{R} \right\}= \left\{ \left(  \widehat{\left( \xi _1 \mathbf{e}_3 - \xi _2 \mathbf{e}_3  \right)}, \xi _1 \mathbf{e}_3 \times \mathbf{x}_0, \mathbf{0} , \xi _2 \mathbf{e}_3 \times \mathbf{p}_0 \right) \mid \xi _1 , \xi _2 \in \mathbb{R}^2\right\}\\
&= \left\{ \left( \left( \xi _1 - \xi _2 \right) \mathbf{e}_3 , \xi _1 r \mathbf{e}_2, \mathbf{0} , - \xi _2 M r\xi _1 ^0 \mathbf{e}_1 \right) \mid \xi _1 , \xi _2 \in \mathbb{R}^2 \right\}.  
\end{align*}
Finally, it can be easily verified that the vector subspace $W\subset \ker T_{{\bf z}_0} \boldsymbol{J}$ 
\begin{align*}
&W:= \Big\{ \left( \delta A , \boldsymbol{\delta}\mathbf{x} , \boldsymbol{\delta}\boldsymbol{\Pi} , \boldsymbol{\delta}\mathbf{p} \right)\mid \boldsymbol{\delta}\boldsymbol{\Pi} \cdot \mathbf{e}_3 =0, \boldsymbol{\delta}\mathbf{x} \cdot \mathbf{p} _0 =0, \boldsymbol{\delta}\mathbf{p} \cdot \mathbf{x} _0 =0, \delta p _2 =- M \xi _1 ^0  \delta x_1 \Big\}  
\end{align*}
is such that 
\begin{equation}
{\rm Ker} T_{{\bf z}_0} \boldsymbol{J} = W  \oplus \mathfrak{t}^2 \cdot {\bf z}_0,
\end{equation}
and hence constitutes a stability space. Moreover, let $\mathbf{u}_1 = \left( \mathbf{0} , \mathbf{e}_1, \mathbf{0} , -M \xi _1 ^0 \mathbf{e}_2  \right)$, $\mathbf{u}_2 = \left( \mathbf{e}_3 ,\mathbf{0} , \mathbf{0} , \mathbf{0} \right)$, $\mathbf{u}_3 =\left( \mathbf{0} , \mathbf{0} , \mathbf{e}_1, \mathbf{0} \right)$, $\mathbf{u}_4 =  \left( \mathbf{0} , \mathbf{0} , \mathbf{e}_2, \mathbf{0} \right)$, $\mathbf{u}_5= \left( \mathbf{0} , \mathbf{0} , \mathbf{0} , \mathbf{e}_3 \right)$, $\mathbf{u}_6= \left( \mathbf{0} , \mathbf{e}_3 , \mathbf{0} , \mathbf{0} \right)$, $\mathbf{u}_7=\left( \mathbf{e}_2, \mathbf{0} , \mathbf{0} , \mathbf{0}  \right)$, and $\mathbf{u}_8=  \left( \mathbf{e}_1, \mathbf{0} , \mathbf{0} , \mathbf{0} \right)$. It can be checked that 
\begin{align}
W={\rm span} \Big\{\mathbf{u}_1,\mathbf{u}_2,\mathbf{u}_3,\mathbf{u}_4,\mathbf{u}_5,\mathbf{u}_6,\mathbf{u}_7,\mathbf{u}_8
 \Big\}.\label{basis regular branch}
\end{align}
The set $\mathcal{B}= \Big\{\mathbf{u}_1,\mathbf{u}_2,\mathbf{u}_3,\mathbf{u}_4,\mathbf{u}_5,\mathbf{u}_6,\mathbf{u}_7,\mathbf{u}_8
 \Big\} $ will be used as a basis of the stability space in order to obtain  matrix expressions for the stability form $\left.\mathbf{d}^2 \left( h - \boldsymbol{J} ^{(\xi _1^0 , \xi _2 )}\right)\left( \mathbf{z}_0  \right)\right| _{W \times W}  $ corresponding to each part of  Theorem~\ref{kozorez relations}.
 
\medskip

\noindent {\bf A stability space for the singular branch ($r=0 $).} Consider now the relative equlibrium ${z} _0 = \left( \left(A_0, \mathbf{x}_0 \right) , \left(\boldsymbol{\Pi}_0,\mathbf{p}_0\right) \right) $ with $A_0=R_{\theta_0}^Z$, $\mathbf{x} _0= \left( 0, 0, 0\right) $, $\boldsymbol{\Pi} _0 = I_3 \left( \xi _1 - \xi _2 \right) \mathbf{e}_3 $, and $\mathbf{p} _0={\bf 0} $, $\xi_1, \xi_2 \in \mathbb{R} $.
In this case, the matrix expression~(\ref{hessian in general}) associated to $\mathbf{d}^2 \left(h- \mathbf{J}^{\boldsymbol{\xi}}\right) (\mathbf{z}) $ is given by:
\begin{equation}\label{Hessian_singular}
\left(
\begin{array}{cccc} 
 - \mu  \widehat{ \mathbf{B} (\mathbf{x}_0 ) } \widehat{\mathbf{e}_3 }  - \xi _1  \widehat{\mathbf{e}_3} \widehat{ \boldsymbol{\Pi_0} }  & 0 & \xi _1 \widehat{ \mathbf{e}_3 }& 0\\ 
0&-\mu {\rm Hess}  \left( {B_z} \right)(\mathbf{x}_0 ) &0&\xi _1 \widehat{\mathbf{e}_3 }\\
- \xi _1 \widehat{\mathbf{e}_3 }&0&\mathbb{I}_{ref} ^{-1}&0\\
0  &- \xi _1 \widehat{ \mathbf{e}_3 }&0&\dfrac{1}{M} \mathbb{I}_{id}
\end{array} 
\right).
\end{equation}
These relative equilibria lay on the singular isotropy type manifold~(\ref{singular orbit type}) and hence by the Bifurcation Lemma (see~\cite[Proposition 4.5.12]{Ortega2004}), the 
kernel of the derivative of the momentum map is necessarily of dimension eleven at those points. Indeed, it can be checked that:
\begin{equation*}
\ker T_{{\bf z}_0} \boldsymbol{J}= \left\{ v \in T_{\mathbf{z}_0 } \left( SE(3) \times \mathfrak{se}(3)^*\right) \mid \boldsymbol{\delta}\boldsymbol{\Pi} \cdot \mathbf{e}_3 =0 \right\}, 
\end{equation*}
and using~(\ref{infinitgenBody}), the tangent space $\mathfrak{t}^2 \cdot {\bf z}_0:=T_{{\bf z}_0} \left(\mathbb{T}^2\cdot {\bf z}_0\right)$ to the toral orbit that goes through the singular relative equilibrium ${\bf z}_0 $ can be characterized as:
\begin{equation*}
\mathfrak{t}^2 \cdot {\bf z}_0= \left\{ \left( \left( \xi _1 - \xi _2 \right) \mathbf{e}_3 , {\bf 0}, \mathbf{0} , {\bf 0} \right) \mid \xi _1 , \xi _2 \in \mathbb{R}^2 \right\}= {\rm span} \{ \left(\mathbf{e}_3, {\bf 0}, {\bf 0}, {\bf 0}\right)\}.  
\end{equation*}
Finally, it can be easily verified that the vector subspace $W\subset \ker T_{{\bf z}_0} \boldsymbol{J}$ 
given by
\begin{equation}
\label{stability space singular}
W={\rm span} \Big\{\mathbf{u}_1,\mathbf{u}_2,\mathbf{u}_3,\mathbf{u}_4,\mathbf{u}_5,\mathbf{u}_6,\mathbf{u}_7,\mathbf{u}_8, \mathbf{u}_9, \mathbf{u}_{10}
 \Big\},
\end{equation}
with 
$\mathbf{u}_1= \left( \mathbf{0} ,{\bf 0}, \mathbf{0} ,  \mathbf{e}_3  \right)$, 
$\mathbf{u}_2= \left( {\bf 0} ,\mathbf{e}_3, \mathbf{0} , \mathbf{0} \right)$, 
$\mathbf{u}_3 =\left( \mathbf{0} , \mathbf{0} , {\bf 0}, \mathbf{e}_2\right)$, 
$\mathbf{u}_4 =  \left( \mathbf{0} , \mathbf{0} , {\bf 0},\mathbf{e}_1 \right)$, 
$\mathbf{u}_5= \left( \mathbf{0} , \mathbf{0} , \mathbf{e}_2 , {\bf 0} \right)$, 
$\mathbf{u}_6= \left( \mathbf{0} , {\bf 0} , \mathbf{e}_1 , \mathbf{0} \right)$, 
$\mathbf{u}_7=\left( {\bf 0}, \mathbf{e}_2 , \mathbf{0} , \mathbf{0}  \right)$, 
$\mathbf{u}_8=  \left( {\bf 0}, \mathbf{e}_1 , \mathbf{0} , \mathbf{0} \right)$,
$\mathbf{u}_9=  \left( \mathbf{e}_2, {\bf 0} , \mathbf{0} , \mathbf{0} \right)$, and
$\mathbf{u}_{10}=  \left(\mathbf{e}_1 , {\bf 0}, \mathbf{0} , \mathbf{0} \right)$
is a $H$--invariant stability space, that is,
\begin{equation}
{\rm Ker} T_{{\bf z}_0} \boldsymbol{J} = W  \oplus \mathfrak{t}^2 \cdot {\bf z}_0,
\end{equation}
and hence constitutes a stability space. We recall that $H:=\left\{\left(e^{i\theta},e^{i\theta}\right)\mid e^{i\theta} \in S ^1\right\} = \mathbb{T}^2_{{\bf z}_0}$  is the isotropy subgroup of the relative equilibrium ${\bf z} _0 $. We will use the set $\mathcal{B}= \Big\{\mathbf{u}_1,\mathbf{u}_2,\mathbf{u}_3,\mathbf{u}_4,\mathbf{u}_5,\mathbf{u}_6,\mathbf{u}_7,\mathbf{u}_8,\mathbf{u}_9,\mathbf{u}_{10}
 \Big\} $ as a basis of the stability space in order to obtain a matrix expression for the stability form~$\left.\mathbf{d}^2 \left( h - \boldsymbol{J} ^{(\xi _1 , \xi _2 )}\right)\left( \mathbf{z}_0  \right)\right| _{W \times W}  $ for the parts {\bf (i)} and {\bf (ii)} of Theorem~\ref{kozorez relations}. 

\medskip

\noindent {\bf Proof of part (i) of the theorem.}

\medskip

\noindent {\bf Stability study for the regular branch.} We start by noting that  the  stability of a relative equilibrium can be determined by using any of the the points that constitute its trajectory in phase space. Hence we can, without loss of generality, use the relative equilibrium point ${\bf z} _0 $ of the form ${\bf z} _0 = \left( (\mathbb{I}_{id}, r \mathbf{e}_1), (I _3 \left( \xi _1 ^0 - \xi _2 \right) \mathbf{e}_3 , M r \xi _1 ^0 \mathbf{e_2}  )\right)$. We recall that the regular relative equilibria  are those for which $r> 0 $ and 
\begin{equation*}
\xi _1^0  = \pm \left( - \dfrac{3h\mu q \mu _0 }{2 \pi MD(\mathbf{x}_0 )^{5/2}}\right)^{1/2}.
\end{equation*}
We now provide the expression of ${\rm Hess}  \left( {B_z} \right)(\mathbf{x}_0 )$ using the same notation as in~(\ref{DBx0regular}) and conclude that 
\[
{\rm Hess}  \left( {B_z} \right)(\mathbf{x}_0 ) = \dfrac{6kh}{D(\mathbf{x}_0 )^{7/2}} \left( \begin{array}{ccc}  D(\mathbf{x}_0 )-5x^2& -5xy&0\\ -5xy&D(\mathbf{x}_0 )-5y^2&0\\ 0&0&3D(\mathbf{x}_0 )-5h^2 
\end{array} \right).
\]
By \eqref{hessian in particular} and \eqref{basis regular branch} we obtain that the stability form $\left.\mathbf{d}^2 \left( h - \boldsymbol{J} ^{(\xi _1^0 , \xi _2 )}\right)\left( \mathbf{z}_0  \right)\right| _{W \times W}  $ can be written as:
\begin{equation}
\!\!\!\!\!\!\!\!\!\!\!\!\!\!\!\!\!\!\!
\left( \begin {array}{cccccccc} 
M \xi_1^0 \dfrac{4 h^2-r^2}{r^2+h^2}
&0
&0
&0
&0
&0
&0
&0\\ 
0
&M {\xi_1^0}^2 r^2
&0
&0
&0
&0
&0
&0\\ 
0
&0
&{\dfrac {1}{I_1}}
&0
&0
&0
&\xi_1^0
&0\\ 
0
&0
&0
&\dfrac{1}{I_1}
&0
&0
&0
&- \xi_1^0 \\ 
0
&0
&0
&0
&\dfrac{1}{M}
&0
&0
&\xi_1^0 r\\ 
0
&0
&0
&0
&0
&M{\xi_1^0}^2 \dfrac{3 r ^2 - 2 h ^2 }{r^2+h^2}
&M { \xi_1 ^0 }^2  r
&0\\
0
&0
& \xi_1^0
&0
&0
&M{\xi_1^0}^2 r
&\dfrac{1}{3}M{\xi_1 ^0 } ^2(r^2+h^2)+\xi _1^0 \Pi _0  
&0\\
0
&0
&0
&- \xi_1^0
& \xi_1^0 r
&0
&0
&\dfrac{1}{3}M{\xi_1 ^0 } ^2 (4 r^2+h^2)+\xi _1^0 \Pi _0
\end {array} \right),
\end{equation}
where 
$\Pi _0 = I_3  \left( \xi _1 ^0 - \xi _2 \right)$. Notice that this matrix is block diagonal and exhibits two blocks of size two and six. The positivity of the block of size two requires that $M {\xi_1^0}^2 r^2 > 0 $ and
$4h ^2 - r ^2 >0$. The first inequality is always satisfied when $\mu q< 0 $ and the second one amounts to
\begin{equation} 
\dfrac{r^2}{h^2} < 4,
\end{equation} 
which yields the right hand side inequality in~(\ref{kozoriez}).
We now study the positivity of the lower six dimensional block of the stability form. As we already pointed out in Remark~\ref{why gaussian elimination}, given that by Sylvester's Law of Inertia the signature of a diagonalizable matrix is invariant with respect to conjugation by invertible matrices, it can hence be read out of the pivots of the matrix obtained by performing Gaussian elimination on this block. Indeed, these pivots are
\begin{equation}
\label{}
\dfrac{1}{I _1 }, \enspace \dfrac{1}{I _1 }, \enspace 
\dfrac{1}{M}, \enspace p _1, \enspace 
p _2, \enspace p _3,
\end{equation}
where 
\begin{eqnarray*}
p _1&:=&M{\xi_1^0} ^2 \dfrac{3 r ^2 - 2 h ^2 }{r^2+h^2},\\ 
p _2&:=&-I_3 \xi_1^0 \xi_2 - {\xi_1^0}^2 \left( \dfrac{2}{3} M \dfrac{ (r^2+h^2) h^2}{3 r^2-2 h^2}+(I_1-I_3) \right), \\ 
p _3&:=&-I_3 \xi_1^0 \xi_2 - {\xi_1^0}^2 \left( -\dfrac{1}{3} M (r^2+h^2)+(I_1-I_3) \right).
\end{eqnarray*} 
The first three are automatically positive. The positivity of $p _1$ is equivalent to
\begin{equation*} \label{radiuscond}
{\dfrac{2}{3}}<\dfrac{r^2}{h^2},
\end{equation*}
which yields the left hand side inequality in~(\ref{kozoriez})
Finally, we study the positivity of the last two pivots $p _2$ and $p _3$.
The simultaneous positivity of $p_2$ and $p_3$ is equivalent to $\min\left\{ p_2,p_3\right\} >0$. 
It is easy to check that $\min\left\{ p_2,p_3\right\} = p_2$, since the condition
\begin{equation}
p_3-p_2=M{\xi_1^0}^2  r^2 \dfrac{r^2+h^2}{3r^2-2h^2}>0
\end{equation}
is always satisfied due to \eqref{radiuscond} and the condition $\mu q<0$. Regarding the positivity of $p _2$ there are two possible cases:\\
1)
$\xi_1^0 >0$, then 
\begin{equation*}
 I_3 (\xi _1^0- \xi _2) >\xi_1^0 \left( I_1+\dfrac{2}{3} M \dfrac{(r^2+h^2)h^2}{3 r^2-2 h^2}\right),
\end{equation*}
2)
$ \xi_1^0<0$, then
\begin{equation*}
 I_3 (\xi _1^0- \xi _2) <\xi_1^0 \left( I_1+\dfrac{2}{3} M \dfrac{(r^2+h^2)h^2}{3 r^2-2 h^2}\right),
\end{equation*}
and hence the positivity of $p_2$ can be summarized as
\begin{equation*}
{\rm sign}(\xi _1^0 ) I_3 \xi _2<-| \xi_1^0  | \left( I_1-I_3+\dfrac{2}{3} M \dfrac{(r^2+h^2)h^2}{3 r^2-2 h^2}\right),
\end{equation*}
which coincides with~(\ref{signxi2}), as required.

\medskip

\noindent
{\bf Stability study for the singular branch.}
 We first notice that the matrix expression associated to $\mathbf{d}^2 \left(h- \mathbf{J}^\xi\right) (\mathbf{z}) $ is given by \eqref{Hessian_singular}, where 
\begin{equation}
\mathbf{B}(\mathbf{x}_0) = -\dfrac{\mu_0 q}{2 \pi h ^2} \mathbf{e}_3 \enspace {\rm and}
\enspace
{\rm Hess}  \left( {B_z} \right)(\mathbf{x}_0 ) = -\dfrac{3\mu_0q\mu}{2 \pi h ^4} \left( \begin{array}{ccc} 1&0&0\\ 0&1&0\\ 0&0&-2 \end{array} \right).
\end{equation}
Consequently, the pivots obtained by Gaussian elimination in the matrix expression of the stability form~$\left.\mathbf{d}^2 \left( h - \boldsymbol{J} ^{(\xi _1 , \xi _2 )}\right)\left( \mathbf{z}_0  \right)\right| _{W \times W}$ are
\begin{equation}
\label{singular pivots}
\frac{1}{M},\frac{3 \mu _0q \mu}{\pi h ^4}, \frac{1}{M}, \frac{1}{M}, \frac{1}{I _1}, \frac{1}{I _1}, p _1, p _1, p _2, p _2,
\end{equation}
where 
\begin{align*}
p _1 =- \frac{3}{2}\frac{\mu _0q \mu}{\pi h ^4}-M \xi_1^2,\quad
p _2= -\frac{ \mu _0q \mu}{2\pi h ^2}- \xi_1(\xi_1 I _1- \Pi_0), \quad \mbox{and} \quad \Pi_0=I _3(\xi_1- \xi_2).
\end{align*}
The formal instability of the singular branch is caused by the fact that the pivots in~(\ref{singular pivots}) cannot simultaneously have all the same sign. Indeed, $1/M $ and $1/I _1 $  are always positive which forces $3 \mu _0q \mu/\pi h ^4>0 $. This is in turn incompatible with $p _ 1 > 0$ because that would require $M \xi _1^2<0 $, which is not possible.

\medskip

\noindent {\bf Proof of part (ii) of the theorem.}

\medskip

\noindent {\bf Stability study for the regular branch.}
In order to prove the second part of Theorem~\ref{kozorez relations}, we follow the same pattern that we used above. Let $f\in C^{\infty}(\mathbb{R} ^2)$ be the function such that $B _z (x,y,z) = f(x ^2 + y ^2,z )$ and $f_0 := f(x ^2+ y ^2, 0)$. Additionally,
\begin{equation*}
f'_1 := \left. \dfrac{\partial f(v,z)}{\partial v}\right|_{v=x ^2+ y ^2, z= 0}, \enspace 
f_{1}'' := \left. \dfrac{\partial ^2 f(v,z)}{\partial v ^2}\right|_{v=x ^2+ y ^2, z= 0}, \quad 
f''_2 :=  \left. \dfrac{\partial ^2 f(v,z)}{\partial z ^2}\right|_{v=x ^2+ y ^2, z= 0},
\end{equation*}
and  we recall that $\xi _1 ^0 = \pm{ \left( -\dfrac{2}{M} \mu f'_1\right)}^{1/2}$. 
 We now compute the components of the matrix $D\mathbf{B} (\mathbf{x} _0) $. Using the equations  \eqref{Bxy_zero}, \eqref{deriv_Bz_zero} and \eqref{Bz_func_sqr_xy} we obtain
 \begin{align*}
\left. \dfrac{\partial B _x}{\partial x }\right|_{\mathbf{x} _0 } = 0, \  
\left. \dfrac{\partial B _x}{\partial y }\right|_{\mathbf{x} _0 } = 0, \  
\left. \dfrac{\partial B _y}{\partial x }\right|_{\mathbf{x} _0 }=0, \  
\left. \dfrac{\partial B _y}{\partial y }\right|_{\mathbf{x} _0 } = 0, \  
\left. \dfrac{\partial B _z}{\partial z }\right|_{\mathbf{x} _0 } = 0, \  
\left. \dfrac{\partial B _z}{\partial x }\right|_{\mathbf{x} _0 } = 2x f'_1, \  
\left. \dfrac{\partial B _z}{\partial y }\right|_{\mathbf{x} _0 } = 2y f'_1.
 \end{align*}
 In order to determine the remaining  two components in $D\mathbf{B} (\mathbf{x} _0) $, we use the Amp\`ere-Maxwell equation $\nabla \times \mathbf{B} = 0$ in the absence of additional currents and time-varying electric fields in the region where the body motion takes place. Indeed, $\nabla \times \mathbf{B} = 0$ implies that 
 \begin{equation*}
 \left. \dfrac{\partial B _x}{\partial z }\right|_{\mathbf{x} _0 }= \left. \dfrac{\partial B _z}{\partial x }\right|_{\mathbf{x} _0 } = 2x f'_1, \enspace
\left. \dfrac{\partial B _y}{\partial z }\right|_{\mathbf{x} _0 }= \left. \dfrac{\partial B _z}{\partial y }\right|_{\mathbf{x} _0 } = 2y f'_1,
\end{equation*}
and hence,
\begin{equation*}
D\mathbf{B}(\mathbf{x} _0 ) = \left(
\begin{array}{ccc} 
0 & 0 & 2x f'_1 \\ 
0&0&2y f'_1 \\
2x f'_1 &2y f'_1 &0
\end{array} 
\right).
\end{equation*}
By expression \eqref{TxF_is_Hess}
\begin{equation}
{T} _{\mathbf{x}_0  } \mathbf{F} (\cdot)(\mathbf{e}_3 ) = {\rm  Hess} \left( {B_z} \right)(\mathbf{x} _0 )=  \left(
\begin{array}{ccc} 
2f'_1+ 4x ^2  f''_1 & 4xy f''_1 & 0\\ 
4xy  f_{1}'' &2f'_1+ 4y ^2  f''_1 &0\\
0&0&f_{2}'' 
\end{array} 
\right).
\end{equation}
\noindent
Using the same argument as in the proof of part {\bf (i)} we use, without loss of generality, the relative equilibrium point ${\bf z} _0 $ of the form ${\bf z} _0 = ( (\mathbb{I}_{id}, r \mathbf{e}_1)$, $(I _3 \left( \xi _1 ^0 - \xi _2 \right) \mathbf{e}_3 , M r \xi _1 ^0 \mathbf{e_2} ) )$, where  $ r>0$. 
The matrix expression of $\left.\mathbf{d}^2 \left( h - \boldsymbol{J} ^{(\xi _1^0 , \xi _2 )}\right)\left( \mathbf{z}_0  \right)\right| _{W \times W}  $ is:

\begin{equation}
\!\!\!\!\!\!\!\!\!\!\!\!\!\!\!\!\!\!\!
\left( \begin {array}{cccccccc} 
-2\mu(f'_1+2r^2f_{1}'' )+3{\xi_1^0}^2M
&0
&0
&0
&0
&0
&0
&0\\ 
0
&M r ^2 {\xi_1^0}^2
&0
&0
&0
&0
&0
&0\\ 
0
&0
&{\dfrac {1}{I_1}}
&0
&0
&0
&\xi_1^0
&0\\ 
0
&0
&0
&\dfrac{1}{I_1}
&0
&0
&0
&- \xi_1^0 \\ 
0
&0
&0
&0
&\dfrac{1}{M}
&0
&0
&\xi_1^0 r\\ 
0
&0
&0
&0
&0
&-\mu f_{2}''
&-2\mu r f_1'
&0\\
0
&0
& \xi_1^0
&0
&0
&-2\mu r f_1'
&\mu f_0+\xi_1^0 \Pi_0
&0\\
0
&0
&0
&- \xi_1^0
& \xi_1^0 r
&0
&0
&\mu f_0+\xi_1^0 (Mr^2\xi_1^0+\Pi_0) 
\end {array} \right),
\end{equation}

\medskip

\noindent where 
$\Pi _0 = I_3  \left( \xi _1 ^0 - \xi _2 \right)$. Notice that this matrix is block diagonal and exhibits two blocks of size two and six. The positivity of the block of size two requires that $\mu f'_1<0 $ and $\mu(2 f_{1}' +r^2f_{1}'' )<0$ which coincide with~(\ref{f11}) and~(\ref{kozorez generalized1}). We now study the positivity of the lower six dimensional block of the stability form. As we did in the proof of part {\bf (i)}, we will  read the signature of this block out of its pivots, which are
$
\dfrac{1}{I _1 }
$,
$
\dfrac{1}{I _1 }
$,
$
\dfrac{1}{M}
$,
$
-\mu  f_{2}''
$,
$
\xi_1^0\left(\Pi_0 - \xi_1^0 I _1\right)+\mu \left(f_0+4  r ^2  \dfrac{{f'_{1}} ^2}{f_{2}'' } \right) 
$,
and 
$\mu f_0+\xi_1^0 \Pi_0-{\xi_1^0}^2 I_1
$.
The first three are automatically positive. The positivity of the fourth requires 
\begin{equation} \label{pivot4}
\mu f_{2}''<0,
\end{equation}
which corresponds to the inequality~(\ref{f22}) in the statement.
Finally, we study the positivity of the last two pivots.
Let 
\begin{equation} \label{b1}
p_1:=\xi_1^0\left(\Pi_0 - \xi_1^0 I _1\right)+\mu \left(f_0+4  r ^2  \dfrac{{f'_{1}} ^2}{f_{2}'' } \right) 
\end{equation}
and
\begin{equation} \label{b2}
p_2:=\mu  f_0  +\xi_1^0 \Pi_0-{\xi_1^0}^2 I_1.
\end{equation}
The simultaneous positivity of $p_1$ and $p_2$ is equivalent to $\min\left\{ p_1,p_2\right\} >0$. 
It is easy to check that $\min\left\{ p_1,p_2\right\} = p_1$, since the condition
\begin{equation}
p_2-p_1=-4\mu r^2 \dfrac{{f'_{1}} ^2}{f_{2}'' }>0
\end{equation}
is satisfied due to \eqref{pivot4}. The positivity of $p_1$ can be summarized as
\begin{equation*}
{\rm sign}(\xi _1^0 ) I_3 \xi _2<-  |\xi_1^0| \left(  (I_1-I_3)+\dfrac{1}{2}M\left(\dfrac{f _0}{f'_1}+4 r ^2\dfrac{{f_{1}' }}{f_{2}'' } \right) \right),
\end{equation*}
which coincides with \eqref{cond3_xi2}, as required.

\medskip

\noindent
{\bf Stability study for the singular branch.}
The matrix expression associated to $\mathbf{d}^2 \left(h- \mathbf{J}^\xi\right) (\mathbf{z}) $ is given by \eqref{Hessian_singular}, where in this case
\begin{equation}
\mathbf{B}(\mathbf{x}_0) =B _z(\mathbf{x}_0) \mathbf{e}_3 = f_0 \mathbf{e} _3\enspace {\rm and}
\enspace
{\rm Hess}  \left( {B_z} \right)(\mathbf{x}_0 ) =\left( \begin{array}{ccc} 2 f' _1&0&0\\ 0&2 f' _1&0\\ 0&0&f''_2 \end{array} \right).
\end{equation}
The pivots obtained by Gaussian elimination in the matrix expression of the stability form~$\left.\mathbf{d}^2 \left( h - \boldsymbol{J} ^{(\xi _1 , \xi _2 )}\right)\left( \mathbf{z}_0  \right)\right| _{W \times W}$ are
\begin{equation*}
p _1 ,p _2 , p _1 , p _1 , p _3 , p _3 , p _4, p _4, p _5, p _5,
\end{equation*}
where 
\begin{align*}
p _1 = \frac{1}{M}, \enspace
p _2 = - \mu f'' _2, \enspace
p _3 = \frac{1}{I _1}, \enspace
p _4 =- 2\mu f'_1-M \xi_1^2,\enspace
p _5= \mu f _0- \xi_1(\xi_1 I _1- \Pi_0), \enspace \mbox{and} \enspace \Pi_0=I _3(\xi_1- \xi_2).
\end{align*}
The pivots $p _1 $ and $p _3$ are automatically positive. The positivity of the pivots $p _2 $, $p _4 $, and $p _5 $ requires that:
\begin{align}
&\mu f'' _2 < 0,\label{f22_gen_singular_cond_proof}\\
&\mu f'_1<0,\label{f11_gen_singular_cond_proof} \\
&\xi_1^2 <   - \dfrac{2}{M}\mu f'_1,\label{xi1_gen_singular_cond_proof}\\
&{\rm sign}(\xi_1)\Pi_0> \dfrac{I _1\xi_1^2- \mu f_0}{|\xi_1|},\label{Pi0_gen_singular_cond_proof}
\end{align}
 which yields the conditions \eqref{f22_gen_singular_cond}--\eqref{Pi0_gen_singular_cond}. 
 We now derive the optimal stability condition \eqref{Pi0_gen_singular_optimal_cond}. Let $g(\xi_1):={\left(I _1 \xi_1 ^2-\mu f _0\right)}/{\xi_1}$ in  \eqref{Pi0_gen_singular_cond}. It is easy to verify that the function $g(\xi_1)$ has a minimum at $\widehat{\xi_1^{+}}=\sqrt{- \mu f _0/I _1}$ and a maximum at $\widehat{\xi_1^{-}}=-\sqrt{- \mu f _0/I _1}$ provided that $\mu f_0<0$. Since the condition \eqref{xi1_gen_singular_cond} has to be satisfied, then $f _0/ f '_1<2I _1/M$ also needs to hold. In that case, the choices $\widehat{\xi_1^{\pm}}= \pm\sqrt{- \mu f _0/I _1}$ and the inequalities
 \begin{equation*}
\Pi_0>\mathop{\rm min}_{\xi_1\in \mathbb{R}^{+}}\left\{g (\xi_1)\right\}=g\left(\widehat{\xi_1^{+}}\right)=2 \sqrt{- \mu f _0 I _1}, \qquad 
\Pi_0<\mathop{\rm max}_{\xi_1\in \mathbb{R}^{-}}\left\{g (\xi_1)\right\}=g\left(\widehat{\xi_1^{-}}\right)=-2 \sqrt{- \mu f _0 I _1}
\end{equation*}
determine the largest possible stability region in the $\Pi_0$ (and consequently the $\xi_2$) variable, as required in~(\ref{Pi0_gen_singular_optimal_cond}). \quad $\blacksquare$

\subsection{Proofs of Propositions~\ref{linear tools for instability} and~\ref{linearization for t stars g}}

\noindent {\bf Proof of Proposition~\ref{linear tools for instability}}

\medskip

\noindent {\bf (i)} It is a consequence of the Witt-Artin decomposition (see for example~\cite[Theorem 7.1.1]{Ortega2004}).

\medskip

\noindent {\bf (ii)} It is a consequence of the fact that the symplectic slice introduced by Marle~\cite{nfm1, nfm}, Guillemin, and Sternberg~\cite{nfgs} can be constructed by Riemannian exponentiation of a symplectic tube. Since we need this construction in the proof of the following parts of the proposition, we briefly recall it using the notation in Chapter 7 of~\cite{Ortega2004}.

The first step is the splitting of the Lie algebra $\mathfrak{g} $ of $G$ into three parts. The first summand is $\mathfrak{g}_\mu:= {\rm Lie} \left(G _\mu\right) $. The equivariance of the momentum map $\mathbf{J} $ implies that $G _m \subset G _\mu  $ and hence $\mathfrak{g}_{m} \subset \mathfrak{g}_{\mu}  $. Hence we can fix an $\mbox{\rm Ad} _{G_{m} }
$-invariant inner product $\langle \cdot , \cdot \rangle $ on $\mathfrak{g}$ (always
available by the compactness of $G_{m} $) and write
\begin{equation}
\label{splitting of lie algebras for slice}
\mathfrak{g}_\mu= \mathfrak{g}_{m}\oplus \mathfrak{m} \quad\text{and}\quad \mathfrak{g}=
\mathfrak{g}_{m}\oplus \mathfrak{m}\oplus \mathfrak{q},
\end{equation}
where $\mathfrak{m}$ is the $\langle \cdot , \cdot \rangle $-orthogonal complement of
$\mathfrak{g}_{m}$ in $\mathfrak{g}_\mu$ and $\mathfrak{q}$ is the $\langle \cdot , \cdot
\rangle $-orthogonal complement of $\mathfrak{g}_\mu$ in $\mathfrak{g}$. The splittings
in~(\ref{splitting of lie algebras for slice}) induce similar ones on the duals 
\begin{equation}
\label{splitting of lie algebras for slice dual}
\mathfrak{g}_\mu^\ast = \mathfrak{g}_{m}^\ast \oplus \mathfrak{m}^\ast  \quad\text{and}\quad
\mathfrak{g}^\ast =
\mathfrak{g}_{m}^\ast \oplus \mathfrak{m}^\ast \oplus \mathfrak{q}^\ast. 
\end{equation}
Each of the spaces in this decomposition should be understood as the set of covectors in
$\mathfrak{g}^\ast $ that can be written as $\langle \xi, \cdot \rangle $, with $\xi $ in
the corresponding subspace. For example, $ \mathfrak{q}^\ast =\{ \langle \xi, \cdot
\rangle\mid  \xi \in \mathfrak{q}\} $.

The second ingredient in the construction of the symplectic tube comes from noting that the compact (by the properness of the action) isotropy subgroup $G _m $ acts linearly and canonically  on $(W, \omega _W)$ with momentum map $\mathbf{J}_W:W \rightarrow \mathfrak{g}_m^\ast  $ given by 
\begin{equation*}
\langle\mathbf{J}_W(w), \eta\rangle= \frac{1}{2}\omega_W \left(\eta_W(w), w\right), \qquad \eta\in \mathfrak{g}_m.
\end{equation*}
It can be shown~\cite[Proposition 7.2.2]{Ortega2004} that there exist $G _m $--invariant neighborhoods $\mathfrak{m}_r ^\ast  $ and $W _r $ of the origin in $\mathfrak{m}^\ast  $ and $W$, respectively, such that  the twisted product $Y _r:=G \times _{G _m } \left(\mathfrak{m}_r ^\ast  \times W _r\right) $ is endowed with a natural symplectic form $\omega_{Y _r} $ whose expression can be found in (7.2.2) of~\cite{Ortega2004}. The Lie group $G $ acts canonically on $\left(Y _r, \omega_{Y _r}\right) $ by $g \cdot \left[h, \eta, w\right]= \left[gh, \eta,w\right] $, for any $g \in  G $ and $[h, \eta, w ]\in Y _r   $, and has a momentum map $\mathbf{J}_{Y _r}: Y _r \rightarrow \mathfrak{g}^\ast $ associated given by the so called Marle--Guillemin--Sternberg normal form:
\begin{equation*}
\mathbf{J}_{Y _r} \left(\left[g, \eta, w\right]\right)= {\rm Ad} ^\ast _{g ^{-1}} \left( \mu+ \eta + \mathbf{J}_W (w)\right), \qquad \left[g, \eta, w\right] \in Y _r.
\end{equation*}
The $G$--symplectic manifold $\left(Y _r, \omega_{Y _r}\right) $ is called a symplectic tube of $(M, \omega) $ at the point $m$. This denomination is justified by the Symplectic Slice Theorem~\cite{nfm1, nfm, nfgs} that proves the existence of a $G$--equivariant symplectomorphism $\phi : U \rightarrow Y _r$ between a $G$--invariant neighborhood $U$ of $m$ in $M$ and $Y _r $ satisfying $\phi(m)=[e,0,0] $. The symplectic slice $S$ in the statement of the proposition is obtained~\cite[Theorem 7.4.1]{Ortega2004} as $S= \phi^{-1} \left(S_{Y _r}\right)$, where $S_{Y _r}:=\left\{[e,0,w]\mid w \in W _r\right\}$ and, more explicitly, as $S= \left\{{\rm Exp} _m (w)\mid w \in W _r\right\} $, with ${\rm Exp} _m  $ the Riemannian exponential associated to a $G _m$--invariant metric. The identity $T _mS=W $ is a consequence of the fact that $T _0 {\rm Exp} _m= {\rm Id} $.

\medskip

\noindent {\bf (iii)} Since $m\in M$ is a relative equilibrium, we have $\mathbf{d} \left(h- \mathbf{J}^\xi\right)(m)=0 $. This implies that $\mathbf{d}h ^\xi_S (m)= \mathbf{d}\left.\left(h- \mathbf{J}^\xi\right)(m)\right|_{T _mS} =0 $ and hence $X_{h _S^\xi}(m) =0 $.

\medskip

\noindent {\bf (iv)} This statement is a consequence of the combination of {\bf (ii)} and {\bf (iii)} with the following lemma. 

\begin{lemma} 
\label{Hessian_Q}
Let $(M, \omega )$ be a symplectic manifold, $h \in C^{\infty}(M)$, and $X_h$ the corresponding Hamiltonian vector field. Suppose that $m_0 \in M$ is an equilibrium point of $X_h$, that is $X_h(m_0) = 0$ and consequently $\mathbf{d}h(m_0) = 0 $. Then, the linearization $X'$ of $X_h$ at $m_0$ is a Hamiltonian vector field on the symplectic vector space $(T_{m_0} M, \omega (m_0))$ with Hamiltonian function $Q \in C^{\infty}(T_{m_0}M)$ given by 
\begin{equation} \label{HessianQ}
Q(v) = \dfrac{1}{2} \mathbf{d} ^2 h(m_0) (v,v).
\end{equation}
\end{lemma}

\noindent\textbf{Proof of the Lemma.\ \ }
Note $V=T_{m_0} M$ and $\omega _V = \omega (m_0)$. Let $v,w \in V$ arbitrary and let $\left\{ c(s) \vert s \in \mathbb{R}\right\} $ be a curve such that $v = \left. \frac{d}{ds } \right|_{s=0} c(s) $. Then if $F_t$ is the flow of $X_h$, we write
\begin{equation} \label{omegaX}
\omega _V \left(X'(v), w\right) = \left.\frac{d}{dt}\right|_{t=0} \omega (m_0) \left( T_{m_0} F_t \cdot v, w\right) = \left.\frac{d}{dt}\right|_{t=0} \omega (m_0) \left(\left.\frac{d}{ds}\right|_{s=0} F_t (c(s)), w\right), \quad v,w \in V. 
\end{equation}
We now take a Darboux chart $(U, \phi )$~\cite[page 75]{Abraham1978} around the point $m_0$. Recall that in Darboux coordinates, the symplectic form $\omega _U$ is constant. Additionally if $\phi : U \longrightarrow \phi (U) \subset \mathbb{R}^n$, let $u \in \mathbb{R}^n$ be such that $T_{m_0} \phi \cdot w = \left( \phi (m_0), u \right) \in \phi (U) \times \mathbb{R}^n = T(\phi (U))$. Now, since $\phi ^\ast \omega _U = \omega \vert_U$, then \eqref{omegaX} can be written as 
\begin{align*}
\left.\frac{d}{dt}\right|_{t=0} \omega (m_0) \left( \left.\frac{d}{ds}\right|_{s=0}F _t (c (s)), w\right) &= \left.\frac{d}{dt}\right|_{t=0} \omega _U \left( \left.\frac{d}{ds}\right|_{s=0} \phi \cdot F _t (c(s)), T_{m_0} \phi \cdot w\right) \\
&=\omega _U \left( \left.\frac{d}{ds}\right|_{s=0} T_{c(s)} \phi \cdot X _h (c(s)), \left( \phi (m_0), u\right) \right) \\
	&= \left.\frac{d}{ds}\right|_{s=0} \omega _U \left( T_{c(s)} \phi \cdot X _{h \circ \phi ^{-1} \circ \phi } ( c(s)) , ( \phi (c(s)), u) \right)   \\
& =\left.\frac{d}{ds}\right|_{s=0} \omega _U \left( X _{h \circ \phi ^{-1}} ( \phi (c(s))) , ( \phi (c(s)), u) \right)  \\
	&= \left.\frac{d}{ds}\right|_{s=0}\mathbf{d} \left( h\circ \phi ^{-1} \right) \left(  \phi (c(s))\right)  \cdot \left(  \phi (c(s)), T_{m_0} \phi \cdot w \right)  \\
& =\mathbf{d^2} \left( h\circ \phi ^{-1} \right) \left(  \phi (m_0)\right) \left(  (\phi (m_0), T_{m_0} \phi \cdot v), (\phi (m_0), T_{m_0} \phi \cdot w) \right)  \\
&= \mathbf{d^2}  h(m_0) (v,w) = \mathbf{d} Q(v) \cdot w.
\end{align*}
Consequently, $\mathbf{i}_{X'} \omega _V = \mathbf{d} Q$, as required. \quad $\blacktriangledown$

\medskip

\noindent {\bf (v)} The hypothesis $T_m \left( G_\mu \cdot m \right) = T_m \left( G \cdot m\right) $ implies that $ \mathfrak{q} \cdot m:=\{ \xi_M(m)\mid \xi\in \mathfrak{q}\}  $; this fact and the construction of the Witt--Artin decomposition (see for example the expression (7.1.11) in~\cite{Ortega2004}) ensure that \eqref{TmM} holds. 
In order to prove \eqref{XQ}, notice that for any $w _1 , w _2 \in W$ 
\begin{align} \label{XQ_proof}
\omega _W (X_Q(w _1 ), w _2 ) &= \mathbf{d} Q (w _1 )\cdot w _2 = \mathbf{d} ^2 h ^ \xi (m) (w _1 , w _2 )=\omega (m) (X_{h^ \xi }'(w _1 ), w _2 )\notag \\
	&= \omega (m) \left( \mathbb{P}_W X_{h^ \xi }'(w _1 ), w _2 \right) + \omega (m) \left( (\mathbb{I} - \mathbb{P} _W) X_{h^ \xi }' (w _1 ), w _2 \right)\notag \\
&= \omega (m) \left( \mathbb{P}_W X_{h^ \xi }'(w _1 ), w _2 \right) =  \omega_W \left( \mathbb{P}_W X_{h^ \xi }'(w _1 ), w _2 \right),\notag
\end{align}
where we used that $ (\mathbb{I} - \mathbb{P} _W) X_{h^ \xi }' (w _1 ) \in W ^\omega $ and hence $\omega (m) \left( (\mathbb{I} - \mathbb{P} _W) X_{h^ \xi }' (w _1 ), w _2 \right)=0 $. Since $w _1, w _2 \in W$ are arbitrary, the equality $\omega _W (X_Q(w _1 ), w _2 )=\omega_W \left( \mathbb{P}_W X_{h^ \xi }'(w _1 ), w _2 \right) $
implies that 
\begin{equation*}
X_Q(w) = \mathbb{P}_W X_{h^ \xi }' (w), \enspace \forall w \in W,
\end{equation*}
which is equivalent to \eqref{XQ}.

\medskip

\noindent {\bf (vi)} Given the local and group invariant  character of this statement, we will prove this statement using the so called reconstruction differential equations~\cite{thesis, Roberts2002, Ortega2004} that determine the Hamiltonian vector field associated to a $G$--invariant Hamiltonian $h \in C^{\infty}(Y _r) $ in the symplectic tube $Y _r $. Consider first $ \pi: G \times  \mathfrak{m} ^\ast _r \times W _r \longrightarrow G \times _{G _m}\left(\mathfrak{m}^\ast _r \times  W _r\right)= Y _r $ the orbit projection; the $G$--invariance of $h$ implies that the composition $h \circ \pi$ can be understood as a $G _m$--invariant function on $G \times  \mathfrak{m} ^\ast _r \times W _r $ that does not depend of the first factor, that is,  $h \circ \pi \in C^{\infty} \left(\mathfrak{m}_r ^\ast  \times  W _r\right)^{G _m} $. The reconstruction equations show that for any $[g, \rho, w] \in Y _r $, 
\[
X _h([g, \rho, w])=T_{(g, \rho, w)} \pi(X_{\mathfrak{m}}( g, \rho,w  ),
X_{\mathfrak{m}^\ast_r}( g, \rho,w  ), X_{W}( g, \rho,w)), 
\] 
where $X_{\mathfrak{m}}( g, \rho,w  )$,
$X_{\mathfrak{m}^\ast_r}( g, \rho,w  )$, and $ X_{W_r}( g, \rho,w) $ are determined by the expressions
\begin{eqnarray}
X_{\mathfrak{m}}( g, \rho,w)&=&T _eL _g(D_{\mathfrak{m}^\ast _{r}}(h\circ\pi)(\rho,w)),
\label{field 1 *}\\
X_{W_r}( g, \rho,w) &=&\omega^\sharp_{W}(D_{W_{r}}(h\circ\pi)(\rho,w)),\label{field 2 *}\\
X_{\mathfrak{m}^\ast _{r}}( g, \rho,w) &=&\mathbb{P}_{\mathfrak{m}^\ast }\Bigl({\rm ad}^\ast _{D_{\mathfrak{m}^\ast _{r}}(h\circ\pi)}
     \rho\Bigr)+{\rm ad}^\ast _{D_{\mathfrak{m}^\ast _{r}}(h\circ\pi)}\mathbf{J}_{W}(w),\label{field 3 *}
\end{eqnarray}
where $\mathbb{P}_{\mathfrak{m}^\ast }: \mathfrak{g}^\ast\rightarrow \mathfrak{m}^\ast  $ is the projection according to the splitting~(\ref{splitting of lie algebras for slice dual}) and $\omega^\sharp_{W}: W ^\ast  \rightarrow W $ is the isomorphism associated to the symplectic form $\omega_W$  in $W $.

We now assume that $X _Q $ is spectrally unstable, which implies by part {\bf (iv) } that the Hamiltonian vector field $X_{h ^\xi_S }  $ on the symplectic slice $S_{Y _r}=\{ [e,0,w]\mid w \in W _r\} $ exhibits an unstable equilibrium at $[e,0,0] $. Notice now that the Hamiltonian vector field $X_{h ^\xi_S }  $ is given by the projection onto $Y _r $ via $T \pi $ of the vector field $(0,0, X_{W _r}^{h _S ^\xi}) $ in $G \times  \mathfrak{m}^\ast _r \times W _r $ determined by
\begin{equation*}
X_{W _r}^{h _S ^\xi}= \omega^\sharp _W(D_{W _r}\left(h \circ \pi\right)(0,w))- \left(\mathbb{P}_{\mathfrak{g}_m}\xi\right)_W(w)=X_{W}^h(e,0,w)-\left(\mathbb{P}_{\mathfrak{h}}\xi\right)_W(w),
\end{equation*}
where $\mathbb{P}_{\mathfrak{g}_m}: \mathfrak{g} \rightarrow \mathfrak{g}_m $ is the projection according to the splitting~(\ref{splitting of lie algebras for slice}), $\left(\mathbb{P}_{\mathfrak{g}_m}\xi\right)_W\in  \mathfrak{X}(W _r)$ is the infinitesimal generator associated to $\mathbb{P}_{\mathfrak{g}_m}\xi \in \mathfrak{g}_m $ using the $G _m$--action on $W _r $, and $ X_{W}^h $ is the  vector field in~(\ref{field 2 *}) that determines the dynamics induced by $h$ on the space $W _r $. The instability of $X_{h _S ^\xi} $ at $[e,0,0] $ implies the same feature for $X_W ^h $ at $(e,0,0) $ and hence the $K $--instability of the relative equilibrium $[e,0,0] $ of $X _h $, for any subgroup $K \subset G $. \quad $\blacksquare$ 

\medskip

\noindent {\bf Proof of Proposition~\ref{linearization for t stars g}}

\medskip

\noindent\textbf{(i)} It is a straightforward consequence of the chain rule as $\mathbf{d}h (g, \mu)=0 $ implies that $\mathbf{d}h ^g(e, \mu)=0 $.

\medskip

\noindent {\bf (ii)} Relation~(\ref{relation quadratic forms}) is a consequence of the following general fact about Hessians: let $m\in M$ and $n\in N$, with $M$ and $N$ smooth manifolds and let $\psi:M\rightarrow N$ be a smooth
map such that $\psi(m)=n$. Let $f\in C^{\infty}(N)$ with $\mathbf{d}
f(n)=0$. Then $\mathbf{d}^2(\psi^\ast  f)(m)=T_m^\ast\psi(\mathbf{d}^2 f(n))$,
that is, for any $v,\,w\in T_mM$:
\[
\mathbf{d}^2(\psi^\ast  f)(m)(v,\,w)=
\mathbf{d}^2 f(n)(T_m\psi\cdot v,\,T_m\psi\cdot w).
\]
In order to establish relation~(\ref{relation linearizations}) note that the map $\Phi _g:  \mathfrak{g}\times \mathfrak{g}^\ast \longrightarrow T_{(g, \mu)}\left(G \times  \mathfrak{g}^\ast\right) $ is a symplectomorphism and hence a Poisson map. Expression~(\ref{relation linearizations}) follows from~\cite[Proposition 4.1.19]{Ortega2004}. 

\medskip

\noindent {\bf (iii)} By Lemma~\ref{Hessian_Q}, the Hamiltonian vector field $X_{Q ^g} $ is determined by the relation $\mathbf{i}_{X_{Q ^g}} \omega(e, \mu)= \mathbf{d} Q ^g  $ or, more explicitly, by
\begin{equation*}
\omega(e, \mu)\left(X_{Q^g}(\xi, \tau), (\eta, \rho)\right)= \mathbf{d} Q ^g(\xi, \tau)\cdot (\eta, \rho), \quad \mbox{for any} \quad (\xi, \tau), (\eta, \rho) \in \mathfrak{g}\times \mathfrak{g}^\ast.
\end{equation*}
Using the expression of the canonical symplectic form  of $T ^\ast G $ in body coordinates (see for instance~\cite[Expression (6.2.5)]{Ortega2004}), this equality can be rewritten as 
\begin{equation}
\label{intermediate 1}
\langle \rho, \pi_{\mathfrak{g}} \left(X_{Q ^g}(\xi, \tau)\right)\rangle- \langle\pi_{\mathfrak{g}^\ast} \left(X_{Q ^g}(\xi, \tau)\right), \eta\rangle+ \langle\mu, {\rm ad}_{\pi_{\mathfrak{g}} \left(X_{Q ^g}(\xi, \tau)\right)}\eta\rangle=\langle{\rm Hess}(\xi, \tau), (\eta, \rho)\rangle.
\end{equation}
Let now ${\rm {\bf pr}}_{\mathfrak{g}}: \mathfrak{g}\times \mathfrak{g}^\ast \rightarrow \mathfrak{g}\times \mathfrak{g}^\ast $ and ${\rm {\bf pr}}_{\mathfrak{g}^\ast }: \mathfrak{g}\times \mathfrak{g}^\ast \rightarrow \mathfrak{g}\times \mathfrak{g}^\ast $ be the maps defined by ${\rm {\bf pr}}_{\mathfrak{g}}(\eta, \rho):=(\eta, 0)$ and ${\rm {\bf pr}}_{\mathfrak{g}^\ast }(\eta, \rho):=(0, \rho)$, for any $(\eta, \rho)\in \mathfrak{g}\times \mathfrak{g}^\ast$, and $\mathbf{i}_{\mathfrak{g}}: \mathfrak{g}\rightarrow \mathfrak{g}\times \mathfrak{g}^\ast $ and $\mathbf{i}_{\mathfrak{g}^\ast }: \mathfrak{g}^\ast \rightarrow \mathfrak{g}\times \mathfrak{g}^\ast $, the canonical injections. Since~(\ref{intermediate 1}) holds for $(\eta, \rho)\in \mathfrak{g}\times \mathfrak{g}^\ast$ arbitrary, we apply it to vectors of the form ${\rm {\bf pr}}_{\mathfrak{g}}(\eta, \rho)=(\eta, 0)$ and ${\rm {\bf pr}}_{\mathfrak{g}^\ast }(\eta, \rho)=(0, \rho)$ and we obtain the following two equalities
\begin{eqnarray*}
\langle{\rm Hess}(\xi, \tau), {\rm {\bf pr}}_{ \mathfrak{g}^\ast }(\eta, \rho)\rangle&=&\langle \pi_{\mathfrak{g}^\ast }(\eta, \rho), \pi_{\mathfrak{g}} \left(X_{Q ^g}(\xi, \tau)\right)\rangle,\\
\langle{\rm Hess}(\xi, \tau), {\rm {\bf pr}}_{ \mathfrak{g}}(\eta, \rho)\rangle&=&- \langle\pi_{\mathfrak{g}^\ast} \left(X_{Q ^g}(\xi, \tau)\right), \pi_{\mathfrak{g}}(\eta, \rho)\rangle+ \langle\mu, {\rm ad}_{\pi_{\mathfrak{g}} \left(X_{Q ^g}(\xi, \tau)\right)}\pi_{\mathfrak{g}}(\eta, \rho)\rangle.
\end{eqnarray*}
Since $(\eta, \rho)\in \mathfrak{g}\times \mathfrak{g}^\ast$ in these expressions are arbitrary, they can be rewritten as:
\begin{eqnarray}
{\rm {\bf pr}}_{ \mathfrak{g}^\ast }^\ast {\rm Hess}(\xi, \tau)&=&\pi_{\mathfrak{g}^\ast }^\ast \left( \pi_{\mathfrak{g}} \left(X_{Q ^g}(\xi, \tau)\right)\right),\label{iw1}\\
{\rm {\bf pr}}_{ \mathfrak{g}}^\ast {\rm Hess}(\xi, \tau)&=&-  \pi_{\mathfrak{g}}^\ast  \left(\pi_{\mathfrak{g}^\ast} \left(X_{Q ^g}(\xi, \tau)\right)\right)+ \pi_{\mathfrak{g}}^\ast \left({\rm ad}_{\pi_{\mathfrak{g}} \left(X_{Q ^g}(\xi, \tau)\right)}^\ast (\mu)\right).\label{iw2}
\end{eqnarray}
We now apply $\pi_{\mathfrak{g}^\ast} $ and $\pi_{\mathfrak{g}} $ to both sides of~(\ref{iw1}) and~(\ref{iw2}), respectively,  and we notice that $\pi_{\mathfrak{g}^\ast} \circ {\rm {\bf pr}}_{ \mathfrak{g}^\ast}^\ast = \pi_{\mathfrak{g}^\ast} $, $\pi_{\mathfrak{g} } \circ {\rm {\bf pr}}_{ \mathfrak{g}} ^\ast  = \pi_{\mathfrak{g} } $, $\pi_{\mathfrak{g}} ^\ast = \mathbf{i}_{\mathfrak{g}}$, $\pi_{\mathfrak{g}^\ast } ^\ast = \mathbf{i}_{\mathfrak{g}^\ast }$,  $\pi_{\mathfrak{g}} \circ \mathbf{i}_{\mathfrak{g}}= {\rm id}_{\mathfrak{g}} $, and $\pi_{\mathfrak{g}^\ast } \circ \mathbf{i}_{\mathfrak{g}^\ast }= {\rm id}_{\mathfrak{g}^\ast } $. We obtain
\begin{eqnarray}
\pi_{\mathfrak{g}} \left(X_{Q ^g}(\xi, \tau)\right) &= &\pi_{ \mathfrak{g}^\ast } {\rm Hess}(\xi, \tau), \\
\pi_{\mathfrak{g}^\ast} \left(X_{Q ^g}(\xi, \tau)\right)&=&-\pi_{\mathfrak{g}} \left({\rm Hess}(\xi, \tau)\right)+ {\rm ad}_{\pi_{ \mathfrak{g}^\ast } {\rm Hess}(\xi, \tau)}^\ast \mu, 
\end{eqnarray}
which is equivalent to~(\ref{expression linearization at e}). \quad $\blacksquare$

\subsection{Linear stability and instability of the standard and generalized orbitron relative equilibria}
\label{Linear stability and instability of the orbitron relative equilibria}

\noindent {\bf The linearization for the regular branches.} The goal in this paragraph is determining the linear Hamiltonian vector fields $X _Q $ in the stability space $W$ used in the proof of Theorem~\ref{kozorez relations} by using the expression~(\ref{XQ}) in Proposition~\ref{linear tools for instability}. Notice that this is indeed possible due to the Abelian character of our symmetry group that ensures that in this situation $G _\mu=G $  and hence the coincidence of the tangent spaces $T _m \left( G_{\mu } \cdot m\right) $ and  $T _m \left( G \cdot m\right) $ that is necessary as a hypothesis in this statement. We start by writing down the decomposition~(\ref{TmM}) that in this case can be achieved by noting that \begin{equation}
W^{\omega } = {\rm span} \left\{ \boldsymbol{\tau} _1, \boldsymbol{\tau} _2 , \boldsymbol{\tau} _3 , \boldsymbol{\tau} _4  \right\} 
\end{equation}
with 
$
\boldsymbol{\tau} _1 = (\mathbf{e}_3 , r \mathbf{e}_2, \mathbf{0}, \mathbf{0}  )\, 
\boldsymbol{\tau} _2 = (-\mathbf{e}_3 , \mathbf{0} , \mathbf{0} , - Mr \xi _1 ^0  \mathbf{e}_1),\,
\boldsymbol{\tau} _3 = (\mathbf{0} , \mathbf{0}, \mathbf{0} , \mathbf{e}_1  )\,
\boldsymbol{\tau} _4 = (\mathbf{0} , \mathbf{e}_1, -Mr \xi _1 ^0 \mathbf{e}_3 , \mathbf{0}  )$.
If we use as a basis for $W$ the vectors introduced in~(\ref{basis regular branch}) and for $W ^\omega $ the ones those that we just described, it is easy to see that the matrix expressions of the inclusion $\boldsymbol{i}_W: W \hookrightarrow T_{\mathbf{z}_0} \left( SE(3) \times \mathfrak{se}(3)^*\right) \simeq \mathbb{R}^{12}$ and  the projection  $\mathbb{P}_W: T_{\mathbf{z}_0} \left( SE(3) \times \mathfrak{se}(3)^*\right) \simeq \mathbb{R}^{12} \longrightarrow W$,   where $\mathbb{R}^{12}$ is endowed with the canonical basis, are given by 
\begin{eqnarray}
\boldsymbol{i}_W &=& \left(\mathbf{u}_1'\vert  \mathbf{u}_2'\vert \mathbf{u}_3'\vert \mathbf{u}_4'\vert \mathbf{u}_5'\vert \mathbf{u}_6'\vert \mathbf{u}_7'\vert \mathbf{u}_8'\right),\label{injection 1}\\
 \mathbb{P}_W &=& \left(\mathbf{e}_8'\vert  \mathbf{e}_7'\vert \mathbf{e}_2'\vert \mathbf{e}_1'\vert -\frac{1}{r}\mathbf{e}_2'\vert \mathbf{e}_6'\vert \mathbf{e}_3'\vert \mathbf{e}_4' \vert -\frac{1}{Mr \xi _1 ^0 } \mathbf{e}_2' \vert \mathbf{0} \vert \mathbf{e}_5' \right), \label{projection 1}
\end{eqnarray}
where the apostrophes stand for the transposition operation and the vertical indicate matrix concatenation. The linearization of the Hamiltonian vector field associated to the augmented Hamiltonian at the relative equilibria of the standard orbitron can be immediately obtained by using the expression~(\ref{expression of linearization}) together with the Hessian already computed in~(\ref{hessian in general}) in the context of the proof of Theorem~\ref{kozorez relations}. The resulting expression is inserted in~(\ref{XQ}) using the injection~(\ref{injection 1}) and the projection~(\ref{projection 1}) and yields the following matrix  for $X _Q $: 
\medskip

\begin{equation}
\label{linearization regular}
{\small 
\!\!\!\!\!
\left( \begin {array}{cccccccc} 
0&-\xi_1^0 r&0&0&0&0&0&0\\
\xi_1^0 r \dfrac{4 h^2-r^2}{r^2+h^2}&0&0&0&0&0&0&0\\
0&0&0&\xi_1^0-\dfrac{\Pi_0}{I_1}&0&0&0&-\frac{1}{3} {\xi_1^0}^2 M (r^2+h^2)\\
0&0&-\xi_1^0+\dfrac{\Pi_0}{I_1}&0&0&-{\xi_1^0}^2 M r&-\frac{1}{3} {\xi_1^0}^2 M (r^2+h^2)&0\\
0&0&-\dfrac{M \xi_1^0 r}{I_1}&0&0&{\xi_1^0}^2 M\dfrac{2 h^2-3 r^2}{r^2+h^2}&-2 {\xi_1^0}^2 M r&0\\
0&0&0&0&\dfrac{1}{M}&0&0&\xi_1^0 r\\
0&0&0&\dfrac{1}{I_1}&0&0&0&-\xi_1^0\\
0&0&\dfrac{1}{I_1}&0&0&0&\xi_1^0&0
\end {array} \right),}
\end{equation}
\medskip
where
$\Pi_0=I_3 (\xi_1^0-\xi_2)$. An analog expression can be obtained for the linearization $X _Q$ at the regular relative equilibria of the generalized orbitron:

\begin{equation*}
{\small 
\!\!\!\!\!
\left( \begin {array}{cccccccc} 
0&
-\xi_1 ^0 r&
0&
0&
0&
0&
0&
0\\
-\dfrac{4 \mu }{Mr \xi _1 ^0 }\left( 2 f' _1 + r ^2 f '' _1 \right) &
0&
0&
0&
0&
0&
0&
0\\
0&
0&
0&
\xi_1^0-\dfrac{\Pi_0}{I_1}&
0&
0&
0&
- \mu f _0 \\
0&
0&
-\xi_1^0+\dfrac{\Pi_0}{I_1}&
0&
0&
2 \mu r  f _1' &
- \mu f _0 &
0\\
0&
0&
-\dfrac{M \xi_1^0 r}{I_1}&
0&
0&
\mu f '' _2 &
4 \mu r  f' _1  &0\\
0&0&0&0&\dfrac{1}{M}&0&0&\xi_1^0 r\\
0&0&0&\dfrac{1}{I_1}&0&0&0&-\xi_1^0\\
0&0&\dfrac{1}{I_1}&0&0&0&\xi_1^0&0
\end {array} \right).}
\end{equation*}
\medskip

\noindent {\bf The linearization for the singular branches.} The same scheme can be reproduced for the singular branches by using the stability space introduced in~(\ref{stability space singular}). In this case, it can be shown that $W ^ \omega ={\rm span} \left\{ \boldsymbol{\tau} _1 , \boldsymbol{\tau} _2 \right\} $, with 
$\boldsymbol{\tau} _1 = (\mathbf{e}_3 , \mathbf{0}, \mathbf{0} , \mathbf{0}  )$ and 
$\boldsymbol{\tau} _2 = (\mathbf{0} , \mathbf{0} , \mathbf{e}_3  , \mathbf{0})$, which yields the following matrix expressions for the inclusion and the projection:
\begin{eqnarray*}
\boldsymbol{i}_W &=& \left(\mathbf{u}_1'\vert  \mathbf{u}_2'\vert \mathbf{u}_3'\vert \mathbf{u}_4'\vert \mathbf{u}_5'\vert \mathbf{u}_6'\vert \mathbf{u}_7'\vert \mathbf{u}_8'\vert \mathbf{u}_9'\vert \mathbf{u}_{10}'\right),\\ 
\mathbb{P}_W &=& \left(\mathbf{e}_{10}'\vert  \mathbf{e}_9'\vert \mathbf{0}\vert \mathbf{e}_8'\vert \mathbf{e}_7'\vert \mathbf{e}_2'\vert \mathbf{e}_6'\vert \mathbf{e}_5' \vert \mathbf{0}' \vert \mathbf{e}_4 ' \vert \mathbf{e}_3' \vert \mathbf{e}_1' \right).
\end{eqnarray*}
Finally, the  matrix expression for $X _Q $ at the singular relative equilibria of the standard orbitron is: 

\begin{equation}
\label{xq singular}
{\small
\left( \begin {array}{cccccccccc} 
0&-3\dfrac{\mu _0 \mu q}{\pi h ^4 }&0&0&0&0&0&0&0&0\\ 
\dfrac{1}{M}&0&0&0&0&0&0&0&0&0\\ 
0&0&0&- \xi _1 &0&0&\dfrac{3}{2}\dfrac{\mu _0 \mu q}{\pi h ^4 }&0&0&0\\ 
0&0&\xi _1 &0&0&0&0&\dfrac{3}{2}\dfrac{\mu _0 \mu q}{\pi h ^4 }&0&0\\ 
0&0&0&0&0&- \xi _1+ \dfrac{\Pi_0}{I _1} &0&0&\dfrac{1}{2}\dfrac{\mu _0 \mu q}{\pi h ^2 }&0\\ 
0&0&0&0&\xi _1-\dfrac{\Pi_0}{I _1} &0&0&0&0&\dfrac{1}{2}\dfrac{\mu _0 \mu q}{\pi h ^2 }\\ 
0&0&\dfrac{1}{M}&0&0&0&0&- \xi _1 &0&0\\ 
0&0&0&\dfrac{1}{M}&0&0&\xi _1 &0&0&0\\ 
0&0&0&0&\dfrac{1}{I _1 }   &0&0&0&0&-\xi _1  \\ 
0&0&0&0&0&\dfrac{1}{I _1 }  &0&0&\xi _1  &0\\ 
\end {array} \right),
}
\end{equation}
\medskip
where
$\Pi_0=I_3 (\xi_1-\xi_2)$. An analog expression can be obtained for the linearization $X _Q$ at the singular relative equilibria of the generalized orbitron:

\begin{equation}
\label{xq singular generalized}
{\small
\left( \begin {array}{cccccccccc} 
0&
\mu f'' _2 &
0&
0&
0&
0&
0&
0&
0&
0\\ 
\dfrac{1}{M}&
0&
0&
0&
0&
0&
0&
0&
0&
0\\ 
0&
0&
0&
- \xi _1 &
0&
0&
2 \mu f' _1 &
0&
0&
0\\ 
0&
0&
\xi _1 &
0&
0&
0&
0&
2 \mu f' _1 &
0&
0\\ 
0&
0&
0&
0&
0&
- \xi _1+ \dfrac{\Pi_0}{I _1} &
0&
0&
- \mu f _0 &
0\\ 
0&
0&
0&
0&
\xi _1-\dfrac{\Pi_0}{I _1} &
0&
0&
0&
0&
- \mu f _0 \\ 
0&0&\dfrac{1}{M}&0&0&0&0&- \xi _1 &0&0\\ 
0&0&0&\dfrac{1}{M}&0&0&\xi _1 &0&0&0\\ 
0&0&0&0&\dfrac{1}{I _1 }   &0&0&0&0&-\xi _1  \\ 
0&0&0&0&0&\dfrac{1}{I _1 }  &0&0&\xi _1  &0\\ 
\end {array} \right).
}
\end{equation}
\medskip

\addcontentsline{toc}{section}{Bibliography}
\bibliographystyle{alpha}
\bibliography{/Users/JP17/JPO_synch/BiblioData/GOLibrary} 
\end{document}